\documentclass{amsart}
\usepackage{fullpage,amssymb}
\usepackage{epsfig,pictex}
\usepackage{url}
\usepackage{color}

\newcommand{\dbar}{{\overline{\partial}}}
\urldef{\mcl}\url{kenmc@email.unc.edu}
\urldef{\millerpd}\url{millerpd@umich.edu}
\urldef{\OPUC}\url{http://www.math.lsa.umich.edu/~millerpd/ResearchAndPublications/CP.html}
\urldef{\DOP}\url{math.CA/0310278}

\newtheorem{lemma}{Lemma}[section]
\newtheorem{assumption}{Assumption}[section]
\newtheorem{prop}{Proposition}[section]
\newtheorem{theorem}{Theorem}[section]
\newtheorem{corollary}{Corollary}[section]
\newtheorem{rhp}{Riemann-Hilbert Problem}[section]
\newtheorem{dbp}{$\dbar$ Problem}[section]

\newenvironment{remark}{$\triangleleft$\,\,{\bf Remark:}}{$\triangleright$}

\title[The $\dbar$ steepest descent method]
{The $\dbar$ steepest descent method and the asymptotic behavior of
polynomials orthogonal on the unit circle with fixed and exponentially
varying nonanalytic weights}
\author{K. T.-R. McLaughlin}
\address{K. T.-R. McLaughlin:  Department of Applied Mathematics\\ CB \#3250 Phillips Hall\\ University of North Carolina at Chapel Hill\\ Chapel Hill, NC 27599-3250\\
Email address: \mcl }
\author{P. D. Miller}
\address{P. D. Miller:  Department of Mathematics\\ University of Michigan\\ East Hall\\ 525 E. University Avenue\\ Ann Arbor, MI 48109-1109\\  Email address: \millerpd}
\date{\today}
\begin{document}
\begin{abstract}
  We develop a new asymptotic method for the analysis of matrix
  Riemann-Hilbert problems.  Our method is a generalization of the
  steepest descent method first proposed by Deift and Zhou; however
  our method systematically handles jump matrices that need not be
  analytic.  The essential technique is to introduce nonanalytic
  extensions of certain functions appearing in the jump matrix, and to
  therefore convert the Riemann-Hilbert problem into a $\dbar$
  problem.  We use our method to study several asymptotic problems of
  polynomials orthogonal with respect to a measure given on the unit
  circle, obtaining new detailed uniform convergence results, and for
  some classes of nonanalytic weights, complete information about the
  asymptotic behavior of the individual zeros.
\end{abstract}
\maketitle
\section{Introduction}
\subsection{Asymptotic analysis of Riemann-Hilbert problems.}
The steepest descent method for asymptotic analysis of matrix
Riemann-Hilbert problems was introduced by Deift and Zhou in 1993
\cite{steepestintro}.  A matrix Riemann-Hilbert problem is specified
by giving a triple $(\Sigma,{\bf v},\mathcal{N})$ consisting of an
oriented contour $\Sigma$ in the complex $z$-plane, a matrix function
${\bf v}:\Sigma \rightarrow SL(N)$ which is usually taken to be
continuous except at self-intersection points of $\Sigma$ where a
certain compatibility condition is required, and a normalization
condition $\mathcal{N}$ as $z\rightarrow\infty$.  If $\Sigma$ is not
bounded, certain asymptotic conditions are required of ${\bf v}$ in
order to have compatibility with the normalization condition.
Consider an analytic function ${\bf
  M}:\mathbb{C}\setminus{\Sigma}\rightarrow SL(N)$ taking continuous
boundary values ${\bf M}_+(z)$ (respectively ${\bf M}_-(z)$) on
$\Sigma$ from the left (respectively right).  The Riemann-Hilbert
problem $(\Sigma,{\bf v},\mathcal{N})$ is then to find such a matrix
${\bf M}(z)$ satisfying the normalization condition $\mathcal{N}$ as
$z\rightarrow\infty$ and the jump condition ${\bf M}_+(z)={\bf
  M}_-(z){\bf v}(z)$ whenever $z$ is a non-self-intersection point of
$\Sigma$ (so the left and right boundary values are indeed
well-defined).  The steepest descent method of Deift and Zhou applies
to certain Riemann-Hilbert problems where the jump matrix ${\bf v}(z)$
depends on an auxiliary control parameter, and is a method for
extracting asymptotic properties of the solution ${\bf M}(z)$ (and
indeed proving the existence and uniqueness of solutions along the
way) when the control parameter tends to a singular limit of interest.

The original method put forth in \cite{steepestintro} bears a
striking resemblance to the well-known steepest descent method or
saddle point method for analyzing contour integrals with
exponential integrands.  A distinguished point on $\Sigma$ is
identified (analogous to a point of stationary phase) and an
explicit change of variables of the form ${\bf N}(z)={\bf
M}(z){\bf t}(z)$ where ${\bf t}(z)$ is a piecewise analytic matrix
is introduced in the vicinity of this point and it is observed
that (i) the matrix ${\bf N}(z)$ satisfies an equivalent
Riemann-Hilbert problem with a new contour $\Sigma_{\bf N}$ and a
new jump matrix ${\bf v}_{\bf N}$, and (ii) the jump matrix ${\bf
v}_{\bf N}$ converges to the identity matrix in the singular limit
of interest for all $z$ bounded away from the stationary phase
point.  One therefore expects that a good approximation to ${\bf
N}(z)$ can be constructed by an explicit local analysis near the
stationary phase point.  With the explicit local approximant
$\dot{\bf N}(z)$ constructed, one uses it in a final change of
variables ${\bf H}(z)={\bf N}(z)\dot{\bf N}(z)^{-1}$ and observes
(i) that ${\bf H}(z)$ satisfies a Riemann-Hilbert problem with a
possibly new contour $\Sigma_{\bf H}$ and a new explicit jump
matrix ${\bf v}_{\bf H}$ and (ii) that the new jump matrix ${\bf
v}_{\bf H}$ is now {\em uniformly} close to the identity matrix in
the limit of interest.  This allows one to construct ${\bf H}(z)$
by iteration of certain singular integral equations that are
equivalent to any given Riemann-Hilbert problem, and to show that
${\bf H}(z)$ is uniformly close to the identity matrix in any
region bounded away from $\Sigma_{\bf H}$. With additional work it
may in some circumstances be shown that ${\bf H}(z)$ is close to
the identity uniformly right up to the contour $\Sigma_{\bf H}$
(this usually requires more detailed information about the jump
matrix).  In this way, one obtains a formula ${\bf M}(z)= {\bf
H}(z)\dot{\bf N}(z){\bf t}(z)^{-1}$ for the solution that can be
used to compute {\em directly} an asymptotic expansion of ${\bf
M}(z)$ valid in the singular limit of interest.

Since the introduction of the steepest descent method for
Riemann-Hilbert problems, there have been several key developments. In
\cite{gfunction} and \cite{steepestextend} a technique was established
in which one makes a change of variables involving a matrix
constructed from a single unknown scalar function $g(z)$ analytic in
${\mathbb C}\setminus\Sigma$.  The transformation modifies the jump
matrix in a way involving the boundary values $g_\pm(z)$ taken on
$\Sigma$.  One then chooses relations between the boundary values of
$g(z)$ such that the transformed Riemann-Hilbert problem becomes
asymptotically simple.  The desired conditions amount to a {\em
scalar} Riemann-Hilbert problem for $g(z)$, which is easily solved in
many circumstances.  A crucial feature of this method is that the
dominant contribution to the solution typically comes from
subintervals of the contour $\Sigma$ of finite length rather than from
isolated points.  Here we therefore see an important difference
between singular limits of matrix Riemann-Hilbert problems and
evaluation of saddle-point integrals.  In the contributing intervals
the transformed jump matrix has a factorization (see
\eqref{eq:targetform}, \eqref{eq:factorizationstrong}, and 
\eqref{eq:factorizationvarying} below) whose factors admit analytic continuation to the left and right
of each such interval. A further change of variables based on this
analytic factorization is carried out in lens-shaped regions
surrounding each contributing interval. Ultimately, a model problem is
solved (typically in terms of Riemann theta functions of genus related
to the number of contributing intervals) and along with local analysis
near the endpoints of the intervals a model for ${\bf M}(z)$ is built
and compared with ${\bf M}(z)$ to obtain a Riemann-Hilbert problem for
the error ${\bf H}(z)$. When the method is successful, the jump matrix
for ${\bf H}(z)$ is uniformly close to the identity and thus ${\bf
H}(z)$ may be constructed via iteration of integral equations.
Significantly, the conditions imposed on the boundary values of $g(z)$
can often be viewed as the Euler-Lagrange conditions for a certain
variational problem (see \cite{steepestKdV} as well as
\cite{op1,op2} and Appendix~\ref{app:potential} of this paper).

A further development emerged from problems in which it was recognized
that no appropriate function $g(z)$ can be found relative to the given
contour $\Sigma$.  In \cite{manifesto} and later in \cite{tail} it was
shown how analyticity of the jump matrix could be exploited to
effectively deform arcs of the contour $\Sigma$ to alternative
locations in the complex plane such that the jump matrix maintains the
same functional form; specific locations of the arcs are determined
such that there exists an appropriate function $g(z)$ as above.  These
selected arcs are the closest relatives in the noncommutative theory
to the paths of steepest descent from saddle points in the asymptotic
theory of contour integrals.  The contour selection principle was also
encoded into a variational problem in \cite{manifesto}.

More recently \cite{manifesto,WKB,dop}, new techniques have been added to the framework of the
steepest descent method that are adapted for determining the asymptotic contribution to
the solution of a coalescence of a large number of poles in the matrix
unknown (this is strictly speaking not a Riemann-Hilbert problem in
the sense described above due to the polar singularities, however the
problem is first converted into a standard Riemann-Hilbert problem by
explicit transformations).  The key idea here is to
exploit certain analytic interpolants of given residues at the poles.

For the fundamentally nonlinear cases in which the dominant
contribution comes from subintervals of a contour, a central
feature is that the analytical methods rely on piecewise
analyticity of the given jump matrix and of the boundary values of
the scalar function $g(z)$.  For some of the cases of long-time
asymptotics of integrable nonlinear partial differential equations
\cite{steepestintro,DIZ}, as well as the recent long time
asymptotic analysis for perturbations of the defocusing nonlinear
Schr\"{o}dinger equation \cite{DZ02}, the dominant contribution
comes from isolated points of the contour $\Sigma$, and while
analyticity is not fundamental, the asymptotic calculations
proceed by an  approximation argument, in which an analytic part
is deformed away, and a (small) residual contribution is handled
by technical and analytical prowess. The approximation argument is
delicate and requires detailed analysis that depends sensitively
on the geometry of the particular contour $\Sigma$.

Far from being a mere pursuit of abstraction, a simple asymptotic
technique that applies to Riemann-Hilbert problems regardless of
whether the jump matrix is analytic or not would have immediate
application in a number of important areas.  For example, a
unified treatment of the asymptotic theory of orthogonal
polynomials on the real line with general nonanalytic weights
would allow the resolution of universality conjectures from random
matrix theory in the most natural and general context (see
\cite{op1,op2} for the analytic case, and \cite{otherpaper}
for an application of the $\dbar$ steepest descent method described in
this paper to the nonanalytic but convex case).  As another example,
if it were possible to treat systematically problems with a large
number of poles that accumulate in a very regular but nonanalytic
fashion, then the important problem of semiclassical asymptotics for
the focusing nonlinear Schr\"odinger equation with general nonanalytic
initial data could begin to be addressed.

In this paper, we present a new generalization of the steepest
descent method for Riemann-Hilbert problems that applies in
absence of analyticity of the jump matrix, and yet does not depend
on an approximation argument for the jump matrix.  While we
believe that the ideas we will develop in this paper are useful in
very general contexts, we have chosen to focus on a particular
application of the steepest descent method in order to demonstrate the
technique.  

\subsection{The essence of the $\dbar$ steepest descent method.}
\begin{figure}[h]
\begin{center}
\input{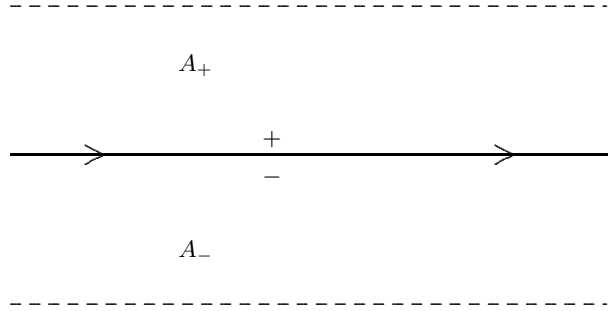}
\end{center}
\caption{\em The complex plane in the vicinity of a contour $\Sigma$
(here coincident with the real axis) supporting a factorized
jump matrix.}
\label{fig:intro}
\end{figure}

As mentioned above, after changing variables using an appropriate
scalar function $g(z)$, the jump matrix is converted into a form
that is well-suited for further asymptotic analysis. A common
``target'' form for the jump matrix in certain arcs of $\Sigma$ is
the following form
\begin{equation}
{\bf v}(x) = \left(\begin{array}{cc} e^{i\kappa(x)/\epsilon} & 1 \\\\
0 & e^{-i\kappa(x)/\epsilon}\end{array}\right) =
\left(\begin{array}{cc} 1 & 0 \\\\ e^{-i\kappa(x)/\epsilon} & 1
\end{array}\right)
\left(\begin{array}{cc} 0 & 1 \\\\-1 & 0\end{array}\right)
\left(\begin{array}{cc} 1 & 0 \\\\e^{i\kappa(x)/\epsilon} & 1
\end{array}\right)
\label{eq:targetform}
\end{equation}
where $x$ is a real parameter along the arc of $\Sigma$ which for
simplicity here we assume lies on the real axis (for a
representation see Figure~\ref{fig:intro}), $\epsilon>0$ is 
the control parameter tending to zero
in the singular limit of interest, and $\kappa(x)$ is a strictly
increasing real function of $x$ that is related to the boundary
values of $g(z)$ on $\Sigma$.
Suppose ${\bf N}(z)$ is the unknown
satisfying ${\bf N}_+(x)={\bf N}_-(x){\bf v}(x)$.  With the
assumption of analyticity of $\kappa(x)$, one may transform the
Riemann-Hilbert problem by introducing as a new unknown a matrix
${\bf O}(z)$ defined in terms of ${\bf N}(z)$ by the following
scheme:  in some region lying on the minus side of the arc
(the region labeled $A_-$ in Figure~\ref{fig:intro}), set
\begin{equation}
{\bf O}(x+iy):={\bf N}(x+iy)\left(\begin{array}{cc} 1 & 0 \\\\ e^{-i\kappa(
x+iy)/\epsilon} & 1\end{array}\right)
\end{equation}
and in some region lying on the plus side of the arc (the region labeled
$A_+$ in Figure~\ref{fig:intro}), set
\begin{equation}
{\bf O}(x+iy):={\bf N}(x+iy)\left(\begin{array}{cc} 1 & 0 \\\\
-e^{i\kappa(x+iy)/\epsilon} & 1\end{array}\right)\,.
\end{equation}
Elsewhere, set ${\bf O}(z)={\bf N}(z)$.  On has thus introduced
two new jump contours, one on either side of the arc (these
are the two dashed lines in
Figure~\ref{fig:intro}). However, the monotonicity of the real analytic
function $\kappa(x)$ implies via the Cauchy-Riemann equations that
the induced jump matrix relating the boundary values of ${\bf
O}(z)$ on these two contours is exponentially close to the
identity matrix in the limit $\epsilon\downarrow 0$.  On the
original arc, the matrix ${\bf O}(z)$ satisfies the constant jump
relation
\begin{equation}
{\bf O}_+(x)={\bf O}_-(x)\left(\begin{array}{cc} 0 & 1 \\\\ -1 & 0
\end{array}\right)\,,
\end{equation}
which can be subsequently dealt with in terms of special functions.
In this paper, we show how this procedure can be carried out
effectively when one does not have the option of extending $\kappa(x)$
from the contour $\Sigma$ because it is not assumed to be an analytic
function.  We choose to extend $\kappa(x)$ in a way that does not
assume any analyticity (see \eqref{eq:extensiondefine} below).  The
price that must be paid is that the analogue of the matrix ${\bf
O}(z)$ above is no longer analytic in the regions to the left and right of
the arc; therefore this matrix cannot be the solution of any
Riemann-Hilbert problem.  It can, however, be the solution of a matrix
$\dbar$ problem (or more generally a mixed Riemann-Hilbert-$\dbar$
problem).  It is into this framework that we extend the steepest
descent method.  This explains the terminology of the ``$\dbar$
steepest descent method''.

As the fundamental contour $\Sigma$ in this paper is the unit circle $S^1$,
we can now be very specific about what we mean by an extension of a nonanalytic function in this context.  Suppose that $f(\theta)$ is a $C^{m-1}(S^1)$ function,
$m=1,2,3,\dots$.  Then, we define an extension operator $E_m:C^{m-1}(S^1)\rightarrow C(\mathbb{R}^2\setminus\{0\})$ as follows:
\begin{equation}
E_mf(r,\theta):=\sum_{p=0}^{m-1}\frac{f^{(p)}(\theta)}{p!}(-i\log(r))^p\,,
\label{eq:extensiondefine}
\end{equation}
where $(r,\theta)$ are the standard polar coordinates for $\mathbb{R}^2$.
Note that this indeed defines a continuous extension to any annulus
$r_+\le r\le r_-$ where $0<r_+<1< r_-<\infty$ since
$E_mf(1,\theta)=f(\theta)$.  Also, since 
$z=re^{i\theta}$ and $\overline{z}=re^{-i\theta}$, 
the fundamental differential operators of complex variable theory
are represented in polar coordinates as 
\begin{equation}
\dbar := \frac{\partial}{\partial\overline{z}} = \frac{e^{i\theta}}{2}\left(
\frac{\partial}{\partial r}+\frac{i}{r}\frac{\partial}{\partial\theta}\right) \hspace{0.2 in}\text{and}\hspace{0.2 in}
\partial:=\frac{\partial}{\partial z} = \frac{e^{-i\theta}}{2}\left(\frac{\partial}{\partial r}-\frac{i}{r}\frac{\partial}{\partial\theta}\right)
\end{equation}
and therefore we see that if $f^{(m-1)}(\theta)$ is Lipschitz, then in particular 
it has a derivative almost everywhere that is uniformly bounded, and we have
\begin{equation}
\dbar E_mf(r,\theta) = \frac{ie^{i\theta}}{2r}\frac{f^{(m)}(\theta)}{(m-1)!}
(-i\log(r))^{m-1}
\label{eq:dbarTaylor}
\end{equation}
and
\begin{equation}
\partial E_mf(r,\theta) = -\frac{ie^{-i\theta}}{2r}
\left(\frac{f^{(m)}(\theta)}{(m-1)!}(-i\log(r))^{m-1} + 2\sum_{p=0}^{m-2}
\frac{f^{(p+1)}(\theta)}{p!}(-i\log(r))^p\right)
\label{eq:dTaylor}
\end{equation}
both holding for all $r\ge 0$ and almost all $\theta\in S^1$ (these
formulae hold at every point of the plane if $f$ is of class
$C^{m}(S^1)$).  It follows that $\dbar E_mf(r,\theta)$ vanishes to
order $m-1$ as $r\rightarrow 1$ uniformly in $\theta$.  In fact, if
$f(\theta)$ is analytic for all $\theta$, then the infinite series
$E_\infty f(r,\theta)$ converges uniformly in some annulus containing
the unit circle $r=1$ and represents the unique analytic extension of
$f(\theta)$.

Generally speaking, Riemann-Hilbert problems with rapidly oscillatory
jump matrices are equivalent to systems of singular integral equations
with Cauchy kernel and rapidly oscillatory densities.  Such equations
can in principle be analyzed asymptotically \cite{varzugin}.  This
approach requires delicate arguments of harmonic analysis.  On the
other hand, the $\dbar$ steepest descent method we will develop in
this paper avoids such complicated reasoning.  Indeed, by extending
contour integration into integration over two-dimensional regions, the
Cauchy kernel becomes less singular, and the analysis becomes
correspondingly more straightforward.

In the analytic case described briefly above, the asymptotic analysis
is in general complicated by the fact that the procedure is valid in
the neighborhood of certain intervals of $\Sigma$, and it turns out
that a different analysis must be carried out in the vicinity of the
endpoints of the intervals.  The same would be expected to be true in
the general nonanalytic case.  In order to have the clearest possible
presentation, we have chosen to describe in this paper the $\dbar$
steepest descent method in the context of a problem where there are
nontrivial cases without endpoint issues, namely the asymptotic
behavior of polynomials orthogonal with respect to weights on the unit
circle.  While the presence of endpoints complicates the analysis,
they do not present an insurmountable obstruction, and the reader is
referred to \cite{otherpaper} for a description of the more general
theory.

\subsection{Polynomials orthogonal on the unit circle.}
Let $\phi(\theta)$ be an integrable $2\pi$-periodic
function satisfying $\phi(\theta)>0$ for almost all $\theta$.  For two
complex-valued functions $f(\theta)$ and $g(\theta)$ there is an
associated inner product
\begin{equation}
\langle f,g\rangle_\phi:=\frac{1}{2\pi}\int_{-\pi}^{\pi} 
f(\theta)\overline{g(\theta)}\phi(\theta)\,d\theta = 
\frac{1}{2\pi i}\oint_\Sigma f(\arg(z))\overline{g(\arg(z))}
\phi(\arg(z))\frac{dz}{z}\,.
\label{eq:innerproduct}
\end{equation}
This inner product leads to a system $\{p_n(z)\}_{n=0}^\infty$ 
of orthogonal polynomials in the
complex variable $z$:  
\begin{equation}
p_n(z)=\gamma_nz^n + \sum_{j=0}^{n-1}c_{n,j}z^j\,,\hspace{0.2 in}
\text{where} \hspace{0.2 in}
\gamma_n>0\,,
\end{equation}
and the defining relation is the orthonormality condition
\begin{equation}
\langle p_m,p_n\rangle_\phi = \delta_{mn}\,,\hspace{0.2 in}
0\le m,n <\infty\,.
\end{equation}
Here the polynomials $p_n(z)$ are considered as complex-valued
functions of $\theta$ by restriction to the unit circle:
$z=e^{i\theta}$.  The constants $\gamma_n$ have the interpretation of
normalization constants, and the corresponding system
$\{\pi_n(z)\}_{n=0}^\infty$ of monic orthogonal polynomials is defined
by rescaling:
\begin{equation}
\pi_n(z)=\frac{1}{\gamma_n}p_n(z)\,,\hspace{0.2 in}
0\le n<\infty\,.
\end{equation}
The orthogonal polynomials satisfy recurrence relations of the form
\begin{equation}
\begin{array}{rcl}
\displaystyle \pi_{n+1}(z)&=&\displaystyle z\pi_n(z) + 
\alpha_{n+1}z^n\overline{\pi_n(1/\overline{z})}\\\\
\displaystyle z^{n+1}\overline{\pi_{n+1}(1/\overline{z})} &=&
\displaystyle \overline{\alpha_{n+1}}z\pi_n(z) + z^n\overline{\pi_n(1/\overline{z})}\,,
\end{array}
\label{eq:recurrence}
\end{equation}
for $n=0,1,2,3,\dots$.  Here, the complex constants
$\{\alpha_n\}_{n=1}^\infty$ are the {\em recurrence coefficients}
associated with the weight $\phi$ (also known
as the {\em Schur parameters} or {\em Verblunsky coefficients}).
By setting $z=0$ in the first equation of
\eqref{eq:recurrence}, it is easy to see that
\begin{equation}
\alpha_n = \pi_n(0)\,.
\label{eq:Verblunskypin0}
\end{equation}

For a general discussion of properties of polynomials orthogonal on
the unit circle, see Chapter XII of Szeg\H{o}'s monograph
\cite{szego}, in which (among other things) the asymptotic behavior of
$\pi_{n}(z)$ for $n \to \infty$ is discussed.  The extraction of
asymptotic formulae for quantities related to the orthogonal
polynomials of large degree with a fixed weight $\phi$ on the unit
circle is the type of asymptotic problem in the theory of general
orthogonal polynomials (that is, beyond particular cases involving
classical special functions) for which results have been known for the
longest time.  This can be traced to the fact that if
$\phi(\theta)^{-1}$ happens to be a positive trigonometric polynomial,
then there is a closed-form expression for the orthonormal polynomial
$p_n(z)$ that is convenient for analysis, as long as $n$ is
sufficiently large compared to the degree of $\phi(\theta)^{-1}$.  See
\cite[\S 11.2]{szego}.  In other words, for certain special fixed
weights $\phi(\theta)$ the asymptotic formulae one obtains become {\em
  exact} as long as $n$ is large enough.  This leads to a general
strategy for asymptotic analysis of orthogonal polynomials on the unit
circle based on approximating an arbitrary given positive function
$\phi(\theta)^{-1}$ by positive trigonometric polynomials.

The asymptotics described in the monograph of Szeg\H{o} are of a
rather general character and hold whenever $\log(\phi(\theta))$ is an
integrable real valued function.  In the years since the origin of
Szeg\H{o}'s methods, there have been many further developments in the
asymptotic theory.  These developments move both in the direction of
generalizing the class of weights for which the Szeg\H{o} asymptotics
are valid (perhaps in a weaker form) and also in the direction of
trading generality of the weight for detail of the asymptotics.  It
seems that certain problems remain difficult to treat by these
methods; in particular, it is difficult to verify convergence in a
uniform sense, and it is difficult to characterize the detailed
asymptotic behavior of zeros.  There is a vast literature on this
subject; we refer the interested reader to the memoir of Nevai
\cite{nevai} and the forthcoming monograph of Simon \cite{simon}.

Beyond being a source of classical information about orthogonal
polynomials on the unit circle, Simon's monograph describes a
different viewpoint of the theory of these polynomials.  Namely, Simon
and his school have made great progress by exploiting the connection
between orthogonal polynomials and spectral theory for operators that
encode the recurrence relations that all orthogonal polynomials
satisfy (see also \cite{GeronimoC79}).  This theory is capable of
establishing a number of very general results relevant to asymptotics
in the limit of large degree.  An important point is that the
hypotheses required to establish results of this kind involve
assumptions about the asymptotic behavior of the sequence
$\{\alpha_n\}_{n=1}^\infty$ of recurrence coefficients.  Indeed, the
fundamental problem of spectral theory in this context is the
construction of the spectral measure $\phi(\theta)\,d\theta$ from the
finite difference operator involving the recurrence coefficients
$\{\alpha_n\}_{n=1}^\infty$.

On the other hand, the recovery of the polynomials
$\{p_n(z)\}_{n=0}^\infty$ and of the recurrence coefficients
$\{\alpha_n\}_{n=1}^\infty$ from the spectral measure
$\phi(\theta)\,d\theta$ is the fundamental problem of inverse spectral
theory.  A general approach to asymptotic problems in the theory of
orthogonal polynomials in which the measure of orthogonality is the
given data therefore involves the translation of the orthogonality
conditions into the conditions making up a Riemann-Hilbert problem for
sectionally analytic matrices.  A Riemann-Hilbert formulation for
polynomials orthogonal with respect to a measure on the unit circle
was described in \cite{longest} and follows closely the well-known
Riemann-Hilbert formulation for polynomials orthogonal with respect to
a measure on $\mathbb{R}$ discovered in \cite{FIK2}.  In
\cite{longest}, and in a number of papers which followed (see, for
example, \cite{longest2,BR1,BR2}), the polynomial of degree $n$
orthogonal with respect to a specific family of weights of the form
$\phi(\theta) = e^{-nV(\theta)}$ where $V(\theta)=\gamma\cos(\theta)$
was studied in the limit $n\rightarrow\infty$.  Note that this is a
joint limit as the degree $n$ of the polynomial in question appears in
the measure of orthogonality as well; see \S~\ref{sec:varying} for a
general discussion of such exponentially varying weights.  For the
large $n$ asymptotics carried out in \cite{longest}, and in subsequent
works with this measure, as well as closely related measures (see, for
example, \cite{tail}), analyticity of the weight $\phi(\theta)$ played
a central role in the analysis.

In \cite{Deift}, Deift used polynomials orthogonal with respect to a
measure on the unit circle to give an example of his theory of
integrable operators.  Specifically, he introduced a one-parameter
family of positive, analytic functions $\phi(\theta;t)$ and related
solutions of Riemann-Hilbert Problem~\ref{rhp:M} in \S~\ref{sec:ops} below
(with $\phi(\theta)$ replaced by $\phi(\theta;t)$) to Toeplitz
determinants.  The $n \to \infty$ asymptotic behavior of the
corresponding solution ${\bf M}^n(z;t)$ to Riemann-Hilbert
Problem~\ref{rhp:M} then yields asymptotics for the associated
Toeplitz determinants. Exploiting the analyticity of $\phi(\theta;t)$,
Deift outlined how one obtains an asymptotic description for ${\bf
  M}^n(z;t)$.  The calculations in we will present in
\S~\ref{sec:strong} may be viewed as complementary to this asymptotic
calculation of \cite{Deift}, in that we will establish asymptotics for
orthogonal polynomials under the much weaker assumption that
$\phi(\theta)$ is a continuous function satisfying a Lipschitz
condition. Furthermore, we show how the error estimates depend on
smoothness properties of $\phi(\theta)$.

\subsection{Outline and summary of results.}
The polynomials orthogonal with respect to a weight given on the
unit circle in the complex plane can be characterized in terms of
the solution of a matrix Riemann-Hilbert problem in which the
contour $\Sigma$ is the unit circle.  In \S~\ref{sec:ops} we
describe this Riemann-Hilbert problem, and then in
\S~\ref{sec:strong} and \S~\ref{sec:varying} we study the singular
limit in which the degree of the polynomials tends to infinity. In
\S~\ref{sec:strong} we consider the weight function to be
held fixed as the degree tends to infinity, while in
\S~\ref{sec:varying} we study the joint limit when the degree
becomes large while the weight function is exponentially varied.
A summary of the relevant logarithmic potential theory
referred to in \S~\ref{sec:varying} is given  
in Appendices~\ref{app:potential} and \ref{app:compare}.

The key results we obtain in the fixed weights case are described in
\S~\ref{sec:strongasymptoticspi}.  To the best of our knowledge, the
uniform nature of the asymptotics we obtain is new to the field, as is
our detailed characterization of the zeros.  While there exist several
classical methods available for the asymptotic analysis of orthogonal
polynomials on the unit circle with fixed weight function
$\phi(\theta)$, with the degree $n$ of the polynomial in question
tending to infinity, this problem is the ideal context in which to
introduce the $\dbar$ steepest descent method.

On the other hand, the asymptotic behavior of polynomials orthogonal
with respect to a varying weight on the unit circle, considered in
\S~\ref{sec:varying}, is more challenging to obtain by more classical
techniques.  The results we obtain with the use of the $\dbar$ steepest
descent method are stated in \S~\ref{sec:asymptoticspids}.
A point we wish to emphasize is that with the use of the $\dbar$
steepest descent method, the analysis in the varying weights case is
no more difficult than in the fixed weights case.  This fact
distinguishes the $\dbar$ steepest descent method from more classical
techniques.

\subsection{Notation.}  
\label{sec:notation}
Throughout, we assume a fixed norm $\|\cdot\|$ on
$2\times 2$ matrices.  For $p=0,1,2,\dots$, we use the following
induced norm on sufficiently smooth matrix functions ${\bf
F}:\mathbb{R}^2\rightarrow\mathbb{C}_{2\times 2}$
\begin{equation}
\||{\bf F}|\|_p:=\sum_{\alpha+\beta\le p}\sup_{\mathbb{R}^2}
\left\|\frac{\partial^{\alpha+\beta}{\bf F}}{\partial x^\alpha\partial y^\beta}
\right\|\,.
\label{eq:pnorm}
\end{equation}
Here $x$ and $y$ are cartesian coordinates in $\mathbb{R}^2$.  For
$p=0$, we may apply this norm on all functions ${\bf F}$ in the space
$L^\infty(\mathbb{R}^2)$, whereas for $p>0$ we may apply this norm to
a subset of functions ${\bf F}$ in the space $C^{p-1,1}(\mathbb{R}^2)$
of functions with Lipschitz continuous mixed partial derivatives of
all orders up to and including $p-1$.  The finiteness of the norm
indicates the uniformity of the Lipschitz condition.  Since for
functions ${\bf F}$ in the class $C^{p-1,1}(\mathbb{R}^2)$ the mixed
partial derivatives of order $p$ exist almost everywhere, the
condition $\||{\bf F}|\|_p<\infty$ can be equivalently expressed as
saying that ${\bf F}$ has all derivatives of total order at most $p$
in the space $L^\infty(\mathbb{R}^2)$.  We also use the notation
$C^{p-1,1}_0(\mathbb{R}^2\setminus\{0\})$ and
$L^\infty_0(\mathbb{R}^2\setminus\{0\})$ for spaces of functions
respectively in $C^{p-1,1}(\mathbb{R}^2)$ and $L^\infty(\mathbb{R}^2)$
that vanish identically for $|\log(x^2+y^2)|$ large enough (that is,
outside some annulus).

For functions $V(\theta)$ defined on the circle $S^1$ (that is, $V$ is
defined for $-\pi\le\theta<\pi$), we also say that $V$ is of class
$C^{k-1,1}(S^1)$ if the periodic extension of $V$ to
$\theta\in\mathbb{R}$ has $k-1$ Lipschitz continuous derivatives, or
equivalently, has $k$ derivatives in $L^\infty(\mathbb{R})$.  A suitable norm
for such functions is given by 
\begin{equation}
\||V|\|_{\circ,k}:=\sup_{-\pi<\theta<\pi}|V(\theta)| + \sup_{-\pi<\theta<\pi}
|V^{(k)}(\theta)|\,,
\end{equation}
since it is easy to establish that for all $m$ satisfying $1\le m \le k-1$, 
\begin{equation}
\sup_{-\pi<\theta<\pi}|V^{(m)}(\theta)|\le (2\pi)^{k-m}\sup_{-\pi<\theta<\pi}|V^{(k)}(\theta)|\,.
\label{eq:cutaway}
\end{equation}

If $V(\theta)$ satisfies a H\"older continuity condition, there exists
a unique function $N(z)$, analytic for $|z|>1$, decaying as
$z\rightarrow\infty$, and taking H\"older continuous boundary values
on $|z|=1$, such that
\begin{equation}
V(\theta)=N(e^{i\theta}) + V_0 + \overline{N(e^{i\theta})}\,,
\label{eq:logphidecompose}
\end{equation}
and $V_0$ is the average
value of $V$.  Thus, $N(e^{i\theta})$ is the negative
frequency component of the Fourier series for $V(\theta)$:
\begin{equation}
V(\theta)=\sum_{j=-\infty}^\infty V_je^{ij\theta}\,,\hspace{0.2 in}
\text{with coefficients}\hspace{0.2 in}
V_j=\frac{1}{2\pi}\int_{-\pi}^\pi V(\theta)e^{-ij\theta}\,d\theta\,,  
\label{eq:FourierSeriesV}
\end{equation}
and we have
\begin{equation}
N(z):=\sum_{j=1}^\infty\frac{V_{-j}}{z^j}\,,\hspace{0.2 in}\text{for
$|z|\ge 1$}
\,.
\label{eq:negfreq}
\end{equation}
We also introduce the function $\Omega:S^1\rightarrow\mathbb{R}$ by
the formula
\begin{equation} \Omega(\theta):=2\Im(N(e^{i\theta}))\,.
\label{eq:Omegadef}
\end{equation}
Note that $\Omega$ and $V$ are functions that are related by the
Cauchy transform.

Throughout the paper we will use a ``bump'' function
$B:\mathbb{R}\rightarrow[0,1]$ with the properties that $B$ is infinitely
differentiable, $B(l)\equiv 1$
for $|l|<\log(2)/2$, and $B(l)\equiv 0$ for $|l|>\log(2)$.

\section
{The Riemann-Hilbert Problem for Polynomials Orthogonal on the Unit
Circle}
\label{sec:ops}
\begin{figure}[h]
\begin{center}
\input{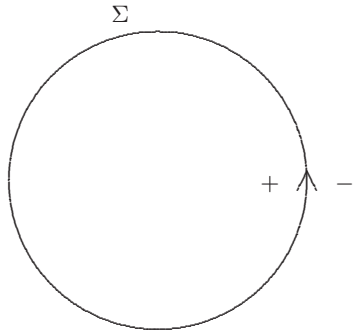}
\end{center}
\caption{\em The contour $\Sigma$ of the Riemann-Hilbert problem for
polynomials orthogonal on the unit circle is the unit circle itself $|z|=1$ oriented in the counterclockwise direction.}
\label{fig:circle}
\end{figure}

Consider the contour $\Sigma$ illustrated in Figure~\ref{fig:circle}.
Let $n$ be a positive integer.
Relative to the contour $\Sigma$ we pose, for each $n=0,1,2,3,\dots$, the
following Riemann-Hilbert problem for a $2\times 2$ matrix ${\bf
M}^n(z)$.
\begin{rhp}
Find a $2\times 2$ matrix ${\bf M}^n(z)$ with the properties:
\begin{itemize}
\item[]{\bf Analyticity.}
${\bf M}^n(z)$ is analytic for $|z|\neq 1$, and takes continuous
boundary values ${\bf M}^n_+(z)$, ${\bf M}^n_-(z)$ as $w$ tends to $z$
with $|z|=1$ and $|w|<1$, $|w|>1$.
\item[]{\bf Jump Condition.}
The boundary values are connected by the relation
\begin{equation}
{\bf M}^n_+(e^{i\theta})={\bf M}^n_-(e^{i\theta})
\left(\begin{array}{cc} 1 &
\phi(\theta)e^{-in\theta}\\\\ 0 & 1
\end{array}\right)\,.
\end{equation}
\item[]{\bf Normalization.}
The matrix ${\bf M}^n(z)$ is normalized at $z=\infty$ as follows:
\begin{equation}
\lim_{z\rightarrow\infty}{\bf M}^n(z)\left(\begin{array}{cc}
z^{-n} & 0 \\\\
 0 & z^n\end{array}\right)=\mathbb{I}\,.
 \label{eq:Mnorm}
\end{equation}
\end{itemize}
\label{rhp:M}
\end{rhp}
\begin{prop}
Suppose that the positive weight function $\phi(\theta)$ satisfies a
uniform H\"older condition $|\phi(\theta_2)-\phi(\theta_1)|\le
K|\theta_2-\theta_1|^\nu$ for some $\nu\in (0,1]$ and with some $K$
independent of $\theta_1$ and $\theta_2$.  Then Riemann-Hilbert
Problem~\ref{rhp:M} has a unique solution for each integer $n\ge 0$,
namely if $n>0$,
\begin{equation}
{\bf M}^n(z)=\left(\begin{array}{cc}
\pi_n(z) & \displaystyle\frac{1}{2\pi i}\oint_\Sigma
\frac{\pi_n(s)s^{-n}}{s-z}\phi(\arg(s))\,ds\\\\
\displaystyle -\gamma_{n-1}^2z^{n-1}\overline{\pi_{n-1}(1/\overline{z})} &
\displaystyle -\frac{\gamma_{n-1}^2}{2\pi i}\oint_\Sigma
\frac{\overline{\pi_{n-1}(s)}s^{-1}}{s-z}\phi(\arg(s))\,ds
\end{array}\right)\,,
\label{eq:Mnotzero}
\end{equation}
and if $n=0$,
\begin{equation}
{\bf M}^0(z)=\left(\begin{array}{cc} 1 & \displaystyle
\frac{1}{2\pi i}\oint_\Sigma\frac{1}{s-z}\phi(\arg(s))\,ds\\\\
0 & 1\end{array}\right)\,.
\label{eq:Mzero}
\end{equation}
In particular, $M^n_{11}(0)=\alpha_n$ and
$M_{21}^n(0)=-\gamma_{n-1}^2$ for $n\ge 1$.  Here,
$\{\pi_n\}_{n=0}^\infty$ is the sequence of monic orthogonal
polynomials with respect to the weight $\phi$ and the inner product
\eqref{eq:innerproduct}, and $\{\alpha_n\}_{n=1}^\infty$ is the
sequence of associated recurrence coefficients (see
\eqref{eq:recurrence}) while $\{\gamma_n\}_{n=0}^\infty$ is the
sequence of associated normalization constants.
\label{prop:M}
\end{prop}
\begin{proof}
If $n=0$, then the Riemann-Hilbert problem is triangular with identity asymptotics and is trivially solved in closed form by a Cauchy integral,
yielding \eqref{eq:Mzero}.  Thus from now on we consider $n\ge 1$.

The uniqueness of the solution for $n\ge 1$ can be seen from the following
argument.  Continuity of the boundary values taken on $\Sigma$ implies that the ratio of any two solutions of Riemann-Hilbert
Problem~\ref{rhp:M} is an entire function of $z$ that tends to the
identity matrix as $z\rightarrow\infty$.  Uniqueness thus follows by Liouville's Theorem.  

To derive \eqref{eq:Mnotzero}, first note that if ${\bf
M}^n(z)$ solves Riemann-Hilbert Problem~\ref{rhp:M}, then
the first column of ${\bf M}^n(z)$ must be analytic throughout the
$z$-plane.  From the normalization condition (\ref{eq:Mnorm}) it is
then clear that $M_{11}^n(z)$ is a monic polynomial of degree $n$
while $M^n_{21}(z)$ is a polymomial of degree at most $n-1$ (the
leading coefficient of $M^n_{21}(z)$ is not determined from the
normalization condition alone).  The jump condition for the second
column reads
\begin{equation}
\begin{array}{rcl}
M^n_{12+}(e^{i\theta})-M^n_{12-}(e^{i\theta})&=&M^n_{11}(e^{i\theta})e^{-in\theta}\phi(\theta)\\\\
M^n_{22+}(e^{i\theta})-M^n_{22-}(e^{i\theta})&=&M^n_{21}(e^{i\theta})e^{-in\theta}\phi(\theta)\,.
\end{array}
\label{eq:seccoljump}
\end{equation}
In particular, since $\phi(\theta)$ satisfies a H\"older condition, we
may express $M^n_{12}(z)$ as a Cauchy-type integral using
(\ref{eq:seccoljump}).  Thus,
\begin{equation}
M^n_{12}(z)=\frac{1}{2\pi i}\oint_\Sigma
\frac{M^n_{11}(s)s^{-n}}{s-z}\phi(\arg(s))\,ds\,,
\label{eq:M12Cauchy}
\end{equation}
and the normalization condition (\ref{eq:Mnorm}) then requires that
this Cauchy integral be of order $z^{-n-1}$ as $z\rightarrow\infty$ for
each fixed $n\in\mathbb{Z}_+$. Expanding the Cauchy kernel in a geometric
series we see that the following conditions must be satisfied
\begin{equation}
\oint_\Sigma M^n_{11}(z) z^{k-n}\phi(\arg(z)) \,dz = 0
\end{equation}
for $k=0,1,2,\dots,n-1$.  Since $d\theta = dz/(iz)$ for an angular
coordinate $\theta$ on the contour $\Sigma $, this proves that
$M^n_{11}(z)$ is orthogonal to the monomials $1$, $z$, $z^2$, $\dots$,
$z^{n-1}$ with respect to the inner product \eqref{eq:innerproduct}.
The existence of such a monic polynomial of degree $n$ follows from
the Gram-Schmidt algorithm.  Thus, $M^n_{11}(z)=\pi_n(z)$, the $n$th
monic orthogonal polynomial with respect to the weight $\phi(\theta)$
on the unit circle.  

A similar argument applies to the second row of
${\bf M}^n(z)$.  Indeed, $M^n_{21}(z)$ is a polynomial of degree at most $n-1$.  Using (\ref{eq:seccoljump}) we may express $M^n_{22}(z)$ as a Cauchy integral:
\begin{equation}
M^n_{22}(z)=\frac{1}{2\pi i}\oint_\Sigma \frac{M^n_{21}(s)s^{-n}}{s-z}\phi(\arg(s))\,ds
\,,
\label{eq:M22Cauchy}
\end{equation}
and then the normalization condition (\ref{eq:Mnorm}) requires that
$M^n_{22}(z) = z^{-n} + O(z^{-n-1})$ as $z\rightarrow\infty$.  Expanding
the Cauchy kernel in a geometric series, one sees that $M^n_{21}(z)$ is required to satisfy the following conditions:
\begin{equation}
\left\langle M^n_{21}(e^{i\theta}),1\right\rangle_\phi = -1 \hspace{0.2 in}\mbox{and}
\hspace{0.2 in}
\left\langle M^n_{21}(e^{i\theta}),e^{ik\theta}\right\rangle_\phi = 0\,,
\hspace{0.2 in}k=1,\dots,n-1\,.
\label{eq:innerproducts}
\end{equation}
Equivalently, these relations may be written in the form
\begin{equation}
\left\langle e^{i(n-1)\theta},e^{i(n-1)\theta}\overline{M_{21}^n(e^{i\theta})}\right\rangle_\phi = -1 \hspace{0.2 in}\mbox{and}
\hspace{0.2 in}
\left\langle e^{ik\theta},e^{i(n-1)\theta}\overline{M_{21}^n(e^{i\theta})}\right\rangle_\phi = 0\,,
\hspace{0.2 in}k=0,\dots,n-2\,.
\end{equation}
Clearly, the degree $n-1$ polynomial $z^{n-1}\overline{M_{21}^n(1/\overline{z})}$ is orthogonal to the monomials $1, z, \dots,z^{n-2}$ with respect to the inner product $\langle,\rangle_\phi$, and the normalization condition then fixes the leading coefficient:  $z^{n-1}\overline{M_{21}^n(1/\overline{z})}=-\gamma_{n-1}^2z^{n-1} + \dots$.  In other words, we have found that $z^{n-1}\overline{M_{21}^n(1/\overline{z})}=-\gamma_{n-1}^2\pi_{n-1}(z)$, or equivalently,
$M_{21}^n(z)=-\gamma_{n-1}^2z^{n-1}\overline{\pi_{n-1}(1/\overline{z})}$.
\end{proof}
Thus, the conditions of Riemann-Hilbert Problem~\ref{rhp:M} serve to
define the orthogonal polynomials as an alternative to other
representations that may be available, possibly including explicit
contour integral formulae (for special families of weights).

Riemann-Hilbert problems like this one frequently arise as a
consequence of the application of the Fourier transform or
$z$-transform to certain types of linear integral equations ({\em
  i.e.} the Wiener-Hopf technique, see \cite{CarrierKP83}).  Reversing
this sort of reasoning, the representation of the orthogonal
polynomials in terms of Riemann-Hilbert Problem~\ref{rhp:M}
immediately yields integral equations for certain auxiliary unknowns.
Of particular interest are the Marchenko equations obtained in
\cite{GeronimoC79}, in which the unknowns are Fourier coefficients of
functions explicitly related to the orthogonal polynomial $\pi_n(z)$,
and the relevant operator is of the form $I-K$, where $K$ is an
integral operator acting in $\ell^2(n+1,n+2,\dots,\infty)$ with a
kernel that depends explicitly on the weight $\phi$ but not otherwise
on the degree $n$ of the polynomial in question.  Furthermore, the
kernel does not depend in any crucial way on the smoothness of the
weight.  Therefore, in principle, this formulation makes possible the
calculation of asymptotics for the polynomials (indeed, this is one of
the applications of the Marchenko equations discussed in
\cite{GeronimoC79}) in a way that is relatively insensitive to the
analyticity properties of the weight function $\phi$.  However, the
correction terms that appear in such a scheme are necessarily in terms
of infinite Fourier series (see, {\em e.g.}, equation (VI.7) of
\cite{GeronimoC79}) that while having known coefficients are not
convenient for detailed analysis of zeros of the polynomials in
regions of the complex plane where these zeros necessarily arise from
a competition between different terms in an expansion.  To provide
details of the asymptotics, it is more advantageous to work
with the Riemann-Hilbert problem directly.

The uniqueness of ${\bf M}^n(z)$ coupled with symmetry of the
Riemann-Hilbert problem under reflection through the unit circle leads
to the following result.
\begin{prop}
The matrix ${\bf M}^n(z)$ satisfying Riemann-Hilbert Problem~\ref{rhp:M}
satisfies the symmetry relation
\begin{equation}
{\bf M}^n(z) = i^{\sigma_3} \overline{{\bf M}^n(0)^{-1}}
\overline{{\bf M}^n(1/\overline{z})}(-iz)^{\sigma_3}\,.
\end{equation}
In particular, by taking $z=0$ above, we obtain the identity
\begin{equation}
\left(\frac{\gamma_{n-1}}{\gamma_n}\right)^2 = 1-|\alpha_n|^2\,.
\label{eq:gammaalpha}
\end{equation}
\label{prop:gammaalpha}
\end{prop}

\section{Fixed Weights}
\label{sec:strong}
\subsection{Asymptotic behavior of orthogonal polynomials and related quantities as $n\rightarrow\infty$.}
\label{sec:strongasymptoticspi}
In this section we describe several asymptotic results in the theory
of orthogonal polynomials with fixed weights on the unit circle that
we will obtain as a fundamental illustration of our method.  It will
be convenient to introduce the real-valued function $V:S^1\rightarrow
\mathbb{R}$ such that
\begin{equation}
\phi(\theta)=e^{-V(\theta)}\,,\hspace{0.2 in}\text{for all $\theta\in S^1$}\,.
\label{eq:phiV}
\end{equation}

The fundamental object of the asymptotic theory for the fixed weight $\phi$ is the so-called
Szeg\H{o} function:
\begin{equation}
S_\phi(z):=\exp\left(-\frac{1}{2\pi i}\oint_\Sigma\frac{V(\arg(s))\,ds}
{s-z}\right)\,,\hspace{0.2 in}z\not\in \Sigma\,.
\label{eq:szegofunc}
\end{equation}
This is a function analytic for $|z|\neq 1$ that decays to zero as
$z\rightarrow\infty$.  Its value at $z=0$ has the interpretation of
the geometric mean of the weight $\phi(\theta)$:
\begin{equation}
S_\phi(0)=\exp\left(-\frac{1}{2\pi}\int_{-\pi}^\pi V(\theta)\,d\theta
\right) = e^{-V_0}\,.
\label{eq:geometricmean}
\end{equation}
If $\phi(\theta)$ satisfies a H\"older continuity condition, then by strict
positivity so does $V(\theta)$.  In this case, by the Plemelj
formula \cite{Musk}, we have
\begin{equation}
\lim_{r\uparrow 1}S_\phi(re^{i\theta}) = \phi(\theta)\lim_{r\downarrow 1}
S_\phi(re^{i\theta})\,.
\end{equation}
Furthermore, recalling the negative frequency component $N(z)$ of $V$ defined
by \eqref{eq:negfreq}, we can obtain the following equivalent characterization
of $S_\phi(z)$.  Since the function
$\overline{N(1/\overline{z})}$ is analytic for $|z|<1$, it follows
that
\begin{equation}
S_\phi(z)=\left\{\begin{array}{ll}e^{N(z)}\,,&\hspace{0.2 in}|z|>1\\\\
e^{-V_0- \overline{N(1/\overline{z})}}\,,&\hspace{0.2 in}|z|<1\,.
\end{array}\right.
\label{eq:SeL}
\end{equation}
Recalling the function $\Omega:S^1\rightarrow\mathbb{R}$ defined
by \eqref{eq:Omegadef}, 
we have 
\begin{equation}
\lim_{r\uparrow 1}S_\phi(re^{i\theta}) \cdot
\lim_{r\downarrow 1}S_\phi(re^{i\theta})  = 
e^{-V_0}e^{i\Omega(\theta)}\,.
\end{equation}

\subsubsection{General theorems.}
The following results hold for weights where $V$ is of class $C^{k-1,1}(S^1)$
with $k\ge 1$.
\begin{theorem}
Suppose that
$\phi(\theta)=e^{-V(\theta)}$ where
$V:S^1\rightarrow\mathbb{R}$ is of class $C^{k-1,1}(S^1)$ with $k\ge
1$.  Then, for
each fixed integer $p$ and for each $\rho>1$ there is a constant
$K_{p,\rho}>0$
such that the estimate
\begin{equation}
\sup_{|z|\ge\rho}
\left|\frac{d^p}{dz^p}\left[\pi_n(z)z^{-n}e^{-N(z)}-1\right]\right|
\le K_{p,\rho}\frac{\log(n)}{n^{2k}}
\label{eq:piasympoutside}
\end{equation}
holds for all $n$ sufficiently large.
\label{thm:pioutside}
\end{theorem}
The constant $K_{p,\rho}$ typically blows up as $\rho\rightarrow 1$,
and only a finite number of derivatives can be controlled.  More
generally,
we have the following result.
\begin{theorem}
Let $p\ge 0$ be a fixed integer.  Suppose that $\phi(\theta)=e^{-V(\theta)}$ where
$V:S^1\rightarrow\mathbb{R}$ is of class $C^{k-1,1}(S^1)$ with $k\ge
2p+1$.  Then there exists a constant $K_p>0$ such that the estimate
\begin{equation}
\sup_{|z| \ge 1}\left|
\frac{d^p}{dz^p}
\left[\pi_n(z)z^{-n}e^{-N(z)}-1\right]\right|\le 
K_p\frac{\log(n)}{n^{k-2p}}
\label{eq:pialloutsideestimate}
\end{equation}
holds for all $n$ sufficiently large.
\label{thm:picircle}
\end{theorem}

\begin{remark}
Note that as a special case of the estimate
\eqref{eq:pialloutsideestimate}
we obtain the following estimate (under the same conditions)
characterizing
the polynomials on the unit circle:
\begin{equation}
\sup_{-\pi<\theta<\pi}\left|
\left(-ie^{-i\theta}\frac{d}{d\theta}\right)^p
\left[\pi_n(e^{i\theta})e^{-in\theta}e^{-N(e^{i\theta})}-1\right]\right|\le 
K_p\frac{\log(n)}{n^{k-2p}}\,.
\label{eq:picircleestimate}
\end{equation}
In fact, the proof of Theorem~\ref{thm:picircle} is to first establish
\eqref{eq:picircleestimate}, from which
the estimate \eqref{eq:pialloutsideestimate} follows 
(with the same constant $K_{p}$) via the
maximum modulus principle.  
\end{remark}

The weakest conditions under which the above theorem provides
large-degree asymptotics are that $\phi$ is a strictly positive weight
that is Lipschitz continuous.  Theorem~\ref{thm:picircle} may be
compared with results reported in the classic monograph of Szeg\H{o}
\cite[\S 12.1]{szego}.  While asymptotics of $\pi_n(z)$ have been established by other methods under weaker
conditions than Lipschitz continuity and strict positivity of the
weight $\phi$, Theorem~\ref{thm:picircle} exhibits clearly the dependence of
the rate of decay of the error on the smoothness of $\phi$, and the number
of derivatives desired.

To our knowledge, the results of Theorem~\ref{thm:picircle} are stronger than those previously known in that they establish the
convergence in a uniform sense.  This leads to the following.
\begin{corollary}
Let $p\ge 0$ be a fixed integer.  Suppose that $\phi(\theta)=e^{-V(\theta)}$
where $V:S^1\rightarrow\mathbb{R}$ is of class $C^{k-1,1}(S^1)$ with
$k\ge 2p+1$.  Then
\begin{equation}
\lim_{n\rightarrow\infty}\frac{1}{n^p}\cdot\frac{\|\pi_n^{(p)}\|_\phi}{
\|\pi_n\|_\phi} = 1\,.
\end{equation}
\end{corollary}
\begin{proof}
This follows directly from Theorem~\ref{thm:picircle}.  Indeed,
upon carrying out the differentiation in \eqref{eq:picircleestimate},
and combining this estimate with its analogues for all smaller values of $p$, we learn that $|\pi_n^{(p)}(e^{i\theta})|/n^p$ converges uniformly to $|\pi_n(e^{i\theta})|$ for $-\pi<\theta<\pi$.  The proof is then complete  since on $S^1$ uniform convergence implies
convergence in $L^2$.
\end{proof}

\begin{remark} Notice that Theorem \ref{thm:picircle} immediately
implies the following formula valid for all $z$ with $|z| < 1$:
\begin{eqnarray}
\pi_{n}(z) = \frac{1}{2 \pi i} \oint_{|s|=1} \frac{s^{n}
e^{N(s)} + h_{n}(s)}{s-z} ds,
\end{eqnarray}
where 
\begin{eqnarray}
\sup_{|s|=1} |h_{n}(s)| \le K_{p} \frac{\log(n)}{n^{k}}.
\end{eqnarray}
While in principle this could be used to compute asymptotics for
$\pi_{n}(z)$ in this region, more detailed analysis gives the
following improved results.
\end{remark}

\begin{theorem}
Suppose that $\phi(\theta)=e^{-V(\theta)}$ where $V:S^1\rightarrow\mathbb{R}$
is of class $C^{k-1,1}(S^1)$ with $k\ge 1$.  Then for each $\rho$ satisfying 
$0<\rho<1$ there are constants $K^\pm_\rho>0$
such that the estimates
\begin{equation}
\sup_{\rho<|z|<1}\left|\pi_n(z)-z^ne^{-V_0-\overline{N(1/\overline{z})}}
e^{E_kV(r,\theta)} \right|\le K^-_\rho\frac{\log(n)}{n^k}\,,
\label{eq:pinearcircleinsideestimate}
\end{equation}
and
\begin{equation}
\sup_{|z|<\rho}\left|\pi_n(z)\right|\le \frac{K^+_\rho}{n^k}
\label{eq:piawaycircleinsideestimate}
\end{equation}
hold for all $n$ sufficiently large.
\label{thm:piinside}
\end{theorem}

An immediate corollary is that there exists an annulus inside the unit
circle that asymptotically contains no zeros.  That the result we are
about to state in this direction is in a sense sharp will be made
clear when we consider more specific weights below in
\S~\ref{sec:jumpweights} (in particular, see
Corollary~\ref{cor:jumpzeroscircle}).

\begin{corollary}[Zero-free regions]
Suppose that $\phi(\theta)=e^{-V(\theta)}$ where $V:S^1\rightarrow\mathbb{R}$
is of class $C^{k-1,1}(S^1)$ with $k\ge 1$.  Let $\delta>0$ be an arbitrarily
small number.  Then  there are no zeros of $\pi_n(z)$ in the region 
\begin{equation}
\left\{z\Bigg| \log(|z|)>-(k-\delta)\frac{\log(n)}{n}\right\}
\end{equation}
as long as $n$ is sufficiently large.
\label{cor:fixedzerofree}
\end{corollary}

\begin{proof}
  This follows immediately from the estimate
  \eqref{eq:pinearcircleinsideestimate}.  Indeed, since
  $\overline{N(1/\overline{z})}$ and $E_kV(r,\theta)$ are bounded for
  $\rho<r<1$, zeros of $\pi_n(z)$ in the region $\rho<|z|<1$
  necessarily arise from a balance between $z^n$ and a term of uniform
  size $\log(n)/n^k$.  However, $z^n$ is large compared with
  $\log(n)/n^k$ in the region where the inequality
  $\log(|z|)>-(k-\delta)\log(n)/n$ holds.
\end{proof}

A second corollary is an immediate consequence of
\eqref{eq:piawaycircleinsideestimate}.
\begin{corollary}[Recurrence coefficients]
Suppose that $\phi(\theta)=e^{-V(\theta)}$ where
$V:S^1\rightarrow\mathbb{R}$ is of class $C^{k-1,1}(S^1)$ with $k\ge
1$.  Then there is a constant
$K>0$ such that the bound
\begin{equation}
\left|\alpha_n\right|\le 
\frac{K}{n^k}
\end{equation}
holds for sufficiently
large $n$.
\label{cor:alpha}
\end{corollary}

\begin{proof}
This follows directly from \eqref{eq:piawaycircleinsideestimate} with
the use of the identity $\alpha_n=\pi_n(0)$.
\end{proof}

Finally, we have the following result concerning the asymptotic
behavior of the normalization constants.
\begin{theorem}
Suppose that $\phi(\theta)=e^{-V(\theta)}$ where
$V:S^1\rightarrow\mathbb{R}$ is of class $C^{k-1,1}(S^1)$ with $k\ge
1$.  Then there is a constant
$K>0$ such that the bound 
\begin{equation}
\left|\gamma_n^2e^{-V_0} - 1\right|\le K
\frac{\log(n)}{n^{2k}}
\label{eq:fixedgammaasymp}
\end{equation}
holds for sufficiently
large $n$.
\label{thm:gamma}
\end{theorem}
We give a direct proof of this theorem based on the identity
$\gamma_{n-1}^2 = -M_{21}^n(0)$ in \S~\ref{sec:fixedinside}.  However,
another proof with a less sharp error estimate may be based upon
Theorem~\ref{thm:picircle} because on $S^1$ uniform convergence
implies convergence in $L^2$.  Thus, since $\|p_n(z)\|_\phi=1$ and
$p_n(z)=\gamma_n\pi_n(z)$,
\begin{equation}
\gamma_n^2 = \left(\frac{1}{2\pi}\int_{-\pi}^\pi 
|\pi_n(e^{i\theta})|^2\phi(\theta)\,d\theta\right)^{-1}\,.
\end{equation}
Using Theorem~\ref{thm:picircle}, one finds that
\begin{equation}
\gamma_n^2 = \left(\frac{1}{2\pi}\int_{-\pi}^\pi\left|
e^{N(e^{i\theta})}
\right|^2\phi(\theta)\,d\theta\right)^{-1} +
O\left(\frac{\log(n)}{n^k}\right)\,.
\label{eq:gamman2half}
\end{equation}
Next, using \eqref{eq:SeL}, we have
\begin{equation}
\left|e^{N(e^{i\theta})}\right|^2 = 
e^{N(e^{i\theta})}e^{\overline{N(e^{i\theta})}} = e^{V(\theta) 
-V_0} = \frac{e^{-V_0}}{\phi(\theta)}\,.
\end{equation}
Substitution into \eqref{eq:gamman2half} completes the alternate proof.

\begin{remark}
  At this point it is important to comment that the $\dbar$ method we
  develop below in \S~\ref{sec:dbarmethodfixed} yields new formulae
  for the polynomial $\pi_{n}(z)$ (see, for example,
  \eqref{eq:pinoutsideexact}).  The formulae are semi-explicit, in
  that they are written in terms of the solution of a $\dbar$ problem
  (or, equivalently, in terms of the solution of an integral
  equation). This $\dbar$ problem is arrived at after a sequence of
  explicit transformations, and we prove that this problem has a
  unique solution, which possesses an asymptotic expansion for
  $n\rightarrow \infty$.  In general, the terms in this expansion can
  be estimated (from above).  Such estimations give rise to the
  general results described in this subsection.  However, in the
  situation that some further information about the weight function
  $e^{-V}$ is known, it is frequently possible to obtain much more
  precise information about the terms in the asymptotic expansion.  To
  illustrate what can be obtained from an analysis of the terms of the
  expansion, we consider in the following subsection a slightly more
  specific family of weights, and present a rather complete
  description of the pointwise asymptotic behavior of the polynomials.
\end{remark}

\subsubsection{More specific weights.}
\label{sec:jumpweights}
While the estimate \eqref{eq:piawaycircleinsideestimate} allows one to
bound the recurrence coefficients, it does not provide an asymptotic
description of the polynomial $\pi_n(z)$ for $z$ bounded within the
unit circle.  In particular, \eqref{eq:piawaycircleinsideestimate} is
insufficient for deducing the location of the zeros.  With further
assumptions on the regularity of $V(\theta)$ we can extract a leading
term that paves the way for further analysis of $\pi_n(z)$ outside the
zero-free region, but within the unit disk.

\begin{theorem}
  Suppose that $\phi(\theta)=e^{-V(\theta)}$ where
  $V:S^1\rightarrow\mathbb{R}$ is of class $C^{k-1,1}(S^1)$ with $k\ge
  2$. Suppose further that $V^{(k)}(\theta)$ is piecewise continuous
  with $\ell<\infty$ jump discontinuities at points
  $-\pi\le\theta_1<\theta_2<\cdots<\theta_\ell<\pi$, of magnitudes
\begin{equation}
\Delta_j^{(k)}:=\lim_{\theta\downarrow\theta_j}V^{(k)}(\theta)-
\lim_{\theta\uparrow\theta_j}V^{(k)}(\theta)\,.
\end{equation}
Let $V^{(k)}(\theta)$ have one Lipschitz continuous derivative between
consecutive jump discontinuities.  Then, for each $\epsilon>0$, $\sigma>0$,
and $\delta>0$,
the estimate
\begin{equation}
\sup_{\log(|z|)<-(k-\sigma)\log(n)/n}
\left|n^{k+1}e^{-\overline{N(1/\overline{z})}}\pi_n(z)-
n^{k+1}z^ne^{-V_0-2\overline{N(1/\overline{z})}}e^{E_kV(r,\theta)}
B(\log(|z|)/\epsilon)
-
f_n(z)\right|\le\delta\,,
\label{eq:pingeneralmodeljumps}
\end{equation}
holds with
\begin{equation}
f_n(z):=\frac{i^{k+1}}{2\pi}\sum_{j=1}^\ell 
\Delta_j^{(k)}e^{i\Omega(\theta_j)}\frac{e^{i(n+1)\theta_j}}
{e^{i\theta_j}-z}\,,
\label{eq:fndef}
\end{equation}
for all $n$ sufficiently large.
\label{thm:piinsidejumps}
\end{theorem}
Note that $f_n(z)$ is a rational function of $z$ with poles at the
$\ell$ points of discontinuity of $V^{(k)}(\theta)$ on the unit
circle, and with $\ell-1$ zeros which may lie anywhere in the complex
plane, and fluctuate about as $n$ is varied.

With this result, we can completely characterize the zeros of
$\pi_n(z)$ under the same assumptions on $V(\theta)$.
The simplest example of orthogonal polynomials on the unit circle is
of course the case $\phi(\theta)\equiv 1$, in which case
$\pi_n(z)=z^n$ for all $n\ge 0$.  Here we see that all zeros of
$\pi_n(z)$ lie exactly at $z=0$.  In particular, the zeros avoid the
unit circle $|z|=1$.  This situation is typical for strictly positive
analytic weights, in which case it is known that the zeros of
$\pi_n(z)$ asymptotically lie within a smaller disk
$|z|\le\rho<1$. Here, the nearest singularity $z_0$ to the unit circle
of the analytic continuation through $|z|=1$ of the function
$S_\phi(z)$ from the domain $|z|>1$ determines the radius $\rho$ by
$\rho=|z_0|$.  However, such confinement of the zeros within the
circle is no longer typical once one leaves the analytic class.  For
example, discontinuities in any derivatives of $\phi(\theta)$ make it
possible for at most a finite number of zeros to be bounded away from
the unit circle while all remaining zeros converge to the unit circle,
as the following corollaries of Theorem~\ref{thm:piinsidejumps} show.
For each $M>0$ let 
\begin{equation}
\begin{array}{rcl}
\displaystyle F_n^+(M)&:=&\displaystyle
\{\text{$z\in\mathbb{C}$ such that 
$\log(|f_n(z)|)>M$}
\}\,,\\\\
\displaystyle F_n^-(M)&:=&\displaystyle
\{\text{$z\in\mathbb{C}$ such that 
$\log(|f_n(z)|)<-M$}
\}\,,\\\\
\displaystyle F_n^0(M)&:=&\displaystyle
\{\text{$z\in\mathbb{C}$ such that 
$|\log(|f_n(z)|)|\le M$}
\}\,.
\end{array}
\end{equation}

\begin{corollary}[Zeros near the unit circle]
\label{cor:jumpzeroscircle}
Assume the same hypotheses as in Theorem~\ref{thm:piinsidejumps}.
Let $A_n(\sigma)$ denote the annulus
\begin{equation}
A_n(\sigma):=\left\{z\Bigg|-(k+1)\frac{\log(n)}{n}-\frac{\sigma}{n}<
\log(|z|)<-(k-\sigma)\frac{\log(n)}{n}\right\}\,.
\end{equation}  
Then, for each $\sigma>0$ there is some $M>0$ such that
the region $A_n(\sigma)\cap F^-_n(M)$ contains no zeros of
$\pi_n(z)$ for sufficiently large $n$.  

For each $M>0$, the zeros of $\pi_n(z)$ in the region
$A_n(\sigma)\cap (F^0_n(M)\cup F^+_n(M))$ satisfy
\begin{equation}
|z|=1-(k+1)\frac{\log(n)}{n}+\frac{1}{n}\log(|f_n(z)|) + o\left(\frac{1}{n}\right)\,,
\label{eq:zerosmod}
\end{equation}
and
\begin{equation}
\theta=-\frac{1}{n}\Omega(\theta)+\frac{1}{n}\arg(f_n(z))+\frac{\pi}{n} +
o\left(\frac{1}{n}\right)
\label{eq:zerosarg}
\end{equation}
modulo $2\pi/n$, where $\theta=\arg(z)$, and in both cases the error
term is uniformly small in the specified region.  It follows that the
angular spacing between neighboring zeros of $\pi_n(z)$ in the
specified region is $\Delta\theta=2\pi/n + o(1/n)$.

For any fixed $M>0$, \eqref{eq:zerosmod} can be rewritten uniformly in
the region $A_n(\sigma)\cap F^0_n(M)$ as 
\begin{equation}
|z|=1-(k+1)\frac{\log(n)}{n}+ O\left(\frac{1}{n}\right)\,,
\end{equation}
and consequently there exists some $\alpha\in (0,\sigma)$ such that
the zeros of $\pi_n(z)$ in the region $A_n(\sigma)\cap F^0_n(M)$
asymptotically lie between the two circles $|z|=1-(k+1)\log(n)/n \pm
\alpha/n$.  

If $M>0$ is sufficiently large, then the region $F^+_n(M)$ is
contained in a disjoint union of small discs centered at the poles
$e^{i\theta_j}$ of the rational function $f_n(z)$.  In this situation,
let $F^+_{n,j}(M)$ denote the component of $F^+_n(M)$ near the pole
$e^{i\theta_j}$.  Then from \eqref{eq:zerosmod} we see that the zeros
of $\pi_n(z)$ in the region $A_n(\sigma)\cap F^+_{n,j}(M)$ satisfy
\begin{equation}
|z|=1-(k+1)\frac{\log(n)}{n}-\frac{1}{n}\log|z-e^{i\theta_j}| + 
O\left(\frac{1}{n}\right)\,,
\end{equation}
where the error term is uniform in the specified region,
which indicates that zeros are attracted to a curve that ``bulges''
outward from the circle $|z|=1-(k+1)\log(n)/n$ in a region of angular
width proportional to $\log(n)/n$ centered at the point
$e^{i\theta_j}$, to a maximum radius defined by the equation 
\begin{equation}
1-|z|=k\frac{\log(n)}{n}+\frac{\log(\log(n))}{n} + O\left(\frac{1}{n}\right)\,.
\end{equation}
Note that this radius is just within the inner boundary of the zero-free
annulus described by Corollary~\ref{cor:fixedzerofree}.
\end{corollary}

\begin{proof}
  The annulus $A_n(\sigma)$ converges toward the unit circle as
  $n\rightarrow\infty$, and therefore if $z$ is a zero of $\pi_n(z)$ and
$n$ is large enough, Theorem~\ref{thm:piinsidejumps} gives
\begin{equation}
n^{k+1}z^ne^{-V_0-2\overline{N(1/\overline{z})}}e^{E_kV(r,\theta)}+f_n(z) = 
o(1)\,,
\label{eq:basiczerofinder}
\end{equation}
as $n\rightarrow\infty$ because for such $z$, $B(\log(|z|)/\epsilon)=1$ and
$\pi_n(z)=0$.  The $o(1)$ error term is uniformly small for all zeros in
$A_n(\sigma)$.
Now, let 
\begin{equation}
c:=\inf_{|z|\le 1}\left|e^{-V_0-2\overline{N(1/\overline{z})}}
e^{E_kV(r,\theta)}\right|
\end{equation}
and note that $c>0$ due to the assumptions in force on $V$.  Since
\begin{equation}
\inf_{z\in A_n(\sigma)} n^{k+1}|z|^n = e^{-\sigma}\,,
\end{equation}
we see that \eqref{eq:basiczerofinder} is inconsistent for large
enough $n$ if $|f_n(z)|<ce^{-\sigma}$.  Therefore, given $\sigma>0$,
$A_n(\sigma)\cap F^-_n(M)$ contains no zeros as $n\rightarrow\infty$
as long as $M>\sigma-\log(c)$.

For any $M>0$ we now consider those zeros $z$ of $\pi_n(z)$ in the region
$A_n(\sigma)\cap (F_n^0(M)\cup F_n^+(M))$, in which case we may divide through
in \eqref{eq:basiczerofinder} by $f_n(z)$ to obtain
\begin{equation}
\frac{n^{k+1}z^n}{f_n(z)}e^{-V_0-\overline{N(1/\overline{z})}}e^{E_kV(r,\theta)}+1 = o(1)\,.
\label{eq:zerofinderrewrite}
\end{equation}
Consistency requires that $n^{k+1}z^n/f_n(z) = O(1)$, and then since
the $n$-independent exponential factors are continuous up to the unit
circle and $A_n(\sigma)$ is converging to the unit circle, we may
replace these factors by their limiting values on the unit circle
without changing the error estimate.  Therefore,
\eqref{eq:zerofinderrewrite} becomes
\begin{equation}
\frac{n^{k+1}z^n}{f_n(z)}e^{i\Omega(\theta)} + 1 = o(1)\,,
\end{equation}
as $n\rightarrow\infty$ uniformly for those zeros $z$ of $\pi_n(z)$
that lie in the annulus $A_n(\sigma)$.  This proves both
\eqref{eq:zerosmod} and \eqref{eq:zerosarg}.  
\end{proof}

\begin{remark}
  Note that the presence of the outward ``bulges'' in the zero curve
  near the points of discontinuity of $V^{(k)}(\theta)$ indicates the
  sharpness of the zero-free region established for more general
  weights in Corollary~\ref{cor:fixedzerofree}.
\end{remark}

While most zeros of $\pi_n(z)$ move toward the unit circle as
$n\rightarrow\infty$ under the hypotheses of
Theorem~\ref{thm:piinsidejumps}, there may be at most $\ell-1$ zeros
further inside the unit circle, which correspond to zeros of $f_n(z)$.
We refer to these as ``spurious zeros''.

\begin{corollary}[Spurious zeros]
  Assume the same hypotheses as in Theorem~\ref{thm:piinsidejumps}.
  For each $M>0$, there exists a $\sigma>0$, such that the zeros of
  $\pi_n(z)$ lying in the disk $\log(|z|)\le-(k+1)\log(n)/n-\sigma/n$
  also lie in the set $F_n^-(M)$ for sufficiently large $n$.
  Moreover, whenever $\epsilon_n$ is a sequence of positive numbers
  such that $\epsilon_n\rightarrow 0$ as $n\rightarrow\infty$, the zeros of
  $\pi_n(z)$ in the disk $|z|\le 1-(k+1)\log(n)/n-1/(n\epsilon_n)$ satisfy
\begin{equation}
f_n(z)=o(1)
\end{equation}
as $n\rightarrow\infty$.  In particular, $\pi_n(z)$ has exactly one
zero for each zero of $f_n(z)$ in this region, making at most $\ell-1$
spurious zeros.
\end{corollary}

\begin{proof}
Define a constant $C>0$ by
\begin{equation}
C:=\sup_{|z|\le 1}\left|e^{-V_0-2\overline{N(1/\overline{z})}}
e^{E_kV(r,\theta)}B(\log(|z|)/\epsilon)\right|\,.
\end{equation}
Therefore,
\begin{equation}
\sup_{\log(|z|)\le -(k+1)\log(n)/n-\sigma/n} \left|n^{k+1}z^n
e^{-V_0-2\overline{N(1/\overline{z})}}
e^{E_kV(r,\theta)}B(\log(|z|)/\epsilon)\right|\le Ce^{-\sigma}\,.
\end{equation}
It then follows easily from \eqref{eq:basiczerofinder} that
all zeros of $\pi_n(z)$ lying in the disk where the inequality
$\log(|z|)\le -(k+1)\log(n)/n-\sigma/n$ holds will also lie in the
set $F_n^-(M)$ for large enough $n$ whenever $\sigma>M+\log(C)$.

The term in \eqref{eq:basiczerofinder} proportional to $n^{k+1}z^n$ is
$o(1)$ as $n\rightarrow\infty$ uniformly for $z$ in the disk delineated by 
the inequality $|z|\le 1-(k+1)\log(n)/n-1/(n\epsilon_n)$ with
$\epsilon_n\rightarrow 0$ as $n\rightarrow\infty$.  Therefore, $f_n(z)=o(1)$
for zeros in this disk, and the one-to-one correspondence of zeros of $f_n(z)$
with spurious zeros of $\pi_n(z)$ in this region follows from the Implicit
Function Theorem.
\end{proof}

\begin{remark}
  Note that the zeros of $f_n(z)$ play an apparently contradictory
  role in the asymptotics.  Indeed, zeros of $f_n(z)$ that occur near
  the unit circle repel zeros of $\pi_n(z)$, while each zero of
  $f_n(z)$ that occurs far enough within the unit circle attracts
  precisely one zero of $\pi_n(z)$.
\end{remark}

\begin{remark}
  A careful reading of the proof of Theorem~\ref{thm:piinsidejumps}
  (see \S~\ref{sec:fixedinside}) shows that the $o(1)$ estimate that
  is stated in \eqref{eq:pingeneralmodeljumps} can be improved to give
  a rate of decay, and that even better estimates can be obtained if
  one does not insist on uniformity.  These simple improvements can
  provide, for example, decay rate information for the error terms in
  the description of the spurious zeros.  We have opted not to give
  these slightly improved estimates in the interest of simplicity of
  presentation.
\end{remark}

\subsubsection{Numerical computation of zeros of $\pi_n(z)$ when
  derivatives of $V$ have jump discontinuities.}  To illustrate the
detailed asymptotic behavior of the zeros of $\pi_n(z)$ explained
above, we have carried out some numerical experiments.  Let us fix
$\ell$ angles of discontinuity
$\{\theta_1,\dots,\theta_\ell\}\subset(-\pi,\pi)$ by the formula
\begin{equation}
\theta_j:=\frac{2\pi }{\ell}\left(j-\frac{1}{2}-\frac{\ell}{2}\right)\,,
\hspace{0.2 in}\text{for $j=1,\dots,\ell$}\,.
\end{equation}
Consider the family of weights $\phi(\theta)$ given by the formula
\begin{equation}
\phi(\theta):=\left\{
\begin{array}{ll}
\displaystyle
1+e^{w_j}\left|\sin\left(\frac{\ell}{2}(\theta-\theta_j)
\right)\right|^k\,,&\hspace{0.2 in}\text{for $\theta_{j-1}<\theta<\theta_j$
with $j=2,3,\dots,\ell$}\,,\\\\
\displaystyle
1+e^{w_1}\left|\sin\left(\frac{\ell}{2}(\theta-\theta_1)\right)
\right|^k\,,&\hspace{0.2 in}\text{for $|\theta|>(\ell-1)\pi/\ell$}\,.
\end{array}
\right.
\label{eq:jumpweight}
\end{equation}
The positive integer $k$ and the real numbers $w_1,\dots,w_\ell$ are free
parameters.  This weight is of the form $\phi=e^{-V}$ where
$V^{(k)}(\theta)$ has jump discontinuities at the points $\theta_j$ of
magnitudes that can be adjusted by choice of the $w_j$.  One advantage
of this family from the point of view of numerical computation is that
the Fourier coefficients of $\phi(\theta)$ can be evaluated
symbolically.  In a package such as Mathematica capable of arbitrary
precision arithmetic, this leads to the possibility of computing the
elements of the Toeplitz matrices (whose minors are assembled to yield
the coefficients of the polynomial $\pi_n(z)$) with sufficient
accuracy for subsequent numerical computation of the zeros when $n$ is
large.  In practice, we computed the coefficients up to an overall
factor by scaling the Fourier coefficients making up the Toeplitz
matrix by $e^{V_0}$.  This is necessary to avoid numerical overflow or
underflow since according to the strong Szeg\H{o} limit theorem the
Toeplitz determinant of $\phi=e^{-V}$ of size $n+1$ scales as
$e^{-nV_0}$.  While the Fourier coefficients can be computed
symbolically, we obtained the coefficients 
$\Delta_j^{(k)}e^{i\Omega(\theta_j)}$
appearing in the rational function $f_n(z)$ defined by
\eqref{eq:fndef} with the help of numerical integration.

\begin{figure}[h]
\begin{center}
\psfig{file=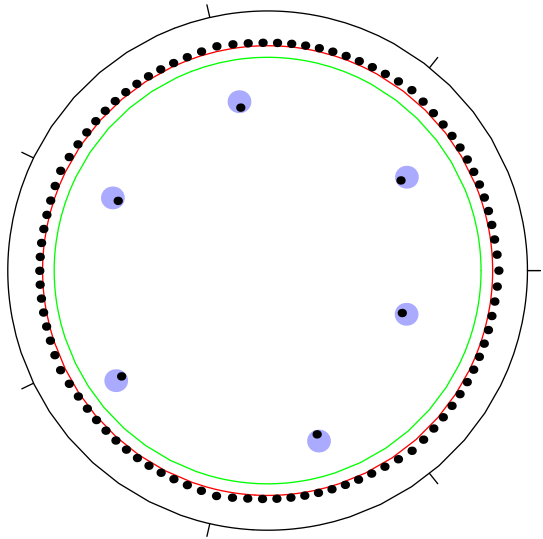,width=3 in}
\end{center}
\caption{\em The output of the Mathematica code for $n=104$, $k=3$,
$\ell=7$, and ${\bf w}^T=(-1,-1/2,-1/4,-1,-1/4,-1,-1/2)$.}
\label{fig:thumbnail}
\end{figure}

The Mathematica code we wrote to carry out these computations is
available from the companion website to this paper:
\OPUC.  The code takes as input the number of jump discontinuities,
$\ell$, the order $k$ of the derivative experiencing the
discontinuites, a vector ${\bf w}$ of length $\ell$ containing the
parameters $w_j$, and the degree $n$ of the polynomial $\pi_n(z)$.
The output is a figure showing the unit circle (black) with exterior
tick marks at the angles $\theta_j$, $j=1,\dots,\ell$, the
zero-attracting circle $|z|=1-(k+1)\log(n)/n$ (green), the inner
boundary circle $|z|=1-k\log(n)/n$ of the zero-free annulus (red), and
the zeros of $f_n(z)$ that occur in a neighborhood of the unit disk
(large lavender dots).  Superimposed on the figure are the zeros of
$\pi_n(z)$ (small black dots).  Sample output from the program is
shown in Figure~\ref{fig:thumbnail}.

\begin{figure}[h]
\psfig{file=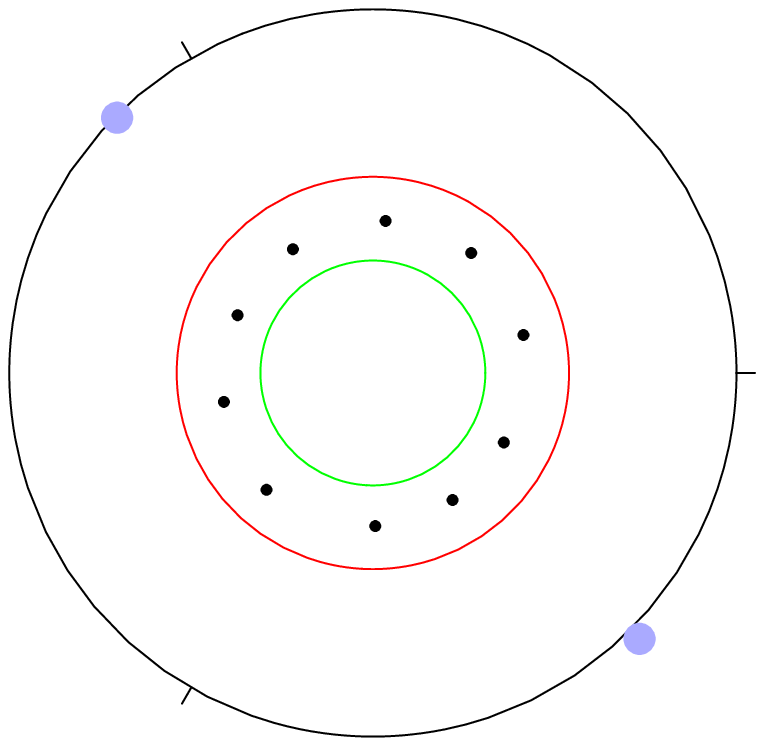,width=1.2 in}
\psfig{file=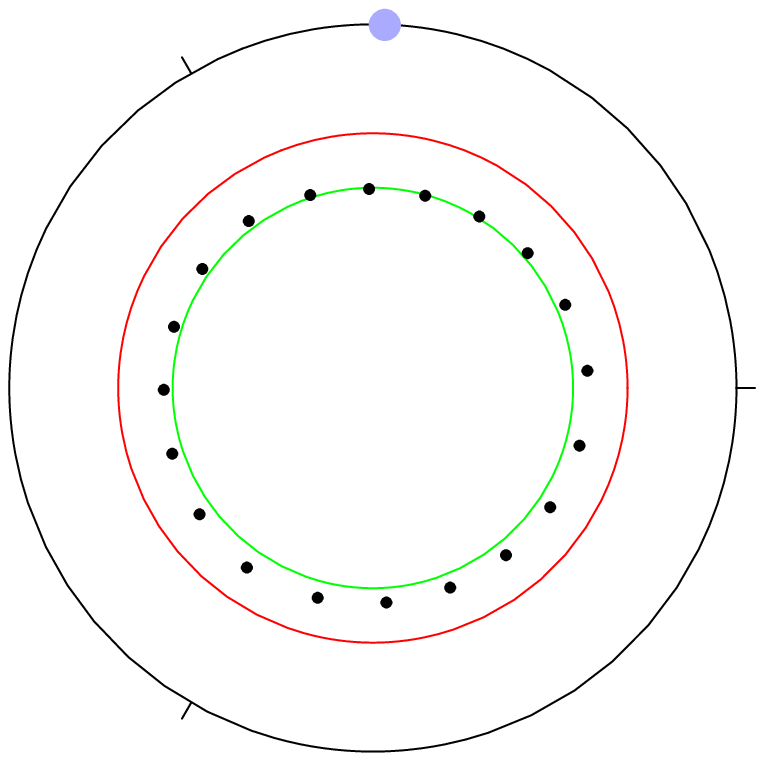,width=1.2 in}
\psfig{file=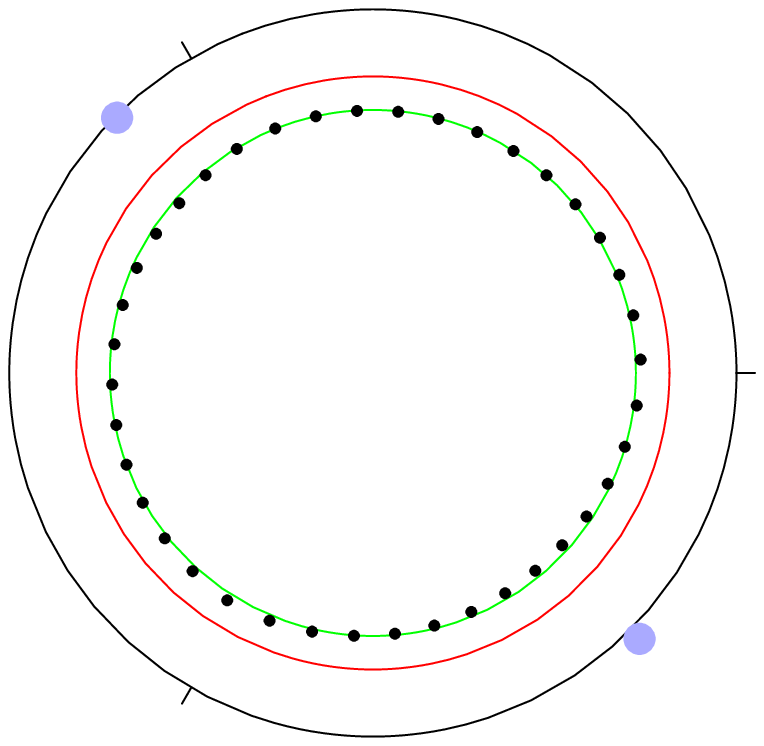,width=1.2 in}
\psfig{file=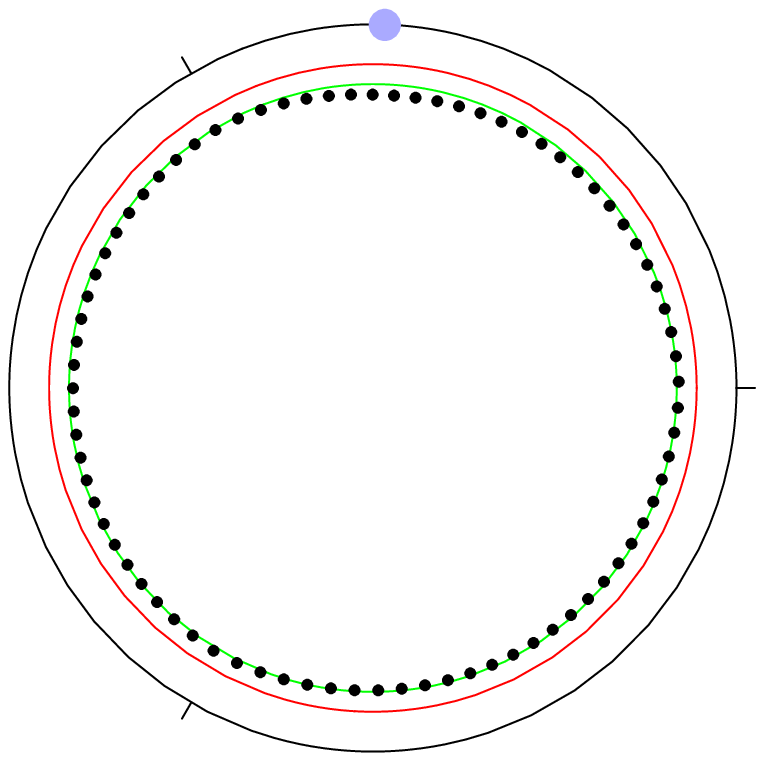,width=1.2 in}
\psfig{file=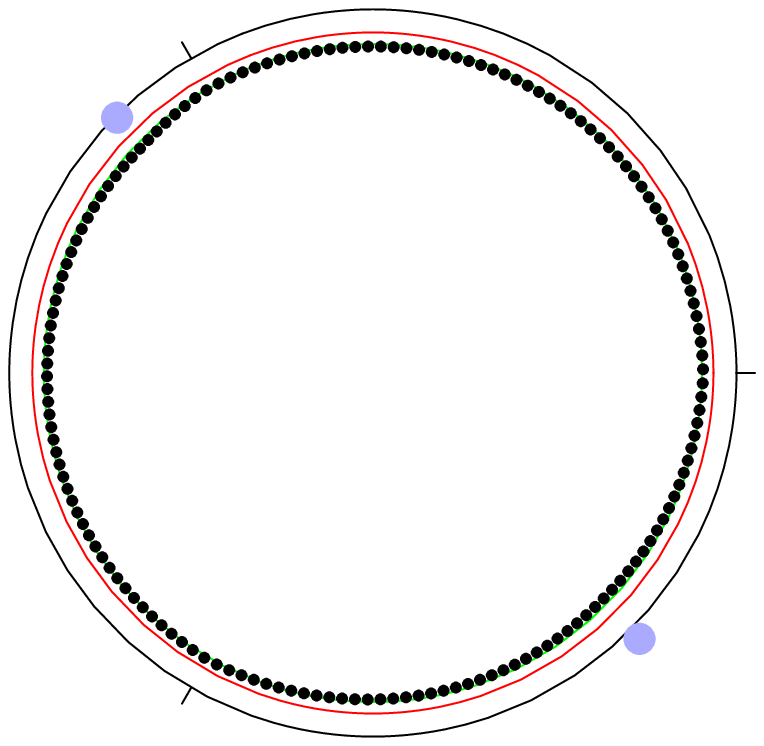,width=1.2 in}
\caption{\em A sequence of plots showing convergence of the zeros of 
  $\pi_n(z)$ toward the curve $|z|=1-(k+1)\log(n)/n$ 
  as $n\rightarrow\infty$.  The parameters are $\ell=3$,
  ${\bf w}^T=(-4,-2,-3)$, and $k=2$.  From left to right, $n=10$,
  $20$, $40$, $80$, and $160$. }
\label{fig:3jumps}
\end{figure}

The first effect we would like to illustrate is the rate of
convergence of the zeros of $\pi_n(z)$ to the unit circle with
increasing $n$.  See Figure~\ref{fig:3jumps}.  The annulus associated
with the inequalities $1-k\log(n)/n<|z|<1$ is asymptotically
zero-free, and the curve $|z|=1-(k+1)\log(n)/n$ asymptotically
attracts the zeros near the unit circle.
The convergence to the zero-attracting circle is clear.  More difficult
to discern from the images is the outward ``bulging'' of the zeros near
the angles $\theta_j$ of discontinuity toward the inner boundary of the
zero-free region.  The imaginative reader can see this effect beginning
in the figure corresponding to $n=160$, but larger values of $n$ (and
a rescaling of the figures near the unit circle) will be necessary to 
resolve the ``bulging'' completely.

\begin{figure}[h]
\mbox{
\psfig{file=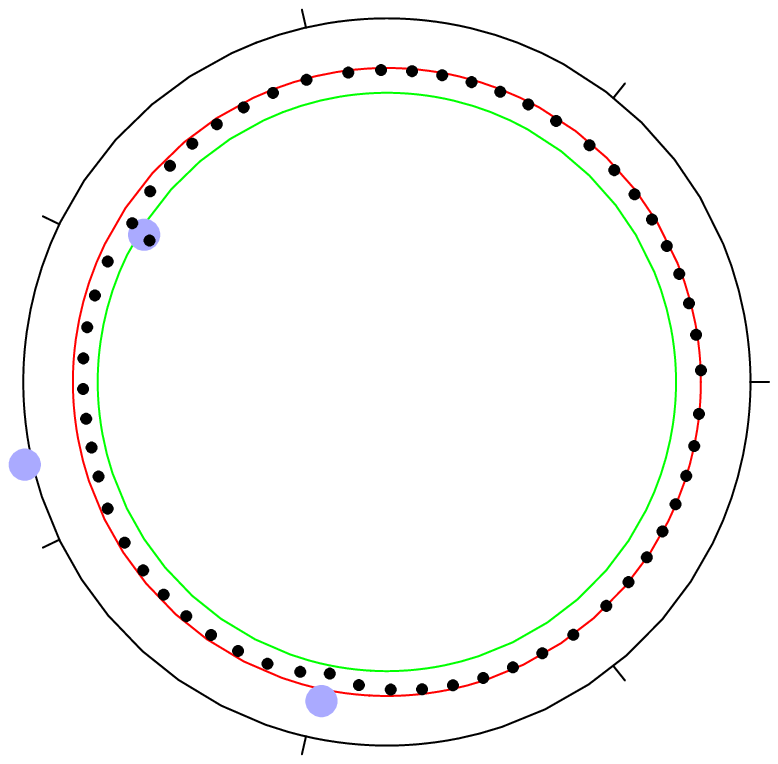,width=0.9 in}
\psfig{file=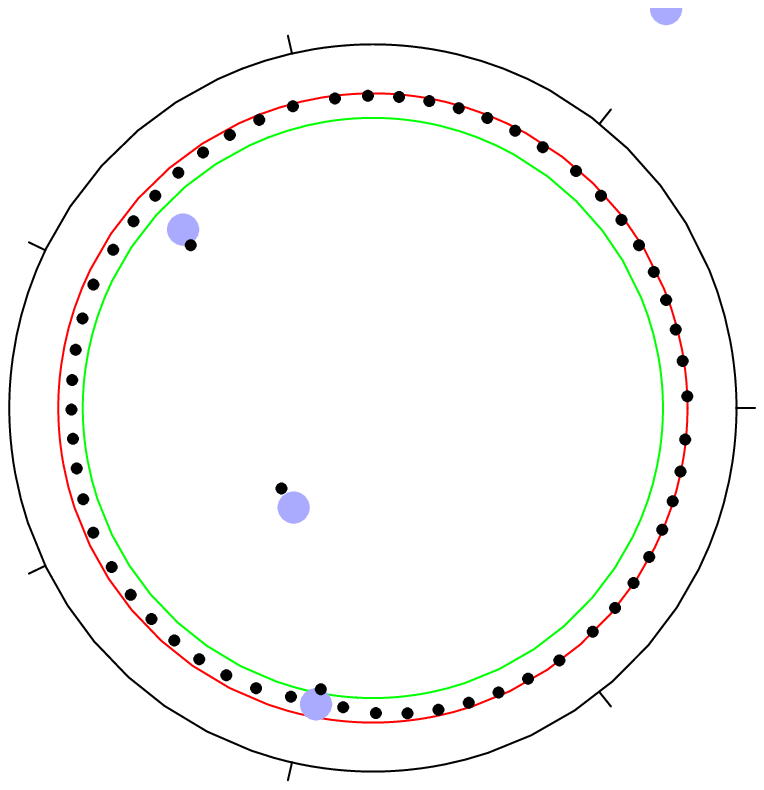,width=0.9 in}
\psfig{file=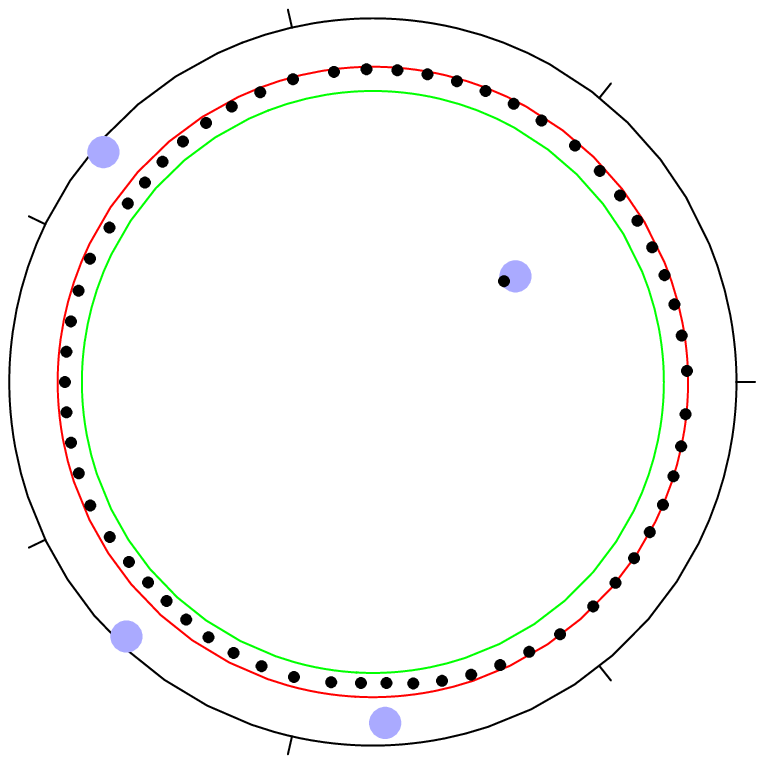,width=0.9 in}
\psfig{file=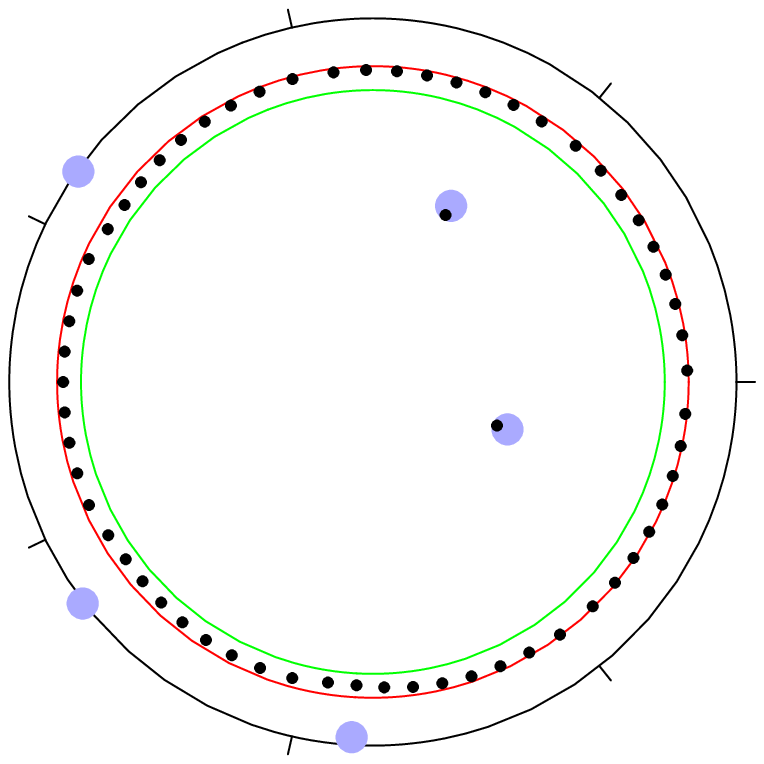,width=0.9 in}
\psfig{file=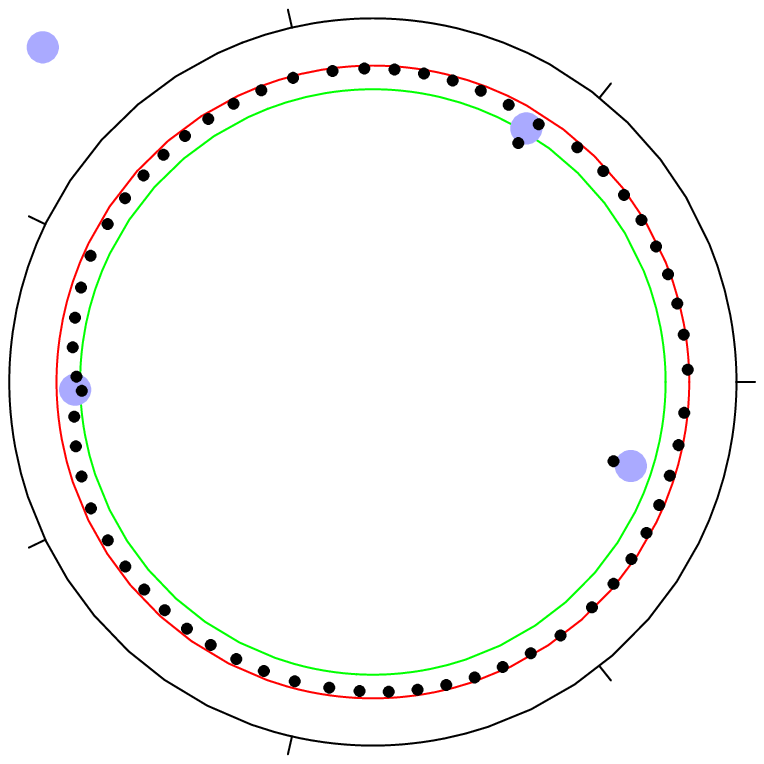,width=0.9 in}
\psfig{file=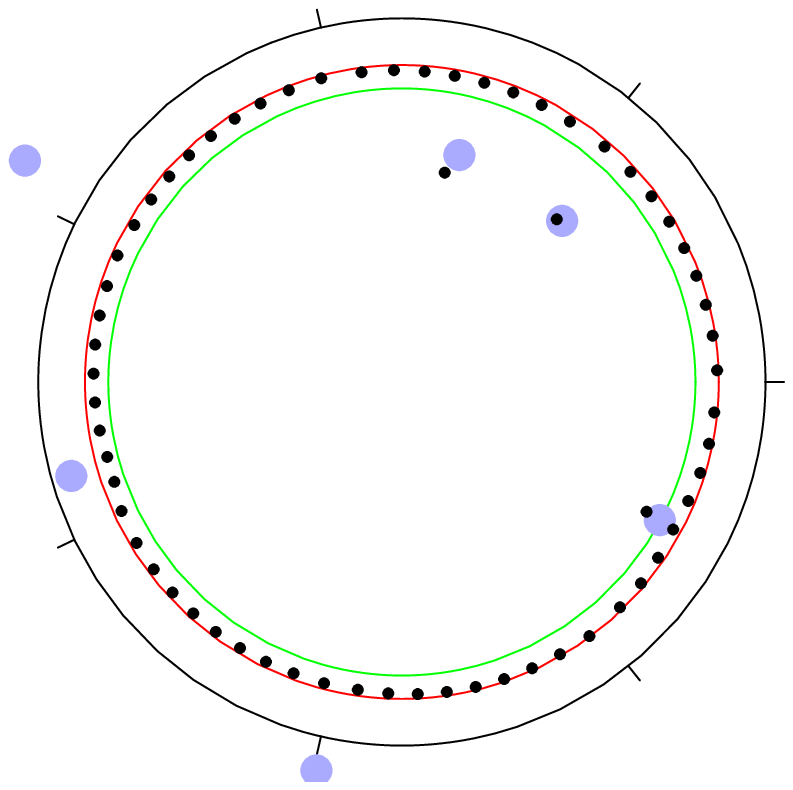,width=0.9 in}
\psfig{file=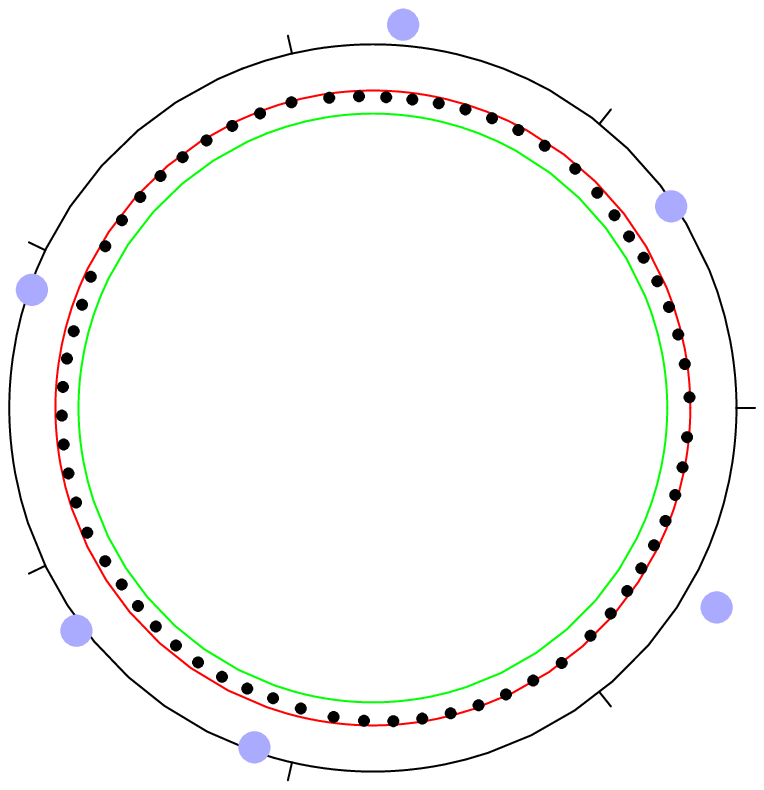,width=0.9 in}
}\\
\mbox{
\psfig{file=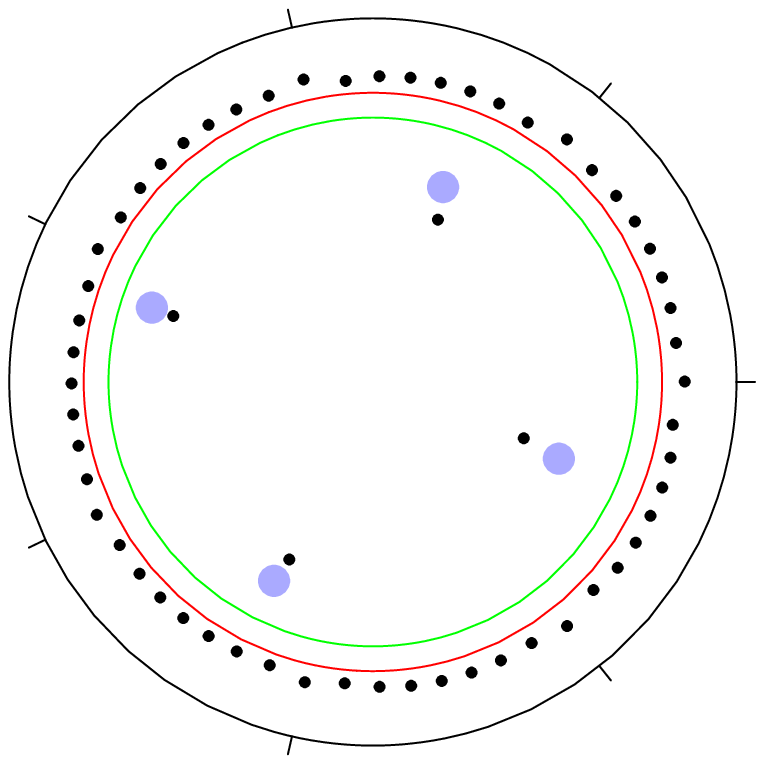,width=0.9 in}
\psfig{file=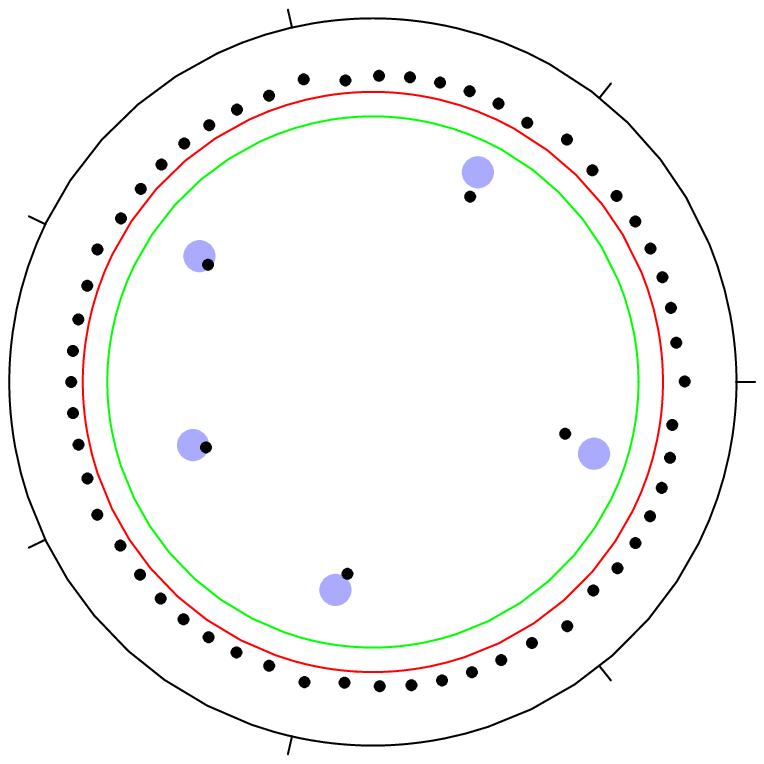,width=0.9 in}
\psfig{file=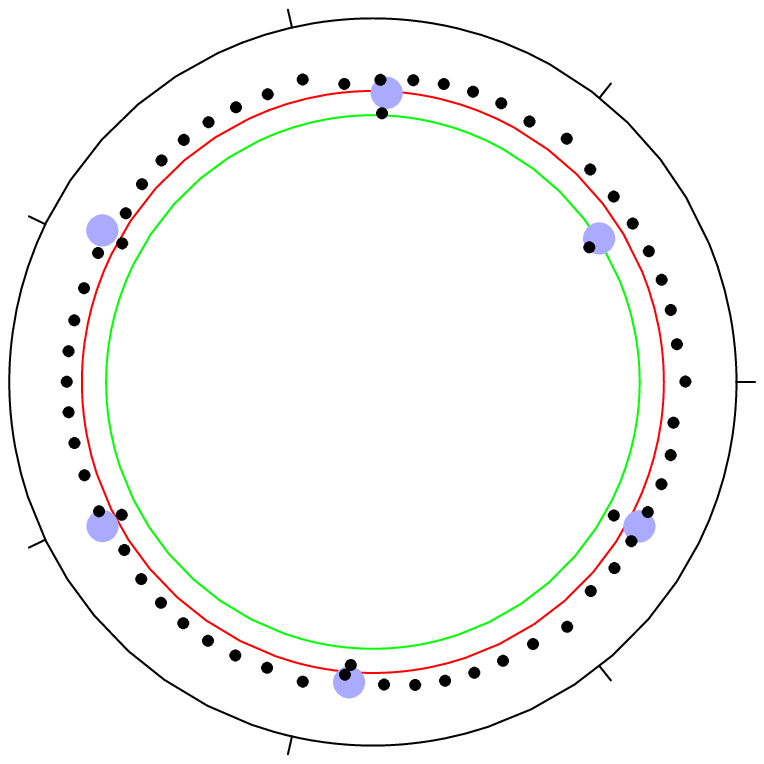,width=0.9 in}
\psfig{file=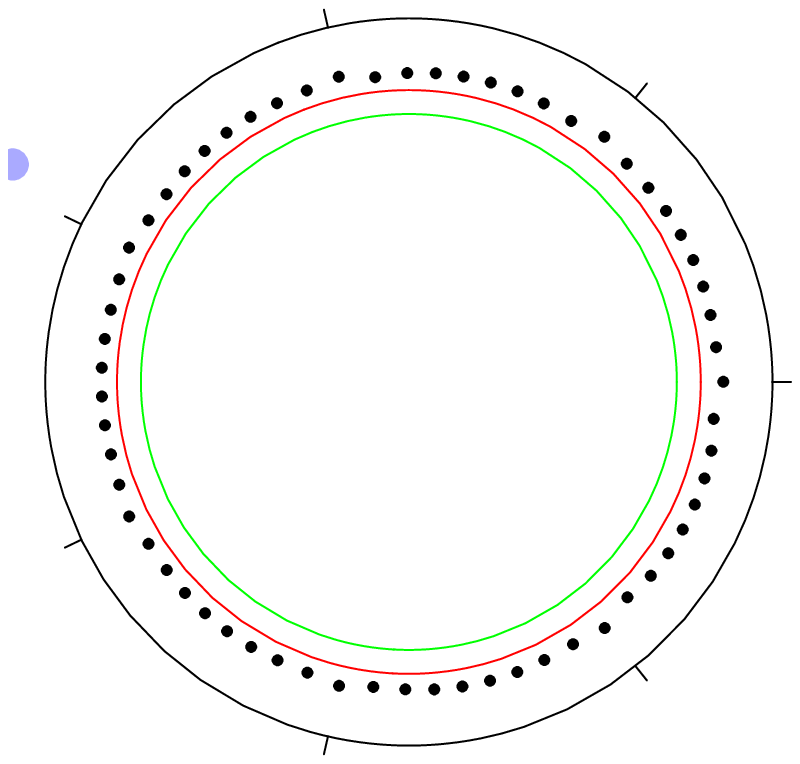,width=0.9 in}
\psfig{file=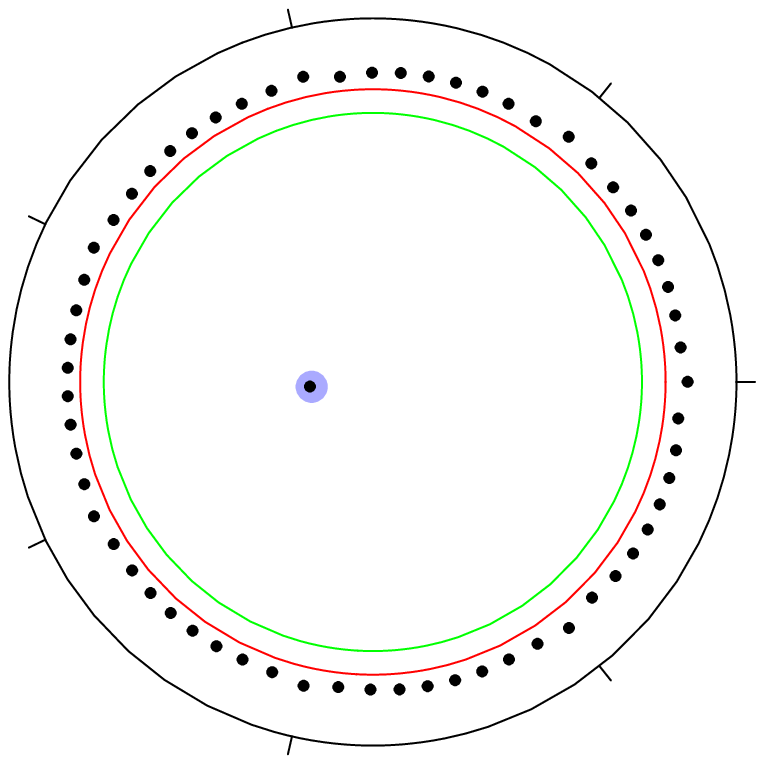,width=0.9 in}
\psfig{file=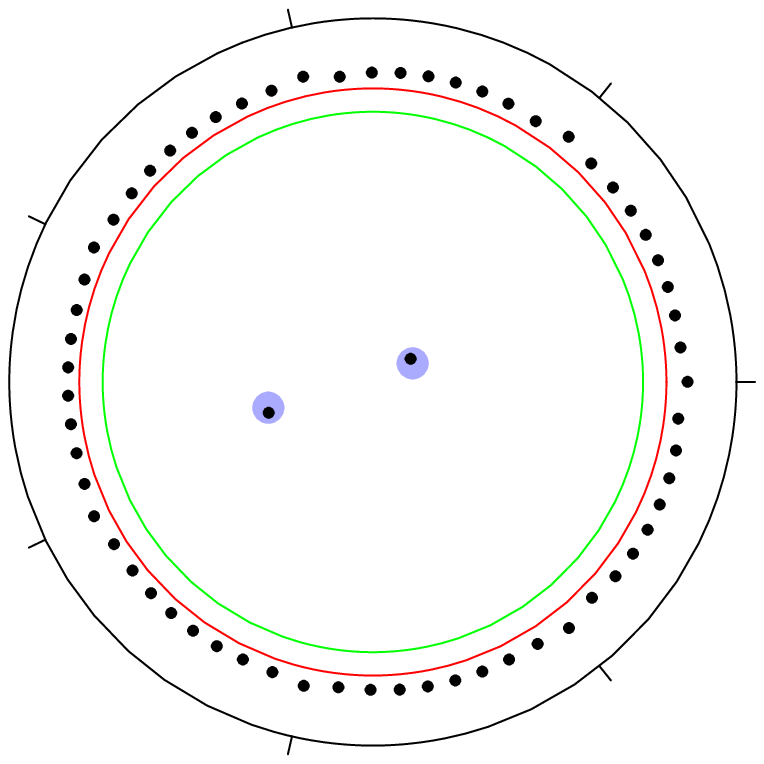,width=0.9 in}
\psfig{file=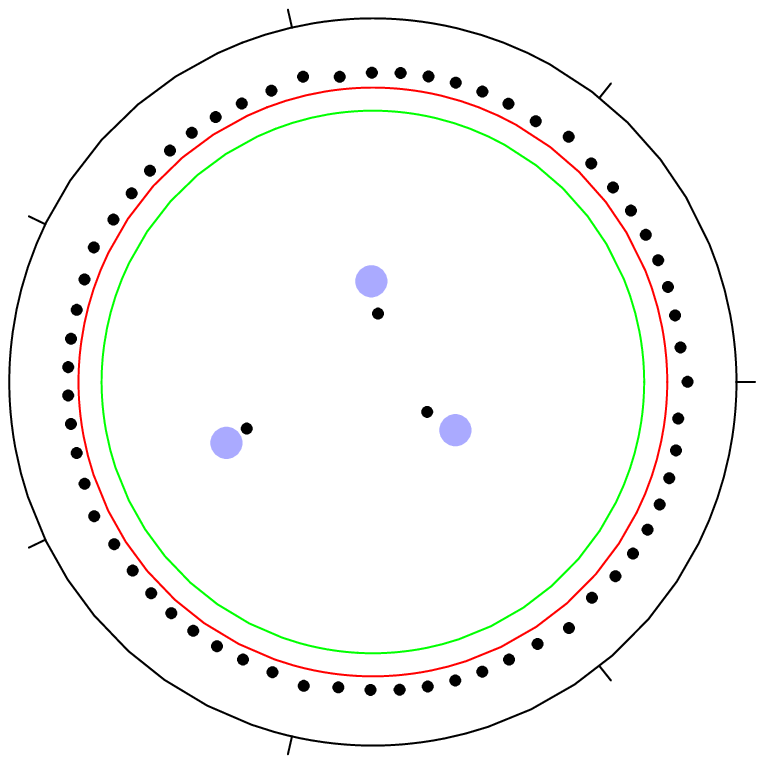,width=0.9 in}
}
\caption{\em Periodic fluctuation of the spurious zeros with respect to $n$.
  Top: $k=2$.  Bottom: $k=3$.  In both cases $n$ varies from $n=60$ to
  $n=66$ from left to right, there are $\ell=7$ equally-spaced points
  of discontinuity, and ${\bf w}^T=(-1,-1/2,-1/4,1/2,-1/4,1/4,-1/2)$.
  Note that while for $k=3$ the majority of zeros appear to lie in the
  zero-free annulus $1-k\log(n)/n<|z|<1$, this is a finite $n$
  effect.}
\label{fig:fluctuation}
\end{figure}

Zeros of $f_n(z)$ play little role for the polynomials whose zeros are
illustrated in Figure~\ref{fig:3jumps}.  Next, we would like to
illustrate the effect zeros of $f_n(z)$ can have on $\pi_n(z)$; this
is the phenomenon of spurious zeros.  Note that in the present case of
equally-spaced angles $\theta_1,\dots,\theta_\ell$, the function
$f_n(z)$ is periodic in $n$ with period $\ell$.  In
Figure~\ref{fig:fluctuation} we present images corresponding to one period
of the function $f_n(z)$ in the case of discontinuities of the second and
third derivatives of $V(\theta)$.
Here it is clear that the $\ell-1$ zeros of $f_n(z)$ fluctuate about
rapidly with $n$, and can be either inside the unit circle or outside.
Each zero of $f_n(z)$ inside the unit circle is an asymptote for exactly
one zero of $\pi_n(z)$, while those outside the unit circle have little effect
on the zeros of $\pi_n(z)$.  The zeros near the unit circle, either inside
or outside, have a repulsive effect on the zeros of $\pi_n(z)$.

\begin{figure}[h]
\mbox{
\psfig{file=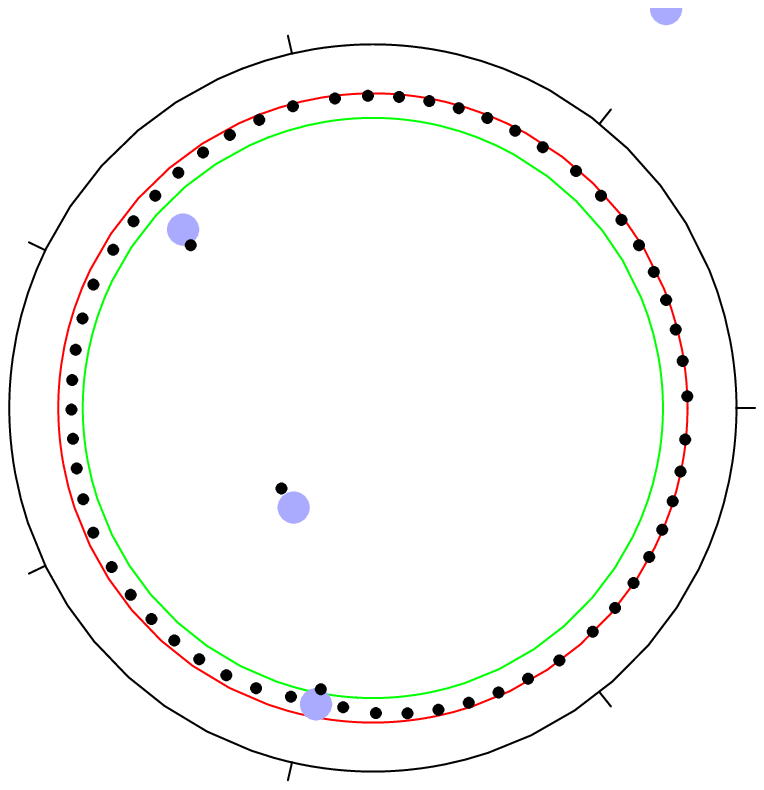,width=0.9 in}
\psfig{file=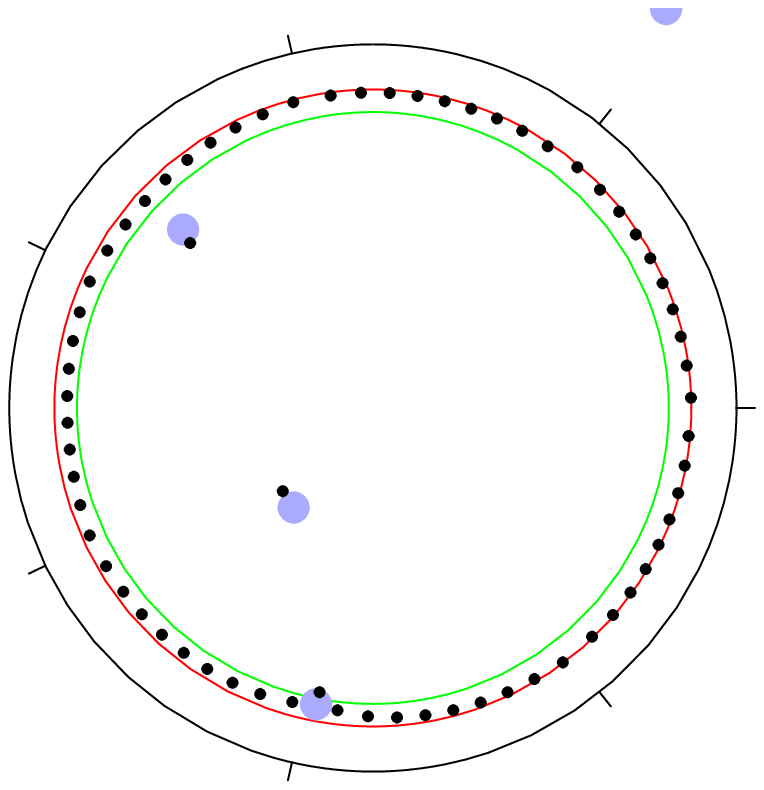,width=0.9 in}
\psfig{file=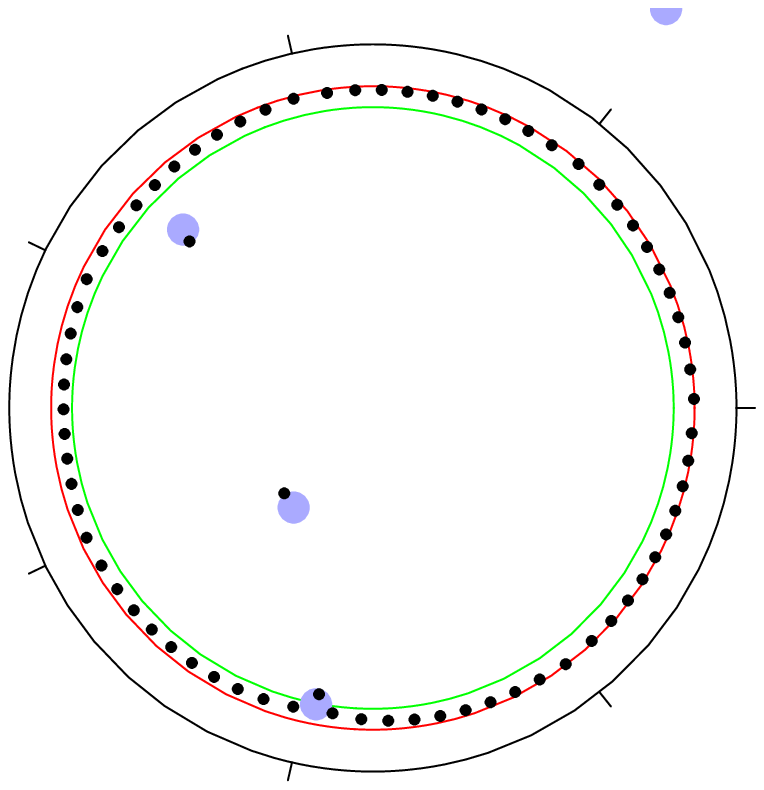,width=0.9 in}
\psfig{file=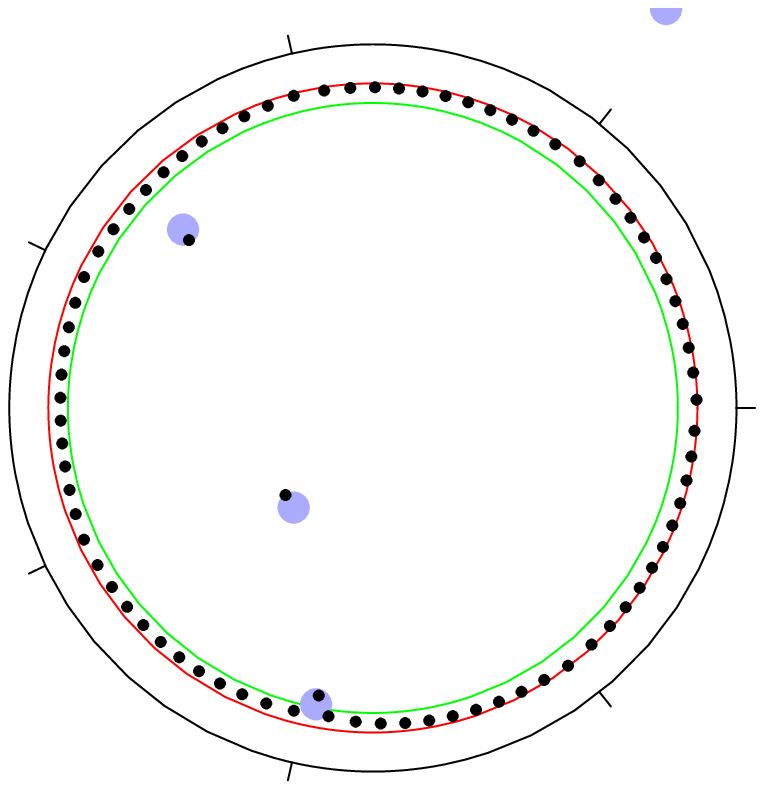,width=0.9 in}
\psfig{file=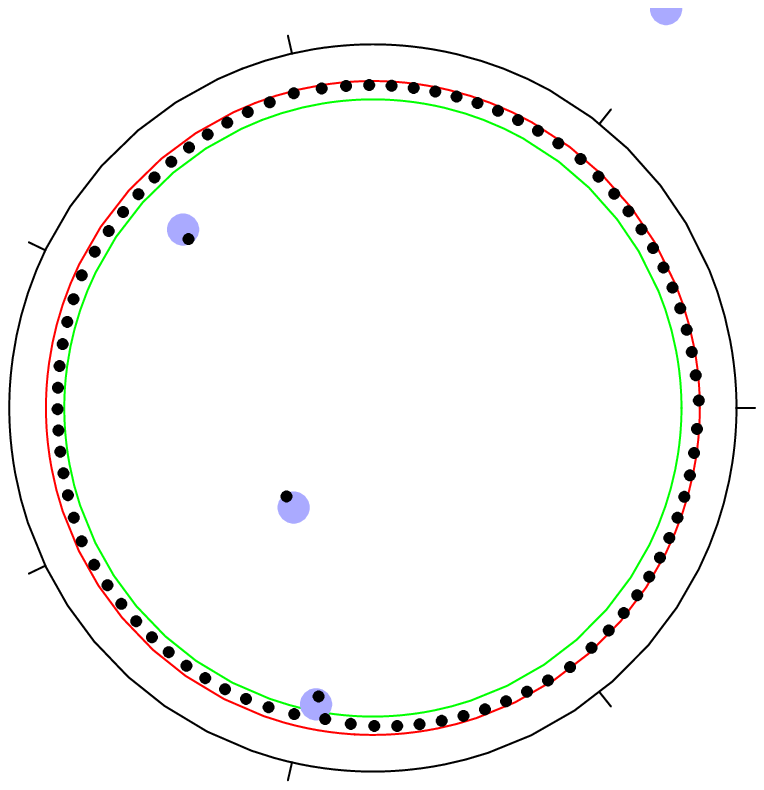,width=0.9 in}
\psfig{file=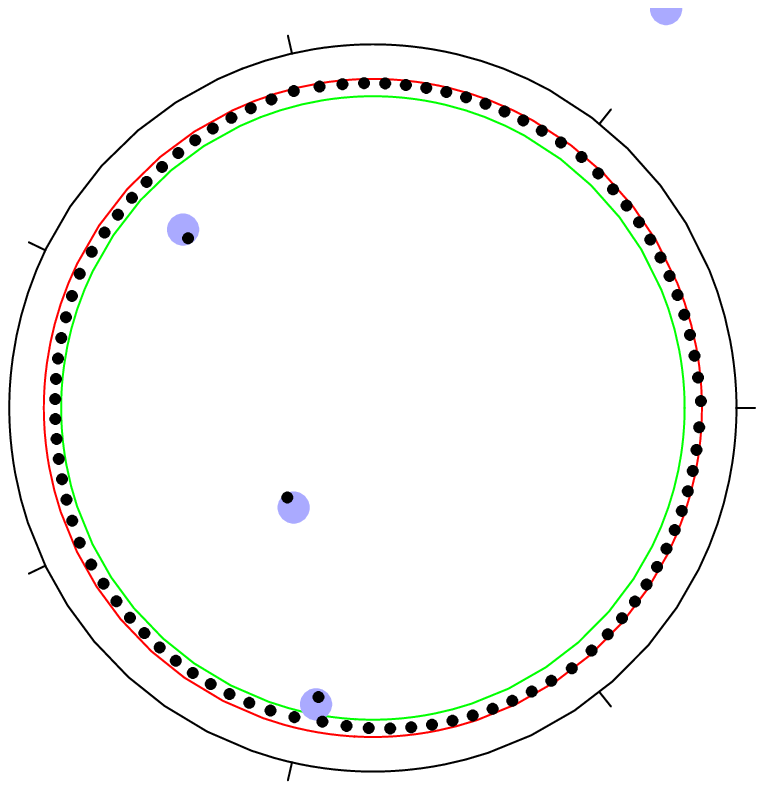,width=0.9 in}
\psfig{file=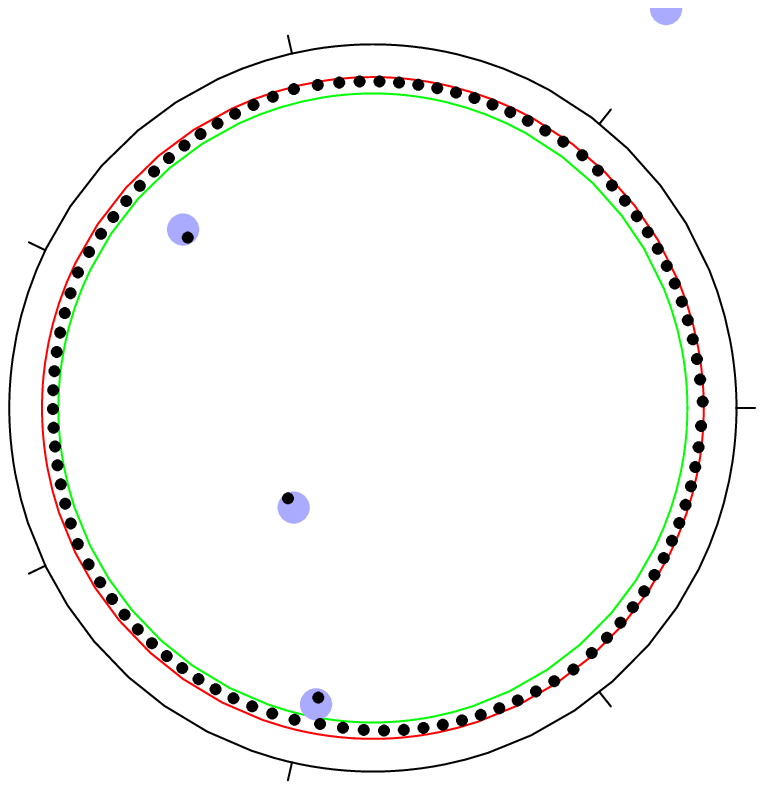,width=0.9 in}}\\
\mbox{
\psfig{file=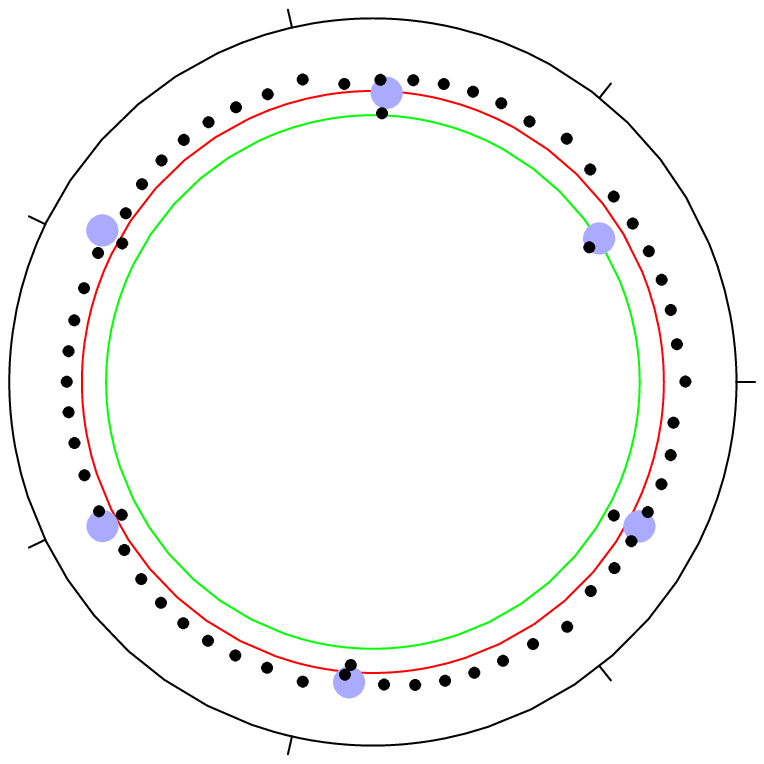,width=0.9 in}
\psfig{file=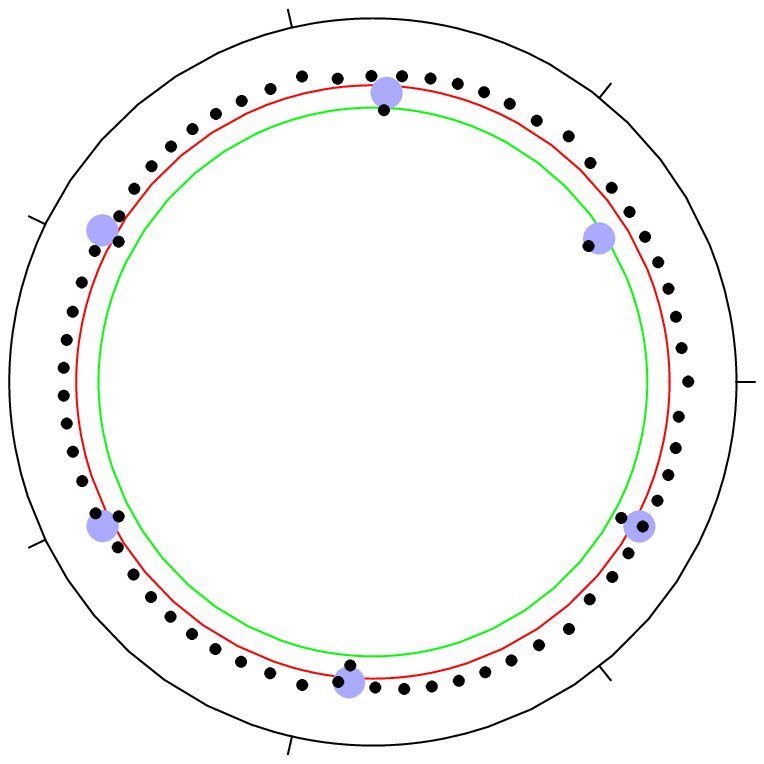,width=0.9 in}
\psfig{file=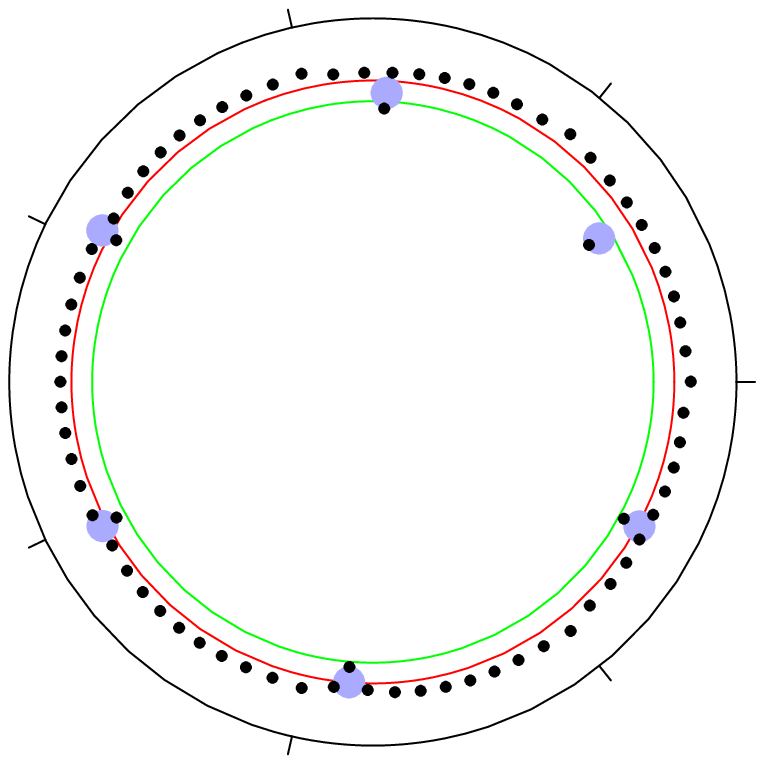,width=0.9 in}
\psfig{file=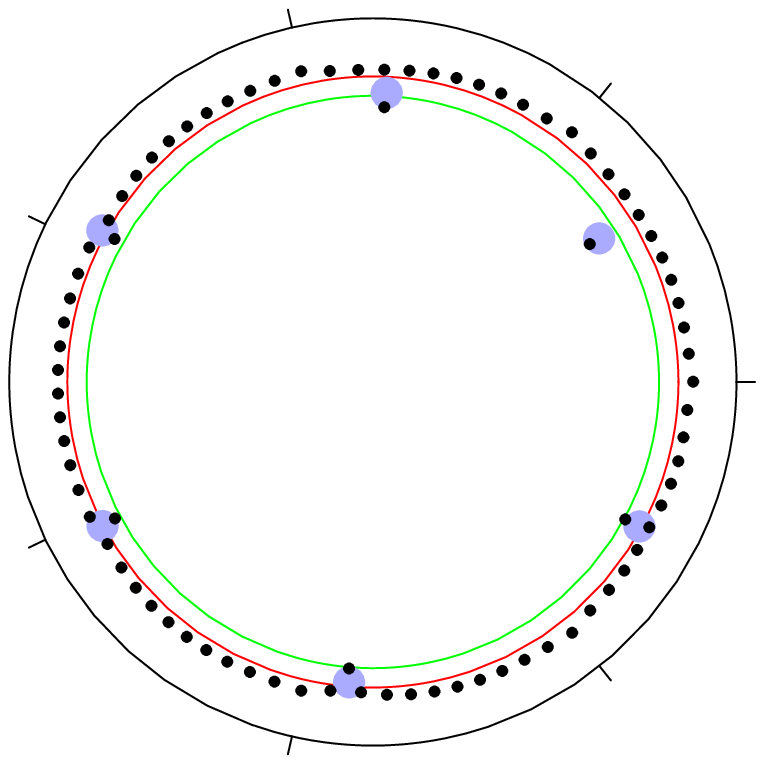,width=0.9 in}
\psfig{file=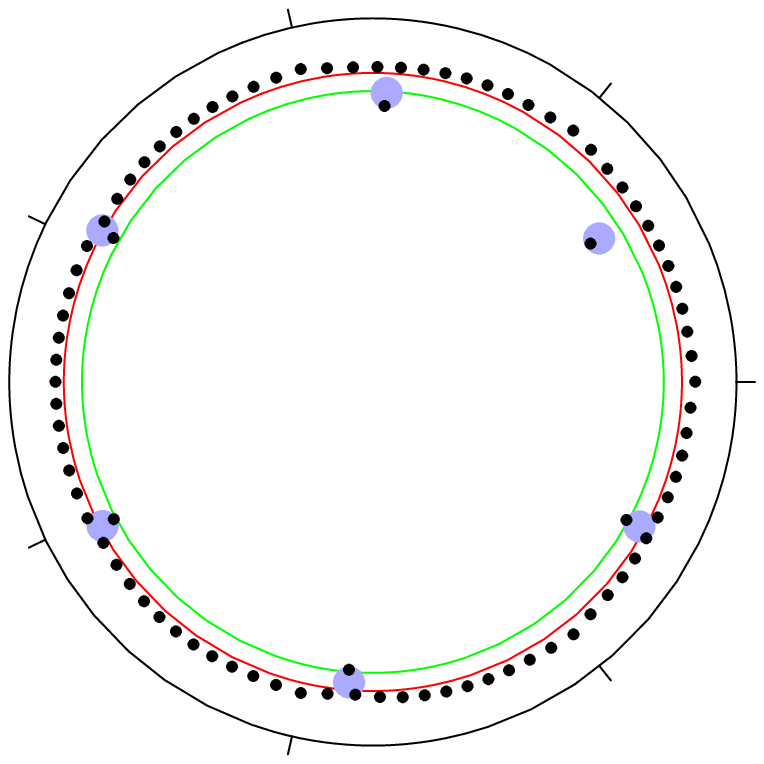,width=0.9 in}
\psfig{file=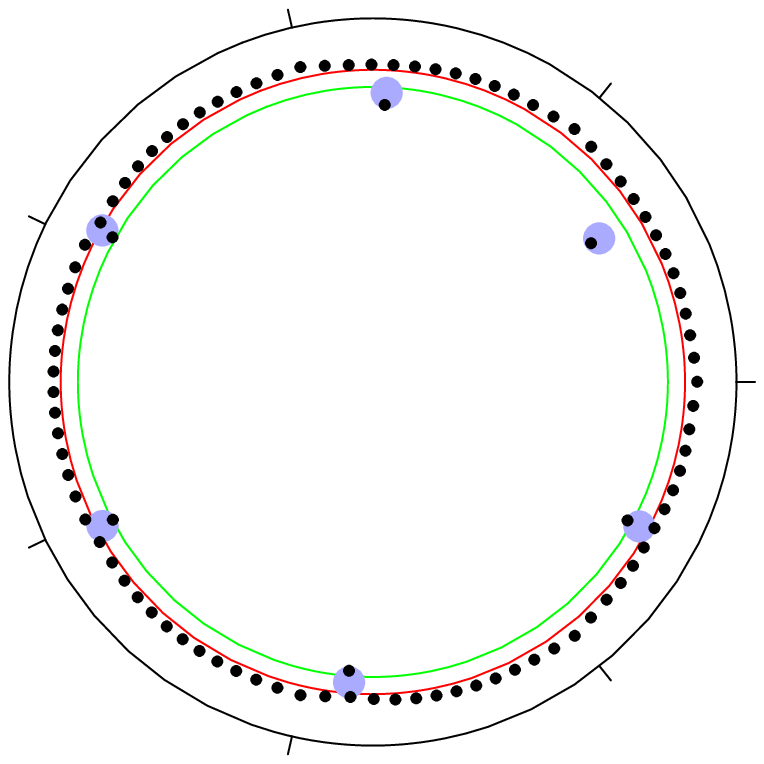,width=0.9 in}
\psfig{file=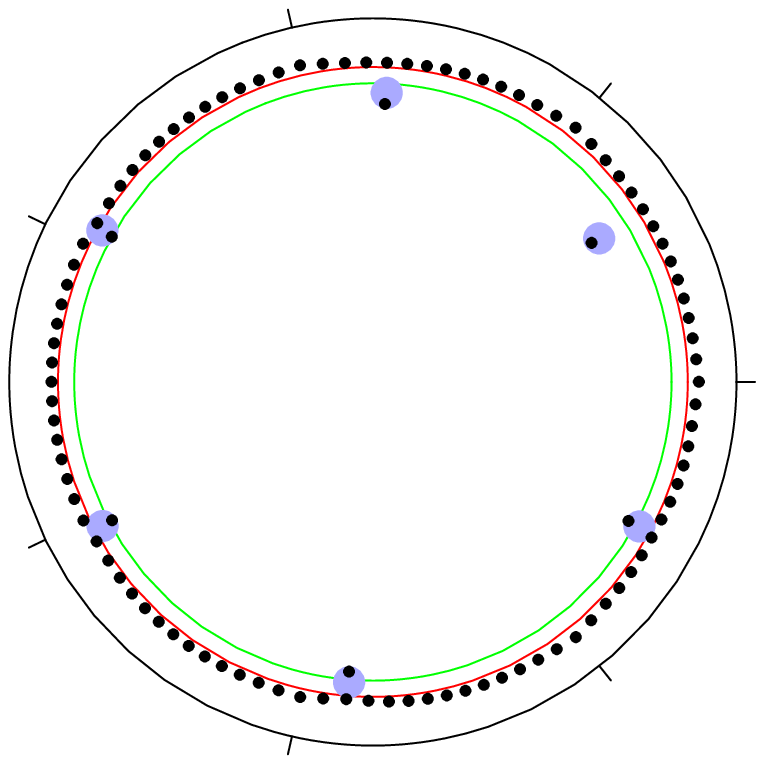,width=0.9 in}
}
\caption{\em Convergence of spurious zeros with increasing $n$.
  Top: $k=2$ and $n=61$, $68$, $75$, $82$, $89$, $96$, and $103$ from
  left to right.  Bottom: $k=3$ and $n=62$, $69$, $76$, $83$, $90$, $97$,
  and $104$ from left to right.  In both cases there are $\ell=7$
  equally-spaced points of discontinuity and ${\bf
    w}^T=(-1,-1/2,-1/4,1/2,-1/4,1/4,-1/2)$.  }
\label{fig:spurious}
\end{figure}

From the images in Figure~\ref{fig:fluctuation} it is not obvious that
the zeros of $f_n(z)$ inside the unit disk attract corresponding
spurious zeros of $\pi_n(z)$ in the limit $n\rightarrow\infty$.  The images
shown in Figure~\ref{fig:spurious} show that this convergence indeed occurs.
Here, we have used the periodicity of $f_n(z)$ to examine the asymptotic
behavior of the zeros of $\pi_n(z)$ along a periodic subsequence of $n$-values
along which the zeros of $f_n(z)$ remain fixed.  

The effect of a zero of $f_n(z)$ upon those of $\pi_n(z)$ is the most
subtle when it occurs near the unit circle.  It should be stressed
that the parameters $w_j$ of the weight under consideration can be
deformed in a continuous manner such that it may always be arranged
that $f_n(z)$ has zeros near the unit circle.  We fixed $n$ and chose
a one-parameter deformation of the $w_j = w_j(t)$ in order to
continuously tune a zero of $f_n(z)$ through the unit circle from
outside to inside.  A movie of this deformation is available at \OPUC,
and several consecutive frames of this movie are shown below in
Figure~\ref{fig:deformation}.
\begin{figure}[h]
  \mbox{ 
    \psfig{file=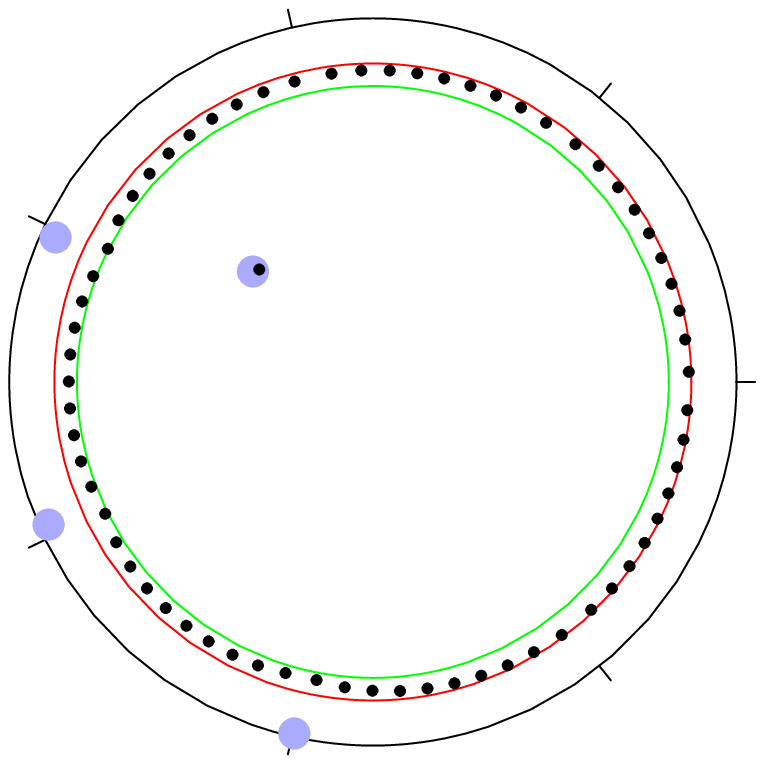,width=1.5 in}
    \psfig{file=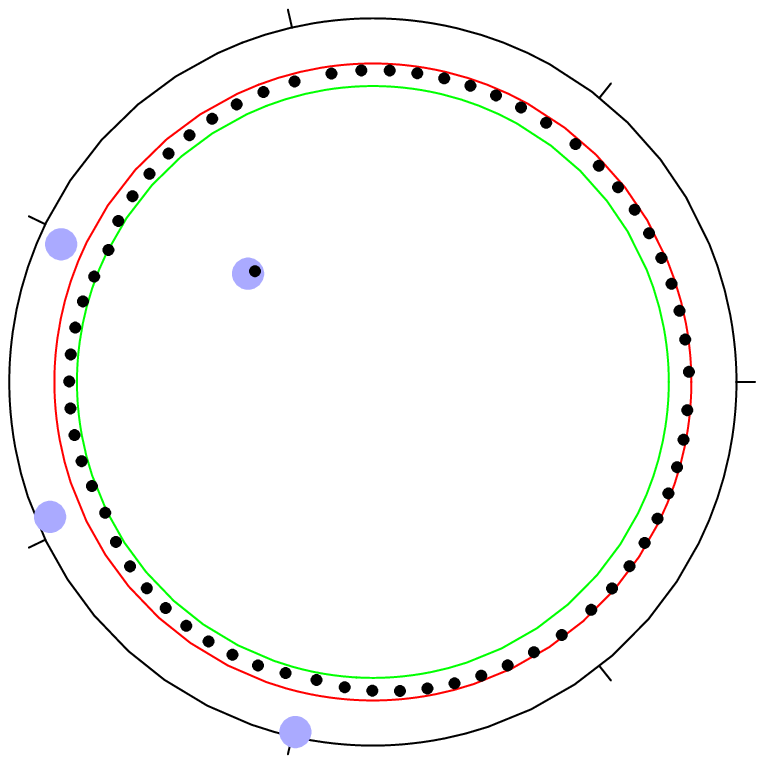,width=1.5 in} 
    \psfig{file=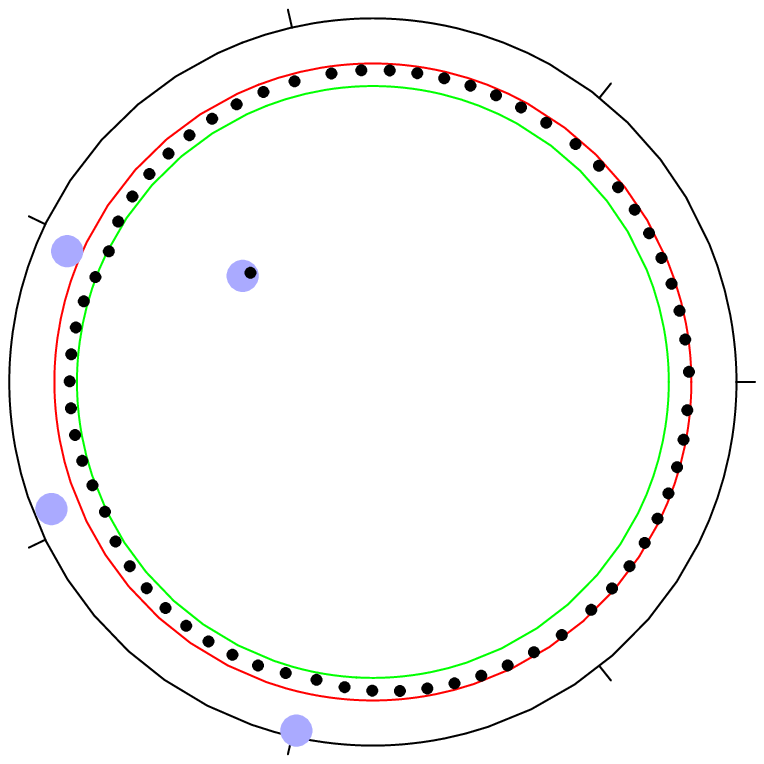,width=1.5 in} 
    \psfig{file=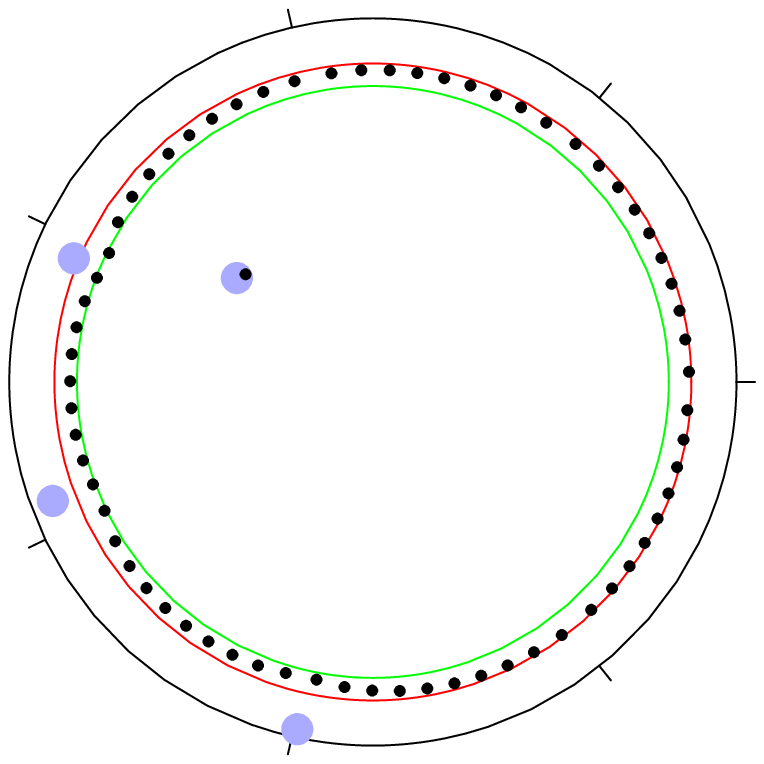,width=1.5 in}
}\\
\mbox{
    \psfig{file=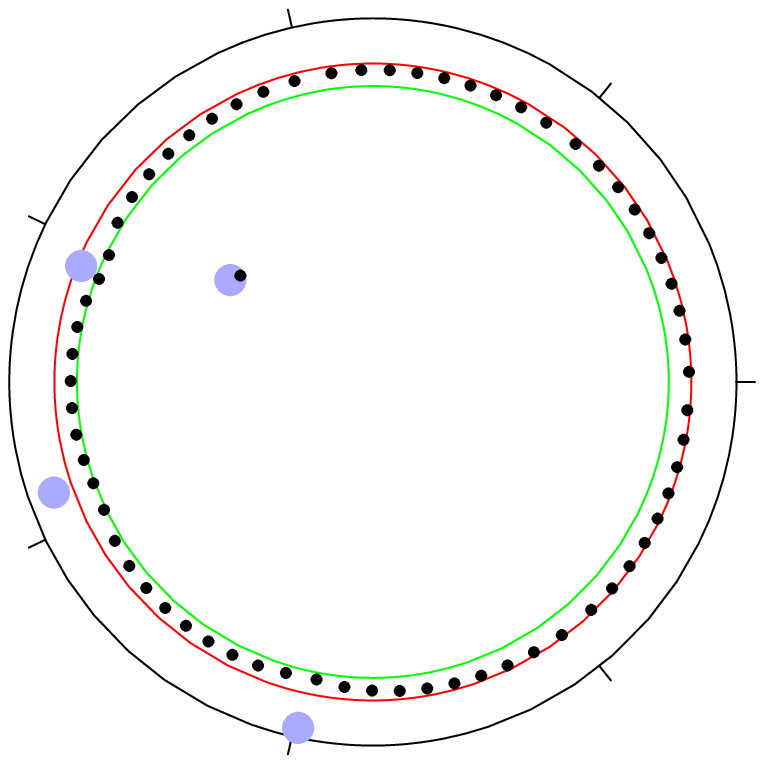,width=1.5 in}
    \psfig{file=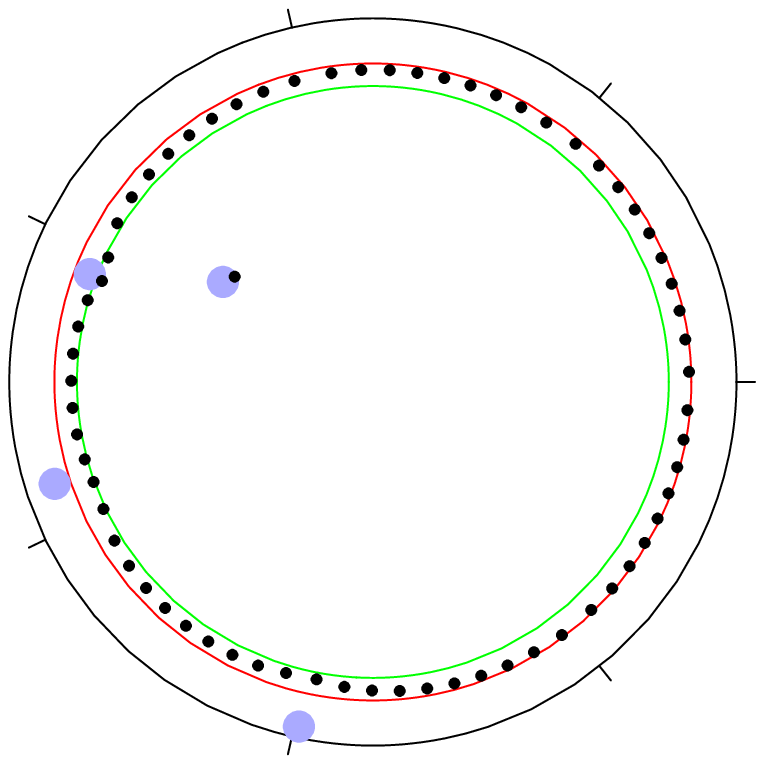,width=1.5 in}
    \psfig{file=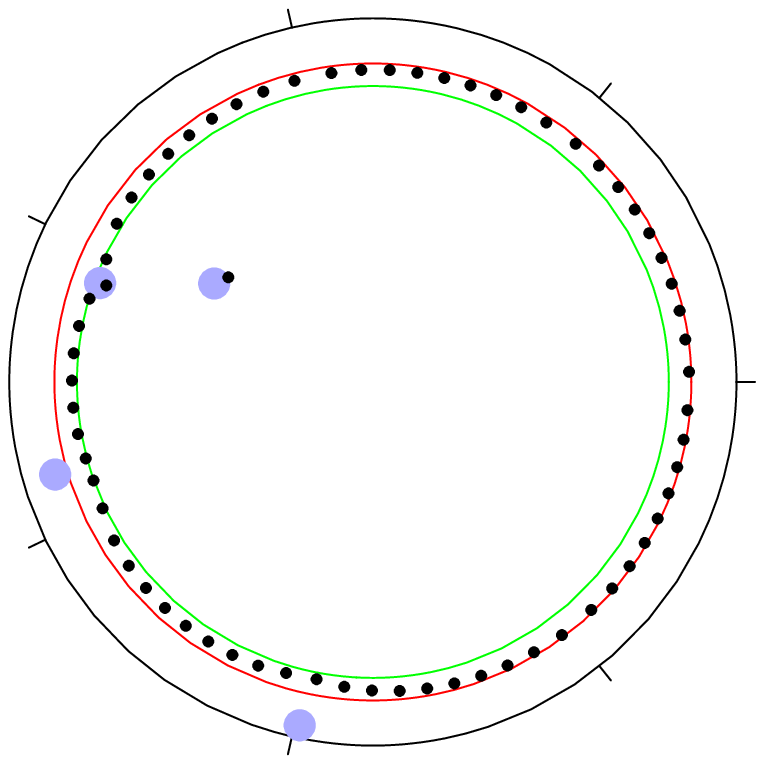,width=1.5 in} 
    \psfig{file=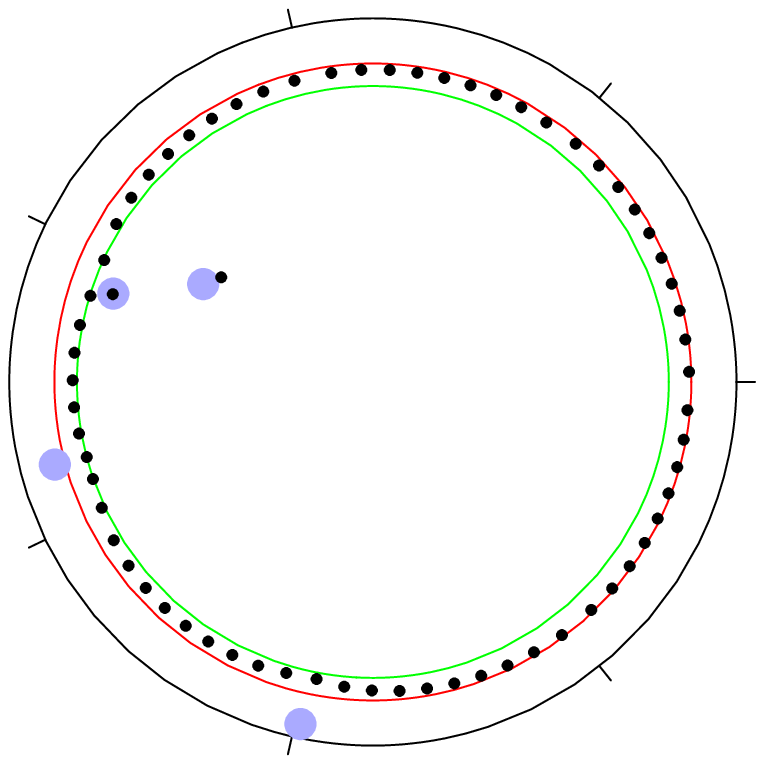,width=1.5 in} 
}\\
\mbox{
    \psfig{file=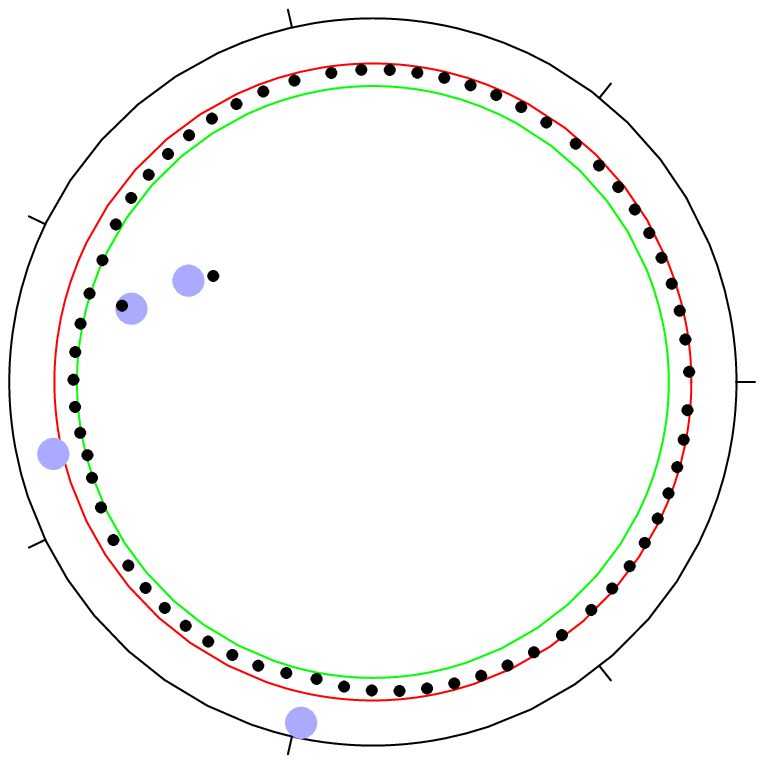,width=1.5 in}
    \psfig{file=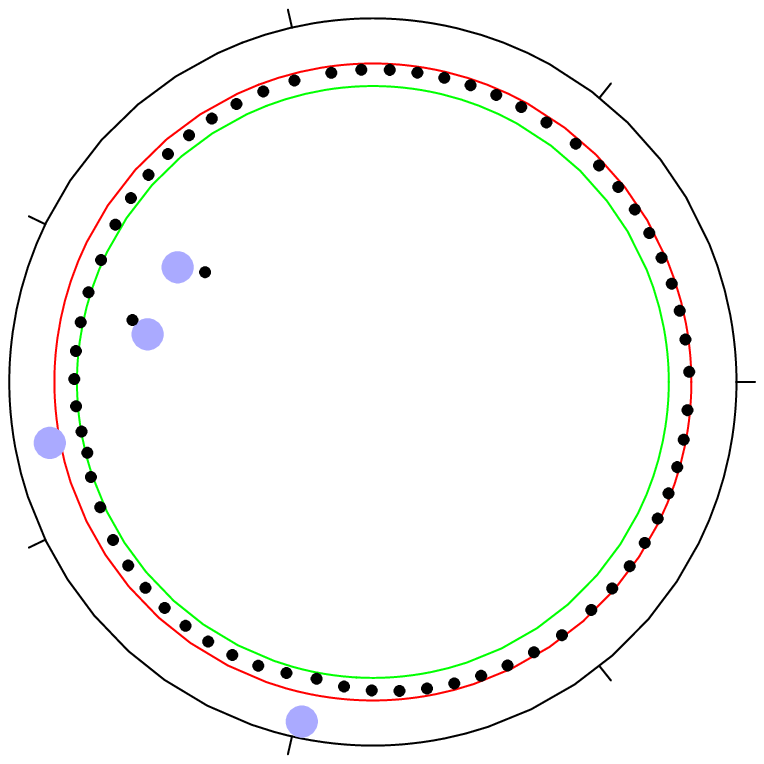,width=1.5 in}
    \psfig{file=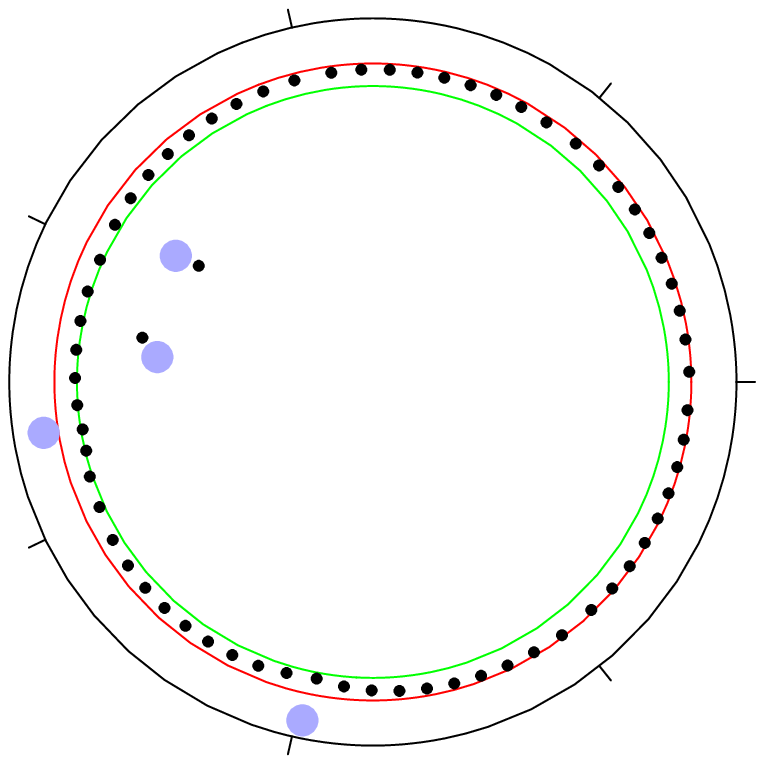,width=1.5 in}
    \psfig{file=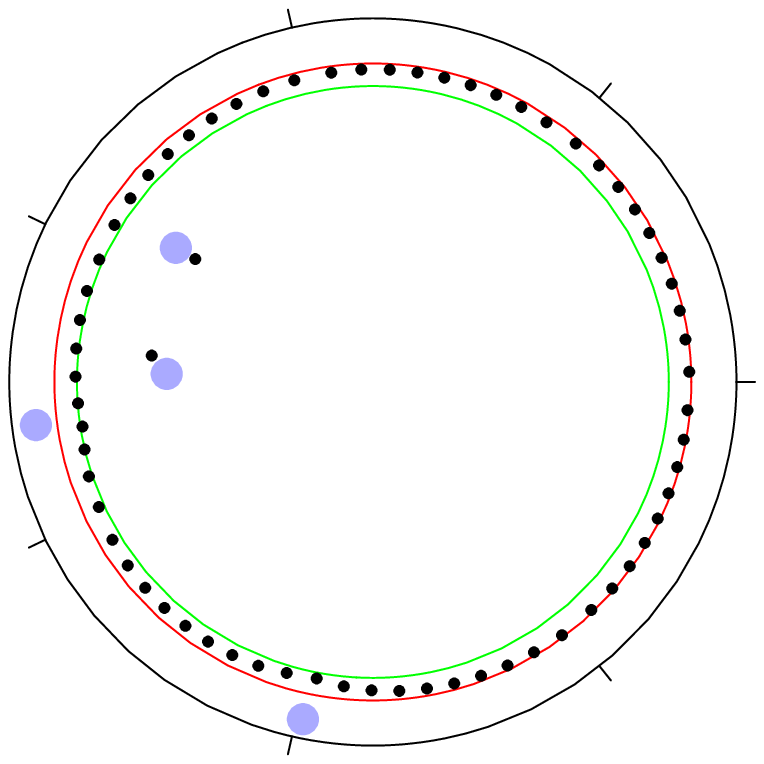,width=1.5 in}
}
\caption{\em A continuous one-parameter deformation of a weight of the form
  \eqref{eq:jumpweight} with $\ell=7$ and $k=2$.  The vector of
  parameters is ${\bf
    w}(t)^T=(0,0,0,1/2,0,1/4,0)-(1,1/2,1/4,0,1/4,0,1/2)t$ and $t$
  varies from $t=1/40$ to $t=13/80$ in steps of $\Delta t= 1/80$.  The
  frames are ordered left to right and top to bottom.  }
\label{fig:deformation}
\end{figure}
Here it can be clearly seen that as a zero of $f_n(z)$ enters the unit
disk, it initially repels the zeros near the attracting circle of
$|z|=1-(k+1)\log(n)/n$ by pushing them inwards.  The zeros along the
attracting curve then move apart to make way for the incoming zero of
$f_n(z)$.  Exactly one zero of $\pi_n(z)$ fails to get out of the way,
however, and instead enters the orbit of the moving zero of $f_n(z)$.
As the zero of $f_n(z)$ moves inside the attracting
circle, it thus draws with it a spurious zero of $\pi_n(z)$.

\subsection{The $\dbar$ steepest descent method for fixed weights.}
\label{sec:dbarmethodfixed}
Here, we begin the task of proving the theorems stated in
\S~\ref{sec:strongasymptoticspi} by analyzing the behavior of the
matrix ${\bf M}^n(z)$ solving Riemann-Hilbert Problem~\ref{rhp:M} for
a fixed weight $\phi$, in the limit $n\rightarrow\infty$.  
We recall the representation \eqref{eq:phiV} of $\phi(\theta)$ in terms
of $V:S^1\rightarrow\mathbb{R}$.  In force,
in order that Riemann-Hilbert Problem~\ref{rhp:M} indeed describes the
orthogonal polynomials with respect to $\phi(\theta)$, is the
following.
\begin{assumption}
$V$ is a real continuous function on the circle
that, for some exponent $\nu\in (0,1]$ and for some constant $K>0$,
satisfies a uniform H\"older continuity condition
$|V(\theta_2)-V(\theta_1)|\le K|\theta_2-\theta_1|^\nu$\,.
\label{ass:positiveHoelder}
\end{assumption}
This guarantees that $\phi(\theta)$ is a strictly positive function
that also satisfies a H\"older continuity condition with the same exponent,
but with a possibly different constant $K$.
\subsubsection
{Conversion to an equivalent $\dbar$ problem.  Solution of the $\dbar$ problem in terms of integral equations.}
\label{sec:strongconvert}
We proceed in several steps.  First, let ${\bf N}^n(z)$ be a new
unknown related to ${\bf M}^n(z)$ as follows:
\begin{equation}
{\bf N}^n(z):=\left\{
\begin{array}{ll}
\displaystyle {\bf M}^n(z)\left(\begin{array}{cc}z^{-n} & 0 \\ \\
0 & z^n\end{array}
\right)\,,&\hspace{0.2 in}\text{for $|z|>1$}\\\\
{\bf M}^n(z)\,,&\hspace{0.2 in}\text{for $|z|<1$}\,.
\end{array}\right.
\label{eq:Ndef}
\end{equation}
It follows from Riemann-Hilbert Problem~\ref{rhp:M} that the new
unknown ${\bf N}^n(z)$ tends to the identity matrix as $z\rightarrow\infty$,
and that ${\bf N}^n(z)$ is analytic for $|z|\neq 1$, with boundary values on
the unit circle related by
\begin{equation}
{\bf N}^n_+(e^{i\theta})={\bf N}^n_-(e^{i\theta})\left(\begin{array}{cc} e^{in\theta} & \phi(\theta)\\\\
0 & e^{-in\theta}
\end{array}\right)\,,
\label{eq:Njump}
\end{equation}
where ${\bf N}^n_+(z)$ (respectively ${\bf N}^n_-(z)$) indicates the
boundary value taken at the point $z$ on the circle from the inside
(respectively outside).

\begin{figure}[h]
\begin{center}
\input{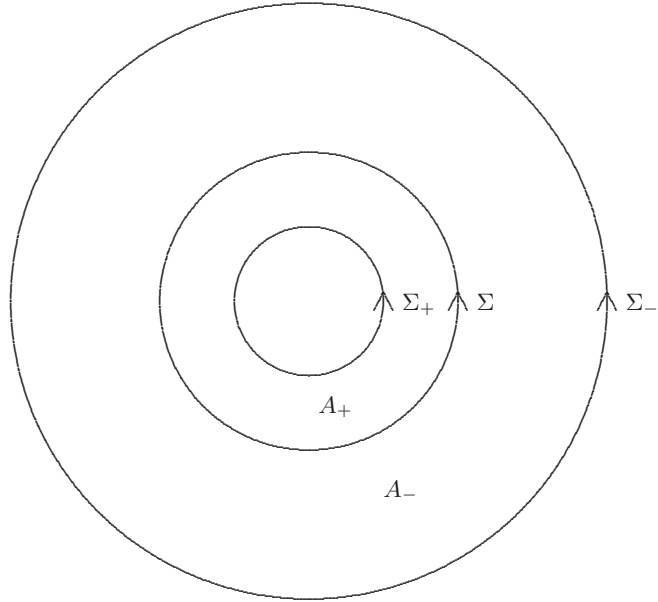}
\end{center}
\caption{\em The annular domains $A_\pm$ of the equivalent $\dbar$ problem for
  polynomials orthogonal on the unit circle. For a given $\epsilon>0$,
  the contour $\Sigma_+$ corresponds to $|z|=2^{-\epsilon}$ and the
  contour $\Sigma_-$ corresponds to $|z|=2^\epsilon$.}
\label{fig:domains}
\end{figure}

Next, observe the following factorization of the jump condition
\eqref{eq:Njump}:
\begin{equation}
{\bf N}^n_+(e^{i\theta})={\bf N}^n_-(e^{i\theta})\left(\begin{array}{cc}
1 & 0 \\ \\
e^{-in\theta}\phi(\theta)^{-1} & 1\end{array}\right)
\left(\begin{array}{cc} 0 & \phi(\theta) \\ \\
-\phi(\theta)^{-1} & 0 \end{array}\right)
\left(\begin{array}{cc} 1 & 0 \\ \\
e^{in\theta}\phi(\theta)^{-1} & 1\end{array}\right)\,.
\label{eq:factorizationstrong}
\end{equation}
To take advantage of this factorization, we introduce two new contours 
$\Sigma_\pm$ which together with $\Sigma$ bound two concentric annular domains 
$A_\pm$ as
shown in Figure~\ref{fig:domains}.

We will now need some extension of the function
$\phi(\theta)^{-1}$ defined for $z$ on the unit circle with
$\arg(z)=\theta$ to the annular domain $\overline{A_+ \cup
A_-}$.   To make use of the family of extensions defined in 
\eqref{eq:extensiondefine}, we now make
the following assumption about the weight $\phi(\theta)$.
\begin{assumption}
The function $V$ is of class $C^{k-1}(S^1)$ for some $k=1,2,3,\dots$.  
\label{ass:Ckm1}
\end{assumption}
Note that when $k=1$ this assumption is contained in
Assumption~\ref{ass:positiveHoelder}, but when $k>1$ it provides new
information.  Recall the ``bump'' function $B$ with the properties
listed in \S~\ref{sec:notation}.  Then, for any integer $m$ in the
range $1\le m \le k$, and for any $\epsilon>0$, we may apply the
extension operator $E_m$ to the function $V$ and therefore define a matrix
${\bf P}^n_{m,\epsilon}(r,\theta)$ as follows.
\begin{equation}
{\bf P}^n_{m,\epsilon}(r,\theta):=
\left\{\begin{array}{ll}
\displaystyle
{\bf N}^n(z)\left(\begin{array}{cc}
1 & 0 \\\\
z^{-n} B(\log(r)/\epsilon)e^{E_mV(r,\theta)} & 1
\end{array}\right)\,, &\hspace{0.2 in}\text{for $z=re^{i\theta}\in A_-$}\,,
\\\\
\displaystyle
{\bf N}^n(z)\left(\begin{array}{cc}
1 & 0 \\\\
-z^n B(\log(r)/\epsilon)e^{E_mV(r,\theta)} & 1
\end{array}\right)\,,&\hspace{0.2 in}
\text{for $z=re^{i\theta}\in A_+$}\,,
\\\\
\displaystyle{\bf N}^n(z)\,,&\hspace{0.2 in}
\text{for $z=re^{i\theta}\not\in\overline{A_+\cup A_-}$}\,.
\end{array}\right.
\label{eq:PN}
\end{equation}
Thus, the factor $e^{E_mV(r,\theta)}$ appearing above is our selected
extension of the function $\phi(\theta)^{-1}$ from the unit circle to
the regions $A_+$ and $A_-$.

Unlike ${\bf M}^n(z)$ and hence ${\bf N}^n(z)$, the matrix ${\bf
P}^n_{m,\epsilon}(r,\theta)$ is not piecewise analytic because the
factors relating ${\bf P}^n_{m,\epsilon}(r,\theta)$ to ${\bf N}^n(z)$
in the domains $A_\pm$ are not analytic. Indeed, in view of
\eqref{eq:dbarTaylor}, the exponent $E_mV(r,\theta)$ is not an
analytic function.  Note however that it follows from
Assumption~\ref{ass:Ckm1} and the analyticity of ${\bf M}^n(z)$ and
hence of ${\bf N}^n(z)$ for $|z|\neq 1$ that the matrix ${\bf
P}^n_{m,\epsilon}(r,\theta)$ is continuous for $r\neq 1$ as long as
$1\le m\le k$.  In particular, the ``bump'' function factor
$B(\log(r)/\epsilon)$ ensures that ${\bf P}^n_{m,\epsilon}(r,\theta)$
is continuous across the circles $\Sigma_\pm$.  At the circle of
discontinuity, $\Sigma$, the boundary values taken by ${\bf
P}^n_{m,\epsilon}(r,\theta)$ satisfy the jump condition
\begin{equation}
\lim_{r\uparrow 1}{\bf P}^n_{m,\epsilon}(r,\theta)=
\lim_{r\downarrow 1}{\bf P}^n_{m,\epsilon}(r,\theta)
\left(\begin{array}{cc} 0 & \phi(\theta)\\ \\
-\phi(\theta)^{-1}
 & 0\end{array}\right)\,.
\end{equation}

\begin{remark}
Note that the approach to this problem taken in \cite{Deift}, 
where analyticity of $\phi$ is assumed, amounts to
replacing $E_mV(r,\theta)$ with the analytic
extension $E_\infty V(r,\theta)$ and omitting the ``bump''
function factor $B(\log(r)/\epsilon)$, with the latter being at the cost of an exponentially near-identity
jump discontinuity across the inner and outer circles $\Sigma_\pm$.  \end{remark}

Next, we may remove the jump discontinuity along the unit circle by
introducing a model matrix $\dot{\bf P}(z)$ that is analytic for
$|z|\neq 1$, tends to the identity matrix $\mathbb{I}$ as
$z\rightarrow\infty$, and that takes on the unit circle continuous
boundary values $\dot{\bf P}_+(z)$ (respectively $\dot{\bf P}_-(z)$)
from the inside (respectively outside) that are related by 
\begin{equation}
\dot{\bf P}_+(e^{i\theta})=\dot{\bf P}_-(e^{i\theta})
\left(\begin{array}{cc} 0 & \phi(\theta)\\\\
-\phi(\theta)^{-1} & 0
\end{array}\right)\,.
\end{equation}
Such a matrix can be found in closed form.  Setting
\begin{equation}
\dot{\bf Q}(z)=
\dot{\bf P}(z)\left(\begin{array}{cc} 0 & -1 \\ \\
1 & 0\end{array}\right)
\end{equation}
for $|z|<1$ and $\dot{\bf Q}(z)=\dot{\bf P}(z)$ for $|z|>1$, one may
equivalently seek a matrix $\dot{\bf Q}(z)$ that is analytic for
$|z|\neq 1$, tends to the identity matrix as $z\rightarrow\infty$, and
that takes on the unit circle continuous boundary values $\dot{\bf
Q}_+(z)$ (respectively $\dot{\bf Q}_-(z)$) from the inside
(respectively outside) that are related by the diagonal jump condition
\begin{equation}
\dot{\bf Q}_+(e^{i\theta})=\dot{\bf Q}_-(e^{i\theta})
\phi(\theta)^{\sigma_3}\,,
\end{equation}
where $\sigma_3$ denotes the Pauli matrix
\begin{equation}
\sigma_3:=\left(\begin{array}{cc} 1 & 0 \\\\ 0 & -1\end{array}\right)\,.
\end{equation}
Clearly, we may seek $\dot{\bf Q}(z)$ as a diagonal matrix.
Assumption~\ref{ass:positiveHoelder} guarantees that
$\log(\phi(\theta))=-V(\theta)$ is a well-defined on the circle that
satisfies a uniform H\"older continuity condition with exponent
$\nu\in (0,1]$, and therefore we may obtain a matrix $\dot{\bf Q}(z)$
with the aforementioned properties in the explicit form $\dot{\bf
Q}(z)=S_\phi(z)^{\sigma_3}$ where $S_\phi(z)$ is the Szeg\H{o}
function associated with the weight $\phi(\theta)$ as defined in
\eqref{eq:szegofunc}.  Going back to $\dot{\bf P}(z)$, the model
matrix we will use to remove the jump discontinuity in ${\bf
P}^n_{m,\epsilon}(r,\theta)$ for $r=1$ is defined by the explicit
formula:
\begin{equation}
\dot{\bf P}(z):=\left\{\begin{array}{ll}
\displaystyle\left(\begin{array}{cc} S_\phi(z) & 0 \\ \\
0 & S_\phi(z)^{-1}\end{array}\right)\,, &
|z|>1\\\\
\displaystyle\left(\begin{array}{cc} 0 & S_\phi(z) \\ \\
-S_\phi(z)^{-1} & 0\end{array}\right)\,, &
|z|<1\,.
\end{array}\right.
\label{eq:Pdot}
\end{equation}

To actually remove the discontinuity, we introduce a new matrix function
${\bf H}^n_{m,\epsilon}(r,\theta)$ defined for $r\neq 1$ by the formula
\begin{equation}
{\bf H}^n_{m,\epsilon}(r,\theta):=
{\bf P}^n_{m,\epsilon}(r,\theta)\dot{\bf P}(z)^{-1}\,.
\label{eq:Edef}
\end{equation}
By Assumption~\ref{ass:positiveHoelder} and Assumption~\ref{ass:Ckm1},
the matrix ${\bf H}^n_{m,\epsilon}(r,\theta)$ is continuous throughout 
the two
regions $r<1$ and $r>1$.  Moreover, a continuous extension to $r=1$ is
possible because ${\bf P}^n_{m,\epsilon}(r,\theta)$ and $\dot{\bf P}(z)$
satisfy the same jump condition at $r=|z|=1$.  Thus, we see that for
$1\le m\le k$, and for any $\epsilon>0$, the matrix function ${\bf H}^n_{m,\epsilon}(r,\theta)$ defined
by \eqref{eq:Edef} may be
viewed as a continuous function on the whole plane with polar
coordinates $-\pi\le\theta<\pi$ and $0\le r<\infty$.

At this point, we can summarize the explicit transformations we have
introduced and relate ${\bf H}^n_{m,\epsilon}(r,\theta)$ directly back to
${\bf M}^n(z)$.  Combining \eqref{eq:Ndef}, \eqref{eq:PN}, \eqref{eq:Pdot},
and \eqref{eq:Edef}, we have by definition
\begin{equation}
{\bf H}^n_{m,\epsilon}(r,\theta):=\left\{
\begin{array}{ll}
\displaystyle {\bf M}^n(z)\left(\begin{array}{cc}
0 & -S_\phi(z)\\\\S_\phi(z)^{-1} & z^nS_\phi(z)B(\log(r)/\epsilon)
e^{E_mV(r,\theta)}
\end{array}\right)\,, &\hspace{0.2 in}0\le r<1\\\\
\displaystyle{\bf M}^n(z)\left(\begin{array}{cc}
z^{-n}S_\phi(z)^{-1} & 0 \\\\
S_\phi(z)^{-1}B(\log(r)/\epsilon)e^{E_mV(r,\theta)} & z^nS_\phi(z)\end{array}
\right)\,,&\hspace{0.2 in}r>1
\end{array}\right.
\label{eq:EM}
\end{equation}
where $z=re^{i\theta}$.

For $r<2^{-\epsilon}$ and $r>2^\epsilon$ we have 
$B(\log(r)/\epsilon)\equiv 0$, and in
these regions the matrix ${\bf H}^n_{m,\epsilon}(r,\theta)$ clearly inherits
analyticity from ${\bf M}^n(z)$; in other words, in these regions
${\bf H}^n_{m,\epsilon}(r,\theta)$ is a smooth function of the combination
$z=re^{i\theta}$.  However, for $z\in\Omega_\pm$, the matrix ${\bf
H}^n_{m,\epsilon}(r,\theta)$ is certainly not analytic.  In order to measure
the deviation from analyticity in the case when $m=k$, we introduce a
further assumption on $\phi(\theta)$.
\begin{assumption}
The function $V$ is of class $C^{k-1,1}(S^1)$, that is, the
function $V^{(k-1)}(\theta)$ is Lipschitz continuous.
\label{ass:ps}
\end{assumption}
Note that since $k\ge 1$, this condition implies in particular that
$V(\theta)$ is Lipschitz, and thus the H\"older
continuity part of Assumption~\ref{ass:positiveHoelder} is subsumed.
With Assumption~\ref{ass:ps} in force, we may compute the
$\dbar$-derivative of ${\bf H}^n_{m,\epsilon}(r,\theta)$ for $r\neq 1$ and for
all integer $m$ in the range $1\le m \le k$.  
Differentiation of \eqref{eq:EM} yields
\begin{equation}
\dbar{\bf H}^n_{m,\epsilon}(r,\theta)=\left\{\begin{array}{ll}
\displaystyle{\bf M}^n(z)\left(\begin{array}{cc}
0 & 0 \\\\
0 & z^nS_\phi(z)\dbar\left[B(\log(r)/\epsilon)e^{E_mV(r,\theta)}\right]
\end{array}\right)\,,&\hspace{0.2 in}0\le r<1\\\\
\displaystyle{\bf M}^n(z)\left(\begin{array}{cc}
0 & 0 \\\\
S_\phi(z)^{-1}\dbar\left[B(\log(r)/\epsilon)e^{E_mV(r,\theta)}\right] 
& 0\end{array}
\right)\,, &\hspace{0.2 in}r>1\,,
\end{array}\right.
\label{eq:dbarEM}
\end{equation}
for almost all $\theta$, and then elimination of ${\bf M}^n(z)$ in
terms of ${\bf H}^n_{m,\epsilon}(r,\theta)$ using \eqref{eq:EM} again gives
\begin{equation}
\dbar{\bf H}^n_{m,\epsilon}(r,\theta)={\bf H}^n_{m,\epsilon}(r,\theta){\bf W}^n_{m,\epsilon}(r,\theta)\,,
\hspace{0.2 in}\text{for $r\neq 1$ and almost all $\theta\in S^1$}\,,
\label{eq:dbarEE}
\end{equation}
where the matrix ${\bf W}^n_{m,\epsilon}(r,\theta)$ is given for $r\neq 1$ and
almost all $\theta$ by
the explicit formula
\begin{equation}
{\bf W}^n_{m,\epsilon}(r,\theta):=\left\{\begin{array}{ll}
\displaystyle
\left(\begin{array}{cc} 0 & z^nS_\phi(z)^{2}\dbar\left[B(\log(r)/\epsilon)
e^{E_mV(r,\theta)}\right]\\\\ 0 & 0\end{array}\right)\,, &\hspace{0.2 in}
0\le r <1\\\\
\displaystyle
\left(\begin{array}{cc} 0 & 0 \\\\
z^{-n}S_\phi(z)^{-2}\dbar\left[B(\log(r)/\epsilon)e^{E_mV(r,\theta)}\right] 
& 0
\end{array}\right)\,,&\hspace{0.2 in}r>1\,.
\end{array}\right.
\label{eq:Wdef}
\end{equation}
In particular, we see that the matrix ${\bf H}^n_{m,\epsilon}(r,\theta)$ is a solution
of the following $\dbar$-problem.
\begin{dbp}
Find a $2\times 2$ matrix ${\bf U}(r,\theta)$ with the properties:
\begin{itemize}
\item[]{\bf Smoothness.}  ${\bf U}(r,\theta)$ is a Lipschitz continuous
function throughout $\mathbb{R}^2$.  
\item[]{\bf Deviation From Analyticity.}
The relation 
\begin{equation}
\dbar{\bf U}(r,\theta)={\bf U}(r,\theta){\bf W}^n_{m,\epsilon}(r,\theta)
\end{equation}
holds for all points in $\mathbb{R}^2$ with the exception of a set of
measure zero.  The matrix ${\bf W}^n_{m,\epsilon}(r,\theta)$ is defined almost
everywhere by
\eqref{eq:Wdef} and is essentially compactly supported.
\item[]{\bf Normalization.}
The matrix ${\bf U}(r,\theta)$ is normalized at $r=\infty$ as follows:
\begin{equation}
\lim_{r\rightarrow\infty}{\bf U}(r,\theta)={\mathbb I}\,.
\end{equation}
\end{itemize}
\label{dbp:E}
\end{dbp}
In writing down this $\dbar$-problem, we have focused on just a few
specific properties of the matrix ${\bf H}^n_{m,\epsilon}(r,\theta)$.  However, it
is important that in doing so, we have not introduced any spurious solutions.
\begin{prop}
Suppose that $\phi=e^{-V}$ where $V:S^1\rightarrow\mathbb{R}$ is of
class $C^{k-1,1}(S^1)$ for some $k=1,2,3,\dots$. Then for all
$n=0,1,2,3,\dots$, for $m=1,2,\dots,k$, and for all $\epsilon>0$, the
matrix ${\bf W}_{m,\epsilon}^n(r,\theta)$ is well-defined almost
everywhere by
\eqref{eq:Wdef} and $\dbar$ Problem~\ref{dbp:E} has a unique solution,
namely ${\bf U}(r,\theta)={\bf H}^n_{m,\epsilon}(r,\theta)$.
\label{prop:dbarEsolve}
\end{prop}
\begin{proof}
The existence of a solution follows from \eqref{eq:EM} and the
existence of ${\bf M}^n(z)$ for $n=0,1,2,\dots$.  To establish the
uniqueness we first consider the determinant of any solution of
$\dbar$ Problem~\ref{dbp:E}.  Clearly, $\det({\bf U}(r,\theta))$ is a
Lipschitz continuous function that tends to $1$ as
$z\rightarrow\infty$.  Moreover, the relation $\dbar\det({\bf
U}(r,\theta)) = {\rm tr}({\bf W}^n_{m,\epsilon}(r,\theta))\det({\bf
U}(r,\theta))$ holds almost everywhere, and thus by
\eqref{eq:Wdef} we see that $\dbar\det({\bf U}(r,\theta))=0$ holds
almost everywhere in the plane.  It follows that $\det({\bf
U}(r,\theta))$ is not only Lipschitz continuous, but is in fact an
entire function of $z=re^{i\theta}$ that tends to $1$ as
$z\rightarrow\infty$. Therefore from Liouville's Theorem we see that
$\det({\bf U}(r,\theta))\equiv 1$.  Next, consider the matrix ratio of
any two solutions ${\bf U}(r,\theta)$ and $\tilde{\bf U}(r,\theta)$ of
$\dbar$ Problem~\ref{dbp:E}; this is the matrix ${\bf R}(r,\theta)$
defined by
\begin{equation}
{\bf R}(r,\theta):={\bf U}(r,\theta)
\tilde{\bf U}(r,\theta)^{-1}\,.
\end{equation}
Since $\det(\tilde{\bf U}(r,\theta))\equiv 1$, it follows that ${\bf
R}(r,\theta)$ is Lipschitz continuous throughout the plane.
By direct calculation, 
we have
\begin{equation}
\begin{array}{rcl}
\displaystyle
\dbar{\bf R}(r,\theta)&=&\displaystyle
\dbar{\bf U}(r,\theta)\cdot
\tilde{\bf U}(r,\theta)^{-1} - 
{\bf U}(r,\theta)\tilde{\bf U}(r,\theta)^{-1}
\dbar\tilde{\bf U}(r,\theta)\cdot\tilde{\bf U}(r,\theta)^{-1}\\\\
&=&\displaystyle
{\bf U}(r,\theta){\bf W}^n_{m,\epsilon}(r,\theta)\tilde{\bf U}(r,\theta)^{-1}
-
{\bf U}(r,\theta){\bf W}^n_{m,\epsilon}(r,\theta)\tilde{\bf U}(r,\theta)^{-1}
\\\\
&=&{\bf 0}\,,
\end{array}
\end{equation}
holding almost everywhere in the plane.
It follows that ${\bf R}(r,\theta)$ is an entire function of
$z=re^{i\theta}$ that tends to the identity matrix as
$z\rightarrow\infty$, so again by Liouville's Theorem we get ${\bf
R}(r,\theta)\equiv \mathbb{I}$, or equivalently $\tilde{\bf
U}(r,\theta)\equiv{\bf U}(r,\theta)$.
\end{proof}
The unique solution of $\dbar$ Problem~\ref{dbp:E} can also be expressed as
a solution of an integral equation with Cauchy kernel.
\begin{prop}
Suppose that $\phi=e^{-V}$ where $V:S^1\rightarrow\mathbb{R}$ is of
class $C^{k-1,1}(S^1)$ for some $k=1,2,3,\dots$.  Then for all
$n=0,1,2,3,\dots$, for $m=1,2,\dots,k$, and for all $\epsilon>0$, the
matrix ${\bf W}^n_{m,\epsilon}(r,\theta)$ is well-defined almost
everywhere by
\eqref{eq:Wdef} and the corresponding solution ${\bf U}(r,\theta)={\bf
H}^n_{m,\epsilon}(r,\theta)$ of $\dbar$ Problem~\ref{dbp:E} satisfies
the integral equation
\begin{equation}
{\bf U}(r,\theta)=\mathbb{I}-\frac{1}{\pi}\int\int \frac{{\bf
U}(r',\theta'){\bf W}^n_{m,\epsilon}(r',\theta')}{z'-z}\,dA'
\label{eq:inteqnE}
\end{equation}
where $z=re^{i\theta}$, $z'=r'e^{i\theta'}$, and $dA'$ is a positive
area element $dA'=r'\,dr'\,d\theta'$.  The integral is taken over
the entire plane.
\label{prop:dbarEinteqn}
\end{prop}
\begin{proof}
Recall that the Cauchy kernel is a fundamental solution
for the $\dbar$ operator.  In the relation \eqref{eq:dbarEE} we may replace
$\dbar{\bf H}^n_{m,\epsilon}(r,\theta)$ by $\dbar[{\bf H}^n_{m,\epsilon}(r,\theta)-\mathbb{I}]$;
multiplying by the Cauchy kernel and integrating over the whole plane gives the
identity
\begin{equation}
-\frac{1}{\pi}\int\int\frac{\dbar[{\bf
H}^n_{m,\epsilon}(r',\theta')-\mathbb{I}]}{z'-z}\,dA' = -\frac{1}{\pi}\int\int
\frac{{\bf H}^n_{m,\epsilon}(r',\theta'){\bf W}^n_{m,\epsilon}(r',\theta')}{z'-z}\,dA'\,.
\end{equation}
On the left-hand side the $\dbar$ operator differentiates with respect
to the primed variables.  Since the Cauchy kernel is absolutely
integrable, we may evaluate the integral on the left-hand side by
replacing the domain of integration by the region $|z'-z|\ge\delta>0$
and subsequently taking the limit $\delta\rightarrow 0$.  For each
positive $\delta$ we may apply Stokes' Theorem and use the facts that
$\dbar[(z'-z)^{-1}]=0$ for $|z'-z|\ge \delta$ and that ${\bf
  H}^n_{m,\epsilon}(r,\theta)-\mathbb{I}$ tends to zero as
$r\rightarrow\infty$ to evaluate the integral over the region
$|z'-z|\ge\delta$ in terms of a line integral over the boundary.  Thus
we have
\begin{equation}
\lim_{\delta\downarrow 0}
\frac{1}{2\pi\delta}
\int_{|z'-z|=\delta}({\bf H}^n_{m,\epsilon}(r',\theta')-\mathbb{I})\,
d\ell' = -\frac{1}{\pi}\int\int
\frac{{\bf H}^n_{m,\epsilon}(r',\theta'){\bf W}^n_{m,\epsilon}(r',\theta')}{z'-z}\,dA'\,,
\end{equation}
where $d\ell'$ is an arc length element.  From the continuity of ${\bf
H}^n_{m,\epsilon}(r,\theta)$ the integral equation
\eqref{eq:inteqnE} with 
${\bf U}(r,\theta)={\bf H}^n_{m,\epsilon}(r,\theta)$ follows.
\end{proof}

\subsubsection
{Asymptotic solution of the integral equation.  Estimates of ${\bf
H}^n_{m,\epsilon}(r,\theta)$ and its derivatives for large $n$.}  
\label{sec:stronginteqn}
Being as knowledge of the matrix ${\bf
H}^n_{m,\epsilon}(r,\theta)$ is equivalent to 
knowledge of
${\bf M}^n(z)$ and hence of the polynomial of degree $n$ in the system
of polynomials orthogonal on the circle with respect to $\phi$, we
would like to use the integral equation
\eqref{eq:inteqnE} to characterize ${\bf H}^n_{m,\epsilon}(r,\theta)$.  
There is a difficulty in that, while existence of solutions for 
\eqref{eq:inteqnE} is not an issue, one does not automatically have 
uniqueness.  However, it turns out that if the parameter $n$ is
sufficiently large, then the integral equation \eqref{eq:inteqnE}
defines a contraction mapping and thus may be solved by iteration
yielding a unique solution in the form of a Neumann series.  In this
connection, we can also obtain from \eqref{eq:inteqnE} asymptotic
information about the matrix ${\bf H}^n_{m,\epsilon}(r,\theta)$, and
consequently of the orthogonal polynomial $\pi_n(z)$, in the limit
$n\rightarrow\infty$.

In order to study \eqref{eq:inteqnE}, it is useful to characterize the
family of matrix functions ${\bf W}^n_{m,\epsilon}(r,\theta)$ more
concretely.
\begin{prop}
Suppose that $V:S^1\rightarrow\mathbb{R}$ is a real function of class
$C^{k-1,1}(S^1)$ for some $k\ge 1$, that $m$ is an integer satisfying
$1\le m\le k$, and that $\epsilon>0$ is fixed.  Let the integer $D$ be
defined as $D:=\min(k-m,m-1)$.  Then, the matrix function ${\bf
W}^n_{m,\epsilon}$ is of class
$C_0^{D-1,1}(\mathbb{R}^2\setminus\{0\})$ if $D>0$, and of class
$L_0^\infty(\mathbb{R}^2\setminus\{0\})$ if $D=0$.  Moreover, if
$\alpha$ and $\beta$ are nonnegative integers such that
$\alpha+\beta\le D$, then there is a constant
$C_{m,\epsilon}^{(\alpha,\beta)}>0$
such that for all $n$ the estimate
\begin{equation}
\left\|\frac{\partial^{\alpha+\beta}}{\partial r^\alpha\partial\theta^\beta}
{\bf W}^n_{m,\epsilon}(r,\theta)\right\|\le
C^{(\alpha,\beta)}_{m,\epsilon}n^\beta
e^{-n|\log(r)|}|\log(r)|^{m-1-\alpha}\sum_{p=0}^{\alpha}n^p|\log(r)|^p
\label{eq:derivWbound}
\end{equation}
holds throughout the region $|\log(r)|\le \epsilon\log(2)$ containing
the essential support of ${\bf W}^n_{m,\epsilon}(r,\theta)$.  
\label{prop:Wbound}
\end{prop}
\begin{proof}
${\bf W}_{m,\epsilon}^n(r,\theta)$ vanishes identically outside of the
annulus $|\log(r)|\le \epsilon\log(2)$.  In the disjoint regions $r<1$
and $r>1$, the matrix function ${\bf W}_{m,\epsilon}^n(r,\theta)$ is
infinitely differentiable with respect to $r$, and the issue is the
continuity of these derivatives at $r=1$.  The relation
\eqref{eq:dbarTaylor} implies that for each fixed $\theta$, ${\bf
W}^n_{m,\epsilon}(r,\theta)$ is proportional to $(\log(r))^{m-1}$ near
$r=1$ (where $B(\log(r)/\epsilon)\equiv 1$ holds), and thus all
derivatives of ${\bf W}^n_{m,\epsilon}(r,\theta)$ with respect to $r$
through order $m-2$ are Lipschitz continuous at $r=1$, and the derivative
$\partial^{m-1} {\bf W}^n_{m,\epsilon}(r,\theta)/\partial r^{m-1}$
remains bounded as $r\rightarrow 1$, but experiences a jump
discontinuity at $r=1$.  

On the other hand, if $r\neq 1$ is fixed, then from
\eqref{eq:extensiondefine} and \eqref{eq:dbarTaylor}, the matrix ${\bf
W}_{m,\epsilon}^n(r,\theta)$ depends analytically on derivatives
$V^{(j)}(\theta)$ for $0\le j \le m$.  Since $V$ is of class
$C^{k-1,1}(S^1)$, all derivatives of ${\bf W}^n_{m,\epsilon}$ with
respect to $\theta$ through order $k-1-m$ will be Lipschitz
continuous, while the derivative $\partial^{k-m} {\bf
W}_{m,\epsilon}^n(r,\theta)/\partial\theta^{k-m}$ will be defined for
almost all $\theta$ and will be uniformly bounded.  

To have all mixed partial derivatives of total order at most $D-1$ to
be Lipschitz continuous, it is sufficient to have both $D\le m-1$ and
$D\le k-m$.  If for $1\le m\le k$ these inequalities force $D=0$, then
no derivatives of ${\bf W}_{m,\epsilon}^n(r,\theta)$ may be taken at
all, but ${\bf W}_{m,\epsilon}^n(r,\theta)$ is uniformly bounded and
compactly supported in the annulus $|\log(r)|\le \epsilon\log(2)$,
that is, ${\bf W}_{m,\epsilon}^n\in
L^\infty_0(\mathbb{R}^2\setminus\{0\})$.  If $D=\min(k-m,m-1)>0$, then
we learn that ${\bf W}_{m,\epsilon}^n\in
C_0^{D-1,1}(\mathbb{R}^2\setminus\{0\})$.

Now there are absolute constants $C_\pm>0$ such that 
\begin{equation}
\left\|\frac{\partial^{\alpha+\beta}}{\partial r^\alpha\partial\theta^\beta}
{\bf W}_{m,\epsilon}^n(r,\theta)\right\| = \left
\{\begin{array}{ll}
\displaystyle
C_+\left|\mathop{\sum_{\alpha'+\alpha''=\alpha}}_{\beta'+\beta''=\beta}
W_{m,\epsilon,\alpha'',\beta''}(r,\theta)
\frac{\partial^{\alpha'+\beta'}}
{\partial r^{\alpha'}\partial\theta^{\beta'}}z^nS_\phi(z)^2
\right|\,,&\hspace{0.2 in} r<1\\\\
\displaystyle
C_-\left|\mathop{\sum_{\alpha'+\alpha''=\alpha}}_{\beta'+\beta''=\beta}
W_{m,\epsilon,\alpha'',\beta''}(r,\theta)
\frac{\partial^{\alpha'+\beta'}}
{\partial r^{\alpha'}\partial\theta^{\beta'}}z^{-n}S_\phi(z)^{-2}
\right|\,,&\hspace{0.2 in}r>1\,,
\end{array}\right.
\end{equation}
where
\begin{equation}
W_{m,\epsilon,\alpha'',\beta''}(r,\theta):=
\frac{\partial^{\alpha''+\beta''}}{\partial r^{\alpha''}
\partial\theta^{\beta''}}\dbar\left[B(\log(r)/\epsilon)e^{E_mV(r,\theta)}
\right]\,.
\end{equation}
Generally, the derivatives indexed by $\alpha'$ and $\beta'$ are
uniformly bounded throughout the regions
$2^{-\epsilon}<r<1$ and $1<r<2^\epsilon$ by a constant multiple of
$n^{\alpha'+\beta'}e^{-n|\log(r)|}$.  Furthermore, using
\eqref{eq:dbarTaylor}, the derivative
$W_{m,\epsilon,\alpha'',\beta''}(r,\theta)$ is uniformly bounded by a
constant multiple of $|\log(r)|^{m-1-\alpha''}$ throughout the region
$|\log(r)|\le \epsilon\log(2)$.
Therefore setting $\alpha''=\alpha-\alpha'$, an inequality of the form
\begin{equation}
\left\|\frac{\partial^{\alpha+\beta}}{\partial r^\alpha\partial\theta^\beta}
{\bf W}^n_{m,\epsilon}(r,\theta)\right\|\le 
\tilde{C}_{m,\epsilon}^{(\alpha,\beta)}e^{-n|\log(r)|}\sum_{\alpha'=0}^{\alpha}
\sum_{\beta'=0}^\beta n^{\beta'}|\log(r)|^{m-1-\alpha}n^{\alpha'}|\log(r)|^{\alpha'}
\end{equation}
holds in the region $|\log(r)|\le \epsilon\log(2)$, 
where $\tilde{C}_{m,\epsilon}^{(\alpha,\beta)}>0$ is a constant.

Finally, since $n^{\beta'}\le n^\beta$, the inequality
\eqref{eq:derivWbound} follows, where $C_{m,\epsilon}^{(\alpha,\beta)} = 
(\beta+1)\tilde{C}_{m,\epsilon}^{(\alpha,\beta)}$.
\end{proof}

An important part of our analysis will be the estimation of certain
two-dimensional Laplace-type integrals with Cauchy kernels.  The main workhorse
in this connection is the following lemma.

\begin{lemma}
Let $\epsilon>0$ and $\nu\ge 1$ be fixed constants.  Then there exists a
corresponding constant $K_{\epsilon,\nu}>0$ such that the estimate
(note that $z'=r'e^{i\theta'}$)
\begin{equation}
\sup_{z\in\mathbb{C}}
\int_{2^{-\epsilon}}^{2^\epsilon}r'\,dr'\,e^{-n|\log(r')|}
|\log(r')|^{\nu-1}\int_{-\pi}^\pi\frac{d\theta'}{|z'-z|}\le K_{\epsilon,\nu}
\frac{\log(n)}{n^\nu}
\label{eq:keybounduniform}
\end{equation}
holds for sufficiently large $n$.  Moreover, for each $\rho>2^\epsilon$, 
there exists a constant $K_{\epsilon,\nu,\rho}>0$ such that the estimate
\begin{equation}
\sup_{|\log(|z|)|\ge\log(\rho)}
\int_{2^{-\epsilon}}^{2^\epsilon}r'\,dr'\,e^{-n|\log(r')|}
|\log(r')|^{\nu-1}\int_{-\pi}^\pi\frac{d\theta'}{|z'-z|}\le 
\frac{K_{\epsilon,\nu,\rho}}{n^\nu}
\label{eq:keyboundaway}
\end{equation}
holds for sufficiently large $n$.  
\label{lem:keybound}
\end{lemma}

\begin{proof}
As $\theta'$ varies over $S^1$, the minimum value of $|z'-z|$
is achieved at $\theta'=\theta$, and we thus have $|z'-z|\ge|r'-r|$, which
implies the inequality
\begin{equation}
\int_{-\pi}^\pi\frac{d\theta'}{|z'-z|}\le\frac{2\pi}{|r'-r|}\,.  
\end{equation}
Let $\mu$ be a positive constant.  Clearly, there is another positive
constant $C_1$ depending on $\mu$ but not on $r$ or $r'$ such that
\begin{equation}
|\log(r')-\log(r)|>\mu \hspace{0.2 in}\text{implies}\hspace{0.2 in}
\int_{-\pi}^\pi\frac{d\theta'}{|z'-z|}\le C_1
\label{eq:thetaintaway}
\end{equation}
because the condition $|\log(r')-\log(r)|>\mu$ also bounds $|r'-r|$
away from zero.  Furthermore, if $\mu$ is sufficiently small, there is
a positive constant $C_2$ depending on $\mu$ but not on $r$ or $r'$
such that
\begin{equation}
\text{$|\log(r')-\log(r)|\le \mu$ and $|\log(r')|\le \epsilon\log(2)$}
\hspace{0.1 in}\text{implies}\hspace{0.1 in}
\int_{-\pi}^\pi\frac{d\theta'}{|z'-z|}\le 
C_2\log\left(n+\frac{1}{|\log(r')-\log(r)|}\right)\,.
\label{eq:thetaintclose}
\end{equation}
To establish \eqref{eq:thetaintclose}, note that for $\mu$
sufficiently small the condition $|\log(r')-\log(r)|\le \mu$ implies
that $|r'-r|<\pi/2$.  Assuming without loss of generality that
$\theta=0$, we use the following estimates of the integrand.  
For $|\theta'|\le |r'-r|$ we use the estimate $|z'-z|\ge |r'-r|$
which follows from the triangle inequality applied to the identity
$z+(z'-z)=z'$.  For $|r'-r|<|\theta'|\le\pi/2$ we use the estimate
$|z'-z|\ge r'|\sin(\theta')|\ge r|\theta'|/2$ which follows from the
diagram shown in Figure~\ref{fig:trigonometry}.  
\begin{figure}[h]
\begin{center}
\input{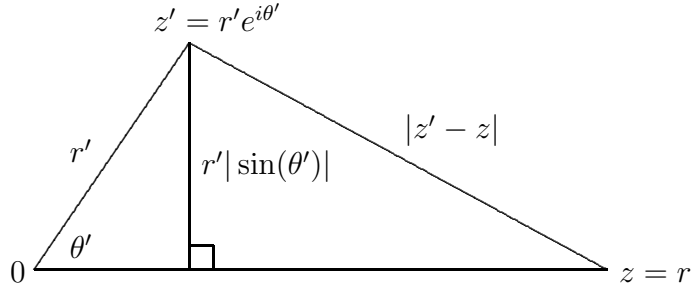}
\end{center}
\caption{\em The estimate $|z'-z|\ge r'|\sin(\theta')|$.  The right-hand
side can be replaced by $r'|\theta'|/2$ for $|\theta'|<\pi/2$.}
\label{fig:trigonometry}
\end{figure}
Finally, for $\pi/2<|\theta'|\le\pi$ we use the estimate $|z'-z|\ge r$
which follows from the Law of Cosines because $\cos(\theta')<0$:
$|z'-z|^2 = r^2 + (r')^2 -2rr'\cos(\theta')\ge r^2 + (r')^2 \ge r^2$
(again, see Figure~\ref{fig:trigonometry}).
Combining these estimates, we have
\begin{equation}
\int_{-\pi}^\pi\frac{d\theta'}{|z'-z|}\le 2 + 
\frac{4}{r'}\int_{|r'-r|}^{\pi/2}\frac{d\theta'}{\theta'} + \frac{\pi}{r}\le 
2 + 2^{2+\epsilon}\int_{|r'-r|}^{\pi/2}\frac{d\theta'}{\theta'} + 
2^\epsilon\pi e^{\mu}\,.
\end{equation}
This is clearly bounded above by a constant multiple of
$\log(|r'-r|^{-1})$ for $\mu$ sufficiently small, and from this
\eqref{eq:thetaintclose} follows as well (adding a constant inside the
logarithm keeps the bound positive away from the singularity as long
as $n\ge 1$, and for later purposes it is convenient to take the
additive constant to be $n$).

Now we estimate the integral over $r'$.  Using $r'\le 2^\epsilon$, along with
\eqref{eq:thetaintaway} and \eqref{eq:thetaintclose}, and changing the
integration variable from $r'$ to $s=n\log(r')$, we get
\begin{equation}
\begin{array}{l}\displaystyle
\int_{2^{-\epsilon}}^{2^\epsilon}r'\,dr'\,e^{-n|\log(r')|}|\log(r')|^{\nu-1}
\int_{-\pi}^\pi
\frac{d\theta'}{|z'-z|} \le 
\frac{2^{2\epsilon}C_1}{n^\nu}\int_{\mathcal{A}}
ds\,e^{-|s|}|s|^{\nu-1} \\\\\displaystyle\hspace{1.5 in}+\,\,\,
\frac{2^{2\epsilon}C_2}{n^\nu}\int_{\mathcal{B}}
ds\,e^{-|s|}|s|^{\nu-1}\log\left(n+\frac{n}{|
s-n\log(r)|}\right)\,.
\end{array}
\end{equation}
where 
\begin{equation}
\begin{array}{rcl}
\mathcal{A}&:=&\{
\text{$s$ such that $|s|\le n\epsilon\log(2)$ and $|s-n\log(r)|> n\mu$}\}\\\\
\mathcal{B}&:=&\{
\text{$s$ such that $|s|\le n\epsilon\log(2)$ and $|s-n\log(r)|\le n\mu$}\}\,.
\end{array}
\end{equation}
Finally, we have
\begin{equation}
\int_{\mathcal{A}}ds\,e^{-|s|}|s|^{\nu-1}\le\int_{-\infty}^\infty ds\,
e^{-|s|}|s|^{\nu-1}\,,
\end{equation}
which is finite and independent of $n$ and $r$ because $\nu\ge 1$, while
\begin{equation}
\begin{array}{rcl}\displaystyle
\int_{\mathcal{B}}ds\,e^{-|s|}|s|^{\nu-1}\left(n+\frac{n}{|s-n\log(r)|}\right)
&\le &\displaystyle
\int_{-\infty}^\infty ds\,e^{-|s|}|s|^{\nu-1}\log
\left(n+\frac{n}{|s-n\log(r)|}\right)\\\\&=&
\displaystyle\log(n)\int_{-\infty}^\infty ds\,e^{-|s|}|s|^{\nu-1}\\\\&&
\displaystyle\hspace{0.2 in}
+\,\,\,\int_{-\infty}^\infty ds\,e^{-|s|}|s|^{\nu-1}
\log\left(1+\frac{1}{|s-n\log(r)|}\right)\,,
\end{array}
\end{equation}
and the last integral is bounded independently of $n$ and $r$, since
by Cauchy-Schwarz,
\begin{equation}
\int_{-\infty}^\infty ds\,e^{-|s|}|s|^{\nu-1}\log\left(1+
\frac{1}{|s-n\log(r)|}\right)\le
\left[\int_{-\infty}^\infty ds\, e^{-2|s|} |s|^{2\nu-2}\right]^{1/2}
\left[\int_{-\infty}^\infty ds 
\left(\log\left(1+\frac{1}{|s|}\right)\right)^2\right]^{1/2}
\end{equation}
with both factors being finite.  Thus, an upper bound for the integral
of interest is proportional to $\log(n)/n^\nu$ in general, which
proves \eqref{eq:keybounduniform}.  If $\rho>2^\epsilon$ and
$|\log(|z|)|\ge\log(\rho)$, then it is not necessary to divide the
integration into sets $\mathcal{A}$ and $\mathcal{B}$, and the bound
\eqref{eq:thetaintaway} can be used over the whole range of
integration in which case the upper bound is then proportional to
$1/n^\nu$ which proves \eqref{eq:keyboundaway}.
\end{proof}

With these results in hand, we can formulate and prove the following.
\begin{prop}
\label{prop:Hbound} 
Suppose that $\phi=e^{-V}$ where $V:S^1\rightarrow\mathbb{R}$ is of
class $C^{k-1,1}(S^1)$ for some $k=1,2,3,\dots$.  Let the integer $m$
lie in the range $1\le m\le k$ and fix $\epsilon>0$. Define the
integer $D:=\min(k-m,m-1)\ge 0$.  Then, for all $n\ge 0$ the matrix
${\bf W}^n_{m,\epsilon}(r,\theta)$ is well-defined almost everywhere
by
\eqref{eq:Wdef}, and for all $n$ sufficiently large, ${\bf
H}^n_{m,\epsilon}(r,\theta)$ is given by a 
Neumann series
\begin{equation}
{\bf H}_{m,\epsilon}^n(r,\theta)=\mathbb{I} +
(\mathcal{W}_{m,\epsilon}^n\mathbb{I})(r,\theta) +
(\mathcal{W}_{m,\epsilon}^n\circ\mathcal{W}_{m,\epsilon}^n\mathbb{I})(r,\theta)
+ \cdots
\end{equation}
which converges in the norm $\||\cdot|\|_D$, where the double integral
operator $\mathcal{W}_{m,\epsilon}^n$ is defined by
\begin{equation}
(\mathcal{W}^n_{m,\epsilon}{\bf
F})(r,\theta):=-\frac{1}{\pi}\int\int\frac{{\bf F}(r',\theta') {\bf
W}^n_{m,\epsilon}(r',\theta')}{z'-z}\,dA'\,.
\label{eq:Cnmdef}
\end{equation}
In particular, if $D=0$ then ${\bf H}^n_{m,\epsilon}$ lies in the
space $L^\infty(\mathbb{R}^2)$, and if $D>0$ then ${\bf
H}^n_{m,\epsilon}$ lies in the space $C^{D-1,1}(\mathbb{R}^2)$ and
$\||{\bf H}^n_{m,\epsilon}|\|_D$ is finite.  
For all integer $p$ in the range $0\le p\le D$, the following estimates hold
for sufficiently large $n$:
\begin{equation}
\||{\bf H}_{m,\epsilon}^n-\mathbb{I}|\|_p\le C_{m,\epsilon}^{(p)}
\frac{\log(n)}{n^{m-p}}\,,
\label{eq:HmI}
\end{equation}
\begin{equation}
\||{\bf H}_{m,\epsilon}^n-\mathbb{I}-\mathcal{W}_{m,\epsilon}^n\mathbb{I}
|\|_p\le C_{m,\epsilon}^{(p)2}\frac{\log(n)^2}{n^{2m-2p}}\,,
\label{eq:HmImWI}
\end{equation}
where $C_{m,\epsilon}^{(p)}>0$ is a constant.  Furthermore, for each
$\rho>2^\epsilon$ and for all integer $p$ in the range $0\le p\le D$,
the following estimates hold for sufficiently large $n$:
\begin{equation}
\sum_{\alpha+\beta\le p}
\mathop{\sup_{-\pi<\theta<\pi}}_{|\log(r)|\ge\log(\rho)}
\left\|\frac{\partial^{\alpha+\beta}}{\partial x^\alpha\partial y^\beta}
\left[{\bf H}^n_{m,\epsilon}(r,\theta)-\mathbb{I}\right]\right\|
\le \tilde{C}^{(p)}_{m,\rho}\frac{1}{n^{m-p}}\,,
\label{eq:HmIaway}
\end{equation}
\begin{equation}
\sum_{\alpha+\beta\le p}
\mathop{\sup_{-\pi<\theta<\pi}}_{|\log(r)|\ge\log(\rho)}
\left\|\frac{\partial^{\alpha+\beta}}{\partial x^\alpha\partial y^\beta}
\left[{\bf H}^n_{m,\epsilon}(r,\theta)-\mathbb{I}-(\mathcal{W}_{m,\epsilon}^n\mathbb{I})(r,\theta)\right]\right\|
\le \tilde{C}^{(p)2}_{m,\rho}\frac{\log(n)}{n^{2m-2p}}\,,
\label{eq:HmImWIaway}
\end{equation}
where $\tilde{C}^{(p)}_{m,\rho}>0$
is a constant.
\end{prop}
\begin{proof}
Fix $k\ge 1$ and $m$ in the range $1\le m\le k$, and set $D=\min(k-m,m-1)$.
Let $p$ be a nonnegative integer satisfying $p\le D$.  If $p>0$, suppose that ${\bf F}(r,\theta)$ is matrix function of class
$C^{p-1,1}(\mathbb{R}^2)$ with all derivatives of total order no
greater than $p$ uniformly bounded in the whole plane, and if $p=0$,
suppose that ${\bf F}(r,\theta)$ is of class $L^\infty(\mathbb{R}^2)$. For
such ${\bf F}(r,\theta)$ we recall the norm \eqref{eq:pnorm}, where the
cartesian coordinates $x$ and $y$ are connected to the polar coordinates $r$
and $\theta$ in the usual way:
$x=r\cos(\theta)$ and $y=r\sin(\theta)$.
Because for $z'\neq z$,
\begin{equation}
\frac{\partial^{\alpha+\beta}}{\partial x^\alpha\partial y^\beta}
\frac{1}{z'-z} = (-1)^{\alpha+\beta}
\frac{\partial^{\alpha+\beta}}{\partial x^{\prime\alpha}
\partial y^{\prime\beta}}
\frac{1}{z'-z}\,, 
\end{equation}
for $\alpha\ge 0$ and $\beta\ge 0$ with $\alpha+\beta\le p$, we have
\begin{equation}
\begin{array}{rcl}\displaystyle
\frac{\partial^{\alpha+\beta}}{\partial x^\alpha\partial y^\beta}
(\mathcal{W}_{m,\epsilon}^n{\bf F})(r,\theta)
&=&\displaystyle
 -\frac{(-1)^{\alpha+\beta}}{\pi}
\int\int{\bf F}(r',\theta'){\bf W}^n_{m,\epsilon}(r',\theta')
\frac{\partial^{\alpha+\beta}}{\partial x^{\prime\alpha}
\partial y^{\prime\beta}}\left[\frac{1}{z'-z}\right]\,dA'
\\\\&=&\displaystyle
-\frac{1}{\pi}\int\int\frac{\partial^{\alpha+\beta}}
{\partial x^{\prime\alpha}\partial y^{\prime\beta}}
\left[{\bf F}(r',\theta'){\bf W}_{m,\epsilon}^n(r',\theta')\right]
\frac{dA'}{z'-z}\,.
\end{array}
\label{eq:ibp1}
\end{equation}
Note that in order to integrate by parts in \eqref{eq:ibp1} for all
$\alpha$ and $\beta$ of interest, we must have $p\le k-m$.  Now, since
${\bf W}_{m,\epsilon}^n(r,\theta)$ is compactly supported in the
annulus $|\log(r)|\le \epsilon\log(2)$,
\begin{equation}
\begin{array}{rcl}
\displaystyle
\left\|\frac{\partial^{\alpha+\beta}}{\partial x^\alpha\partial y^\beta}
\left[{\bf F}(r,\theta){\bf W}^n_{m,\epsilon}(r,\theta)\right]\right\|
&\le &\displaystyle
\mathop{\sum_{\alpha'+\alpha''=\alpha}}_{\beta'+\beta''=\beta}
\left[\sup_{\mathbb{R}^2}\left\|\frac{\partial^{\alpha'+\beta'}}
{\partial x^{\alpha'}\partial y^{\beta'}}{\bf F}(r,\theta)\right\|\right]
\cdot \left\|\frac{\partial^{\alpha''+\beta''}}{\partial x^{\alpha''}
\partial y^{\alpha''}}{\bf W}_{m,\epsilon}^n(r,\theta)\right\|\\\\
&\le &\displaystyle
\||{\bf F}|\|_p\cdot
\sum_{\alpha''=0}^\alpha\sum_{\beta''=0}^\beta
\left\|\frac{\partial^{\alpha''+\beta''}}{\partial x^{\alpha''}
\partial y^{\beta''}}{\bf W}_{m,\epsilon}^n(r,\theta)\right\|\\\\
&\le &\displaystyle
\||{\bf F}|\|_p\sum_{\alpha'+\beta'\le \alpha+\beta} K_{m,\epsilon}^{(\alpha',\beta')}\left\|\frac{\partial^{\alpha'+\beta'}}
{\partial r^{\alpha'}\partial \theta^{\beta'}}{\bf W}^n_{m,\epsilon}(r,\theta)\right\|\,,
\end{array}
\end{equation}
where $K_{m,\epsilon}^{(\alpha',\beta')}$ are some positive constants.
Using Proposition~\ref{prop:Wbound}, we then find that
\begin{equation}
\left\|\frac{\partial^{\alpha+\beta}}{\partial x^\alpha\partial y^\beta}
\left[{\bf F}(r,\theta){\bf W}_{m,\epsilon}^n(r,\theta)\right]\right\|
\le 
K_{m,\epsilon}^{(\alpha+\beta)}\cdot
\||{\bf F}|\|_pe^{-n|\log(r)|}|\log(r)|^{m-1-(\alpha+\beta)}\sum_{j=0}^{\alpha+\beta}n^j|\log(r)|^j
\end{equation}
where $K_{m,\epsilon}^{(l)}>0$ is
some constant.  Finally, we arrive at the estimate
\begin{equation}
\left\|\left|\mathcal{W}_{m,\epsilon}^n{\bf F}\right|\right\|_p \le
\frac{\tilde{K}_{m,\epsilon}^{(p)}}{\pi}
\cdot\||{\bf F}|\|_p\cdot
\sum_{j=0}^pn^j\int_{2^{-\epsilon}}^{2^\epsilon}
r'\,dr'\,e^{-n|\log(r')|}|\log(r')|^{m+j-1-p}\int_{-\pi}^\pi\frac{d\theta'}
{|z'-z|}\,,
\label{eq:curlyWFintermediate}
\end{equation}
where $\tilde{K}_{m,\epsilon}^{(p)}>0$ 
is another constant.

Set $\nu:=m+j-p$, and note that $\nu\ge 1$ since as $j$ ranges from
$0$ to $p$, $\nu$ ranges from $m-p$ to $m$, and we have $p\le D\le
m-1$.  Lemma~\ref{lem:keybound} may thus be applied to each integral
on the right-hand side of \eqref{eq:curlyWFintermediate}, with the
result that
\begin{equation}
\left\|\left|\mathcal{W}_{m,\epsilon}^n{\bf F}\right|\right\|_p\le
\frac{C_{m,\epsilon}^{(p)}\log(n)}{n^{m-p}}\||{\bf F}|\|_p\,,
\label{eq:Cnmbound}
\end{equation}
for some constant $C_{m,\epsilon}^{(p)}>0$ and $n$ sufficiently large.
We note in passing that in order for \eqref{eq:Cnmbound} to provide
control of the operator $\mathcal{W}_{m,\epsilon}^n$, we need to have
$p<m$.  The two restrictions in force on $p$, namely $p<m$ and $p\le
k-m$, have been expressed in the statement of the proposition as the
inequality $p\le D$.  If one restricts attention to those
$z=re^{i\theta}$ for which $|\log(r)|\ge \log(\rho)>\epsilon\log(2)$,
then Lemma~\ref{lem:keybound} implies the inequality
\begin{equation}
\sum_{\alpha+\beta\le p}\mathop{\sup_{-\pi<\theta<\pi}}_{|\log(r)|\ge
\log(\rho)}\left\|\frac{\partial^{\alpha+\beta}}{\partial x^\alpha
\partial y^\beta}(\mathcal{W}_{m,\epsilon}^n{\bf F})(r,\theta)\right\|
\le \frac{\tilde{C}^{(p)}_{m,\rho}}{n^{m-p}}
\||{\bf F}|\|_p\,,
\label{eq:Cnmboundaway}
\end{equation}
where $\tilde{C}_{m,\rho}^{(p)}>0$
is a constant and $n$ is sufficiently large.

From \eqref{eq:Cnmbound} it is clear that if $n$ is sufficiently
large, the double integral operator $\mathcal{W}^n_{m,\epsilon}$
defined by the formula \eqref{eq:Cnmdef} and acting in the integral
equation \eqref{eq:inteqnE} thus defines a contraction mapping in the
space $C^{p-1,1}(\mathbb{R}^2)$ equipped with the norm
$\||\cdot|\|_p$, or in the space $L^\infty(\mathbb{R}^2)$ if $p=0$.
This implies that there is a unique solution of \eqref{eq:inteqnE} in
this space that may be found by iteration resulting in the
$\||\cdot|\|_p$-convergent Neumann series
\begin{equation}
{\bf U}(r,\theta):=\mathbb{I} + (\mathcal{W}^n_{m,\epsilon}\mathbb{I})(r,\theta) +
(\mathcal{W}^n_{m,\epsilon}\circ\mathcal{W}^n_{m,\epsilon}\mathbb{I})(r,\theta)+\cdots\,.
\label{eq:Neumann}
\end{equation}
In particular, choosing $p=0$ one sees that the Neumann series
\eqref{eq:Neumann} furnishes a unique solution of the integral
equation
\eqref{eq:inteqnE} in the space $L^\infty(\mathbb{R}^2)$.  
Since Proposition~\ref{prop:dbarEinteqn} guarantees that the matrix
${\bf H}^n_{m,\epsilon}(r,\theta)$ is a known solution of the integral
equation \eqref{eq:inteqnE} that is (in particular) uniformly bounded
in the plane, we may identify it with the Neumann series
\eqref{eq:Neumann} for $n$ sufficiently large.
From this point forward in our proof, we assume that $n$ is indeed
large enough for this to be the case. Since the same Neumann series
\eqref{eq:Neumann} also converges in the norm $\||\cdot|\|_p$ where
$p$ may be taken to be as large as $D$, we also learn that if $D>0$
then ${\bf H}_{m,\epsilon}^n(r,\theta)$ lies in the space
$C^{D-1,1}(\mathbb{R}^2)$ and that $\||{\bf H}^n_{m,\epsilon}|\|_D$ is
finite.

From \eqref{eq:Neumann}, taking $\||\cdot|\|_p$ norms
and using  \eqref{eq:Cnmbound}, we see that
\begin{equation}
\||{\bf H}_{m,\epsilon}^n-\mathbb{I}|\|_p\le
\|\mathbb{I}\|\sum_{j=1}^\infty\left(\frac{C^{(p)}_{m,\epsilon}
\log(n)}{n^{m-p}}\right)^j
\le 2\|\mathbb{I}\|\frac{C^{(p)}_{m,\epsilon}
\log(n)}{n^{m-p}}
\label{eq:crudebound}
\end{equation}
if $n$ is large enough that $n^{m-p}$ exceeds
$2C^{(p)}_{m,\epsilon}\log(n)$.  Since according
to the integral equation \eqref{eq:inteqnE} satisfied by ${\bf
H}^n_{m,\epsilon}(r,\theta)$ we have
\begin{equation}
{\bf H}^n_{m,\epsilon}(r,\theta)-\mathbb{I}-(\mathcal{W}_{m,\epsilon}^n\mathbb{I})(r,\theta) = (\mathcal{W}_{m,\epsilon}^n({\bf H}^n_{m,\epsilon}-\mathbb{I}))(r,\theta)
\label{eq:inteqnEagain}
\end{equation}
we may take norms and use \eqref{eq:Cnmbound} and \eqref{eq:crudebound} to learn that
\begin{equation}
\||{\bf H}_{m,\epsilon}^n-\mathbb{I}-\mathcal{W}_{m,\epsilon}^n\mathbb{I}
|\|_p\le 2\|\mathbb{I}\|
\frac{C^{(p)2}_{m,\epsilon}\log(n)^2}{n^{2m-2p}}
\label{eq:crudeboundtwoterms}
\end{equation}
holds for sufficiently large $n$.  The proof of the estimates
\eqref{eq:HmI} and \eqref{eq:HmImWI} is complete upon 
appropriate redefinition of the constant 
$C_{m,\epsilon}^{(p)}$.  

Note that \eqref{eq:Cnmboundaway} implies that the upper bounds in
\eqref{eq:crudebound} and \eqref{eq:crudeboundtwoterms} can be reduced
by a factor of $\log(n)$ if in each case the supremum on the left-hand
side is taken over only those values of $r$ satisfying
$|\log(r)|\ge\log(\rho)$ for a fixed $\rho>2^\epsilon$.  This
completes the proof of the estimates \eqref{eq:HmIaway} and
\eqref{eq:HmImWIaway} upon appropriate redefinition of the constant
$\tilde{C}_{m,\rho}^{(p)}$.
\end{proof}

\begin{remark}
  To uniformly control $p$ derivatives of ${\bf
    H}_{m,\epsilon}^n(r,\theta)$, Proposition~\ref{prop:Hbound}
  requires that $m$ should lie in the range $1+p\le m\le k-p$, and
  therefore in order for there to exist suitable values of $m$,
  $V:S^1\rightarrow\mathbb{R}$ should be of class $C^{k-1,1}(S^1)$ for
  some $k\ge 2p+1$.  Also, note that the utility of the estimates
  \eqref{eq:HmImWI} and \eqref{eq:HmImWIaway} is that the matrix ${\bf
    W}_{m,\epsilon}^n(r,\theta)$ is off-diagonal, so the diagonal
  matrix elements of ${\bf H}_{m,\epsilon}^n(r,\theta)-\mathbb{I}$
  experience more rapid decay than do the off-diagonal elements.
\end{remark}

\subsection{Proofs of theorems stated in \S~\ref{sec:strongasymptoticspi}.}
If we solve \eqref{eq:EM} for ${\bf M}^n(z)$ in terms of ${\bf
  H}_{m,\epsilon}^n(r,\theta)$, then since $M_{11}^n(z)=\pi_n(z)$, the
monic polynomial of degree $n$ in the system of polynomials orthogonal
with respect to the inner product $\langle\cdot,\cdot\rangle_\phi$
defined by \eqref{eq:innerproduct}, we can easily obtain from
Proposition~\ref{prop:Hbound} asymptotic formulae for $\pi_n(z)$ and
its derivatives, valid for large $n$, with uniform error estimates.

\subsubsection{Asymptotic behavior of $\pi_n(z)$ for $|z|>1$.  Proof
  of Theorem~\ref{thm:pioutside}.}
In the region $r>1$ we have for each
$\epsilon>0$ and for each $m=1,\dots,k$ the exact representation
\begin{equation}
\pi_n(z)=H_{m,\epsilon,11}^n(r,\theta)z^ne^{N(z)}-
H_{m,\epsilon,12}^n(r,\theta)
B(\log(r)/\epsilon)e^{E_mV(r,\theta)-N(z)}\,.
\label{eq:pinoutsideexact}
\end{equation}
Here we have used \eqref{eq:SeL} to write
$S_\phi(z)=e^{N(z)}$ for $|z|\ge 1$.  Equivalently,
\begin{equation}
\pi_n(z)z^{-n}e^{-N(z)}-1 = \left[H_{m,\epsilon,11}^n(r,\theta)-1\right] -
z^{-n}H_{m,\epsilon,12}^n(r,\theta)B(\log(r)/\epsilon)e^{E_mV(r,\theta)
-2N(z)}\,,\hspace{0.2 in}r\ge 1\,.
\label{eq:pioutsideequation}
\end{equation}
If $\rho>1$ is fixed, then we may choose $\epsilon>0$ small enough that 
$B(\log(r)/\epsilon)\equiv 0$ whenever $|z|\ge \rho$.  In this case,
we have simply
\begin{equation}
\pi_n(z)z^{-n}e^{-N(z)}-1 = 
H_{m,\epsilon,11}^n(r,\theta)-1\,,\hspace{0.2 in}
|z|\ge\rho>1\,.
\end{equation}
The best decay estimate comes from taking $m=k$.  In this case, using
Proposition~\ref{prop:Hbound} (specifically recalling the estimate
\eqref{eq:HmImWIaway} and the fact that ${\bf
  W}_{m,\epsilon}^n(r,\theta)$ is an off-diagonal matrix), we see that
for some constant $K_{\rho}>0$,
\begin{equation}
\mathop{\sup_{-\pi<\theta<\pi}}_{r\ge\rho>1}
\left|H_{m,\epsilon,11}^n(r,\theta)-1\right|\le 
K_{\rho}\frac{\log(n)}{n^{2k}}\,.
\label{eq:H11boundoutside}
\end{equation}
Now as the combination $\pi_n(z)z^{-n}e^{-N(z)}-1$ is a function of
$z$ that is analytic in the region $|z|>1$ and decaying to zero as
$z\rightarrow\infty$, 
we may express arbitrary derivatives of it as Cauchy integrals:
\begin{equation}
\frac{d^p}{dz^p}\left[\pi_n(z)z^{-n}e^{-N(z)}-1\right] = 
-\frac{p!}{2\pi i}\oint_{|s|=\rho}\frac{\pi_n(s)s^{-n}e^{-N(s)}-1}{(s-z)^p}
\,ds\,,
\end{equation}
where the contour of integration is oriented in the counterclockwise
direction, and $|z|>\rho$.  Using \eqref{eq:H11boundoutside} to bound
the integrand then gives a uniform bound of the same order of
magnitude for derivatives over regions bounded away from the circle
$|z|=\rho$, which can be taken arbitrarily close to the unit circle.
This proves \eqref{eq:piasympoutside}, and completes the proof of
Theorem~\ref{thm:pioutside}.

\subsubsection{Asymptotic behavior of $\pi_n(z)$ for $|z|= 1$.  
  Proof of Theorem~\ref{thm:picircle}.}  Theorem~\ref{thm:picircle}
follows from the estimate \eqref{eq:picircleestimate} by noting that
the error is an analytic function of $z$ in the exterior domain
$|z|>1$ that decays as $z\rightarrow\infty$, and therefore
\eqref{eq:picircleestimate} implies the more general result stated in
Theorem~\ref{thm:picircle} via the maximum modulus principle.  To
prove \eqref{eq:picircleestimate}, we fix any positive value of
$\epsilon$ and consider $1\le r\le 2^{\epsilon/2}$ in which case
$B(\log(r)/\epsilon)\equiv 1$ and therefore \eqref{eq:pioutsideequation}
implies the following formula:
\begin{equation}
\frac{d^p}{dz^p}\left[\pi_n(z)z^{-n}e^{-N(z)}-1\right] = 
\partial^p\left[H_{m,\epsilon,11}^n(r,\theta)-1\right]
-\partial^p\left[z^{-n}H_{m,\epsilon,12}^n(r,\theta)
e^{E_mV(r,\theta)-2N(z)}\right]\,,\hspace{0.2 in}
1\le r\le 2^{\epsilon/2}\,.
\label{eq:intermediate1}
\end{equation}
We remind the reader that estimates on derivatives like
\eqref{eq:intermediate1} are valid for $1+p\le m\le k-p$ (see the
remark at the end of \S~\ref{sec:dbarmethodfixed}). Using the estimate
\eqref{eq:HmImWI} from Proposition~\ref{prop:Hbound} and noting that
${\bf W}^n_{m,\epsilon}(r,\theta)$ is an off-diagonal matrix, we see
that for some constant $K_p>0$,
\begin{equation}
\sup_{\mathbb{R}^2}\left|\partial^p\left[H_{m,\epsilon,11}^n(r,\theta)-1
\right]\right|\le K_p\frac{\log(n)^2}{n^{2m-2p}}\,,
\end{equation}
if $m\le k-p$.  On the other hand, the dominant contributions actually
come from those terms in the second member of the right-hand side of
\eqref{eq:intermediate1} in which none of the $p$ derivatives fall on
the exponential factor $e^{E_mV(r,\theta)-2N(z)}$ (which has $k-m+1$
uniformly bounded derivatives).  Since $|z|\ge 1$, it suffices to
estimate $n^j\partial^{p-j} H_{m,\epsilon,12}^n(r,\theta)$ with the
use of the inequality \eqref{eq:HmI} in Proposition~\ref{prop:Hbound}.
Taking $m=k-p$ for the best possible decay estimate then gives
\begin{equation}
\mathop{\sup_{-\pi<\theta<\pi}}_{1\le r\le 2^{\epsilon/2}}
\left|\frac{d^p}{dz^p}\left[\pi_n(z)z^{-n}e^{-N(z)}-1\right]\right|
\le K_p\frac{\log(n)}{n^{k-2p}}\,,
\end{equation}
where $K_p>0$ is a constant.  This proves \eqref{eq:picircleestimate},
upon taking the limit $r\downarrow 1$ and writing the $z$ derivatives
in terms of $\theta$ (differentiation commutes with the limit
process).

\subsubsection{Asymptotic behavior of $\pi_n(z)$ for $|z|<1$
and of $\gamma_n$ .  Proof of Theorem~\ref{thm:piinside},
  Theorem~\ref{thm:gamma}, and Theorem~\ref{thm:piinsidejumps}.}
\label{sec:fixedinside}
Using \eqref{eq:EM} and the fact (see \eqref{eq:SeL}) that
$S_\phi(0)=e^{-V_0}$, we have the exact representation:
\begin{equation}
\gamma_{n-1}^2 = -M_{21}^n(0)=H_{m,\epsilon,22}^n(0,\theta)e^{V_0}\,,
\label{eq:gammaMformula}
\end{equation}
and, whenever $|z|<1$,
\begin{equation}
\pi_n(z)=M_{11}^n(z)=z^ne^{-V_0-\overline{N(1/\overline{z})}}
B(\log(r)/\epsilon)e^{E_mV(r,\theta)}H^n_{m,\epsilon,11}(r,\theta)
-e^{V_0+\overline{N(1/\overline{z})}}H^n_{m,\epsilon,12}(r,\theta)\,,
\label{eq:piinsideMformula}
\end{equation}
where \eqref{eq:SeL} has been used, and $z=re^{i\theta}$.

To prove Theorem~\ref{thm:gamma}, we simply apply
\eqref{eq:HmImWIaway} from Proposition~\ref{prop:Hbound} in the case
$p=0$ and $m=k$ to the identity \eqref{eq:gammaMformula}. This
immediately yields \eqref{eq:fixedgammaasymp} and completes the proof.

The proof of Theorem~\ref{thm:piinside} is based on a similar analysis
of \eqref{eq:piinsideMformula}.  Recalling that $\rho\in (0,1)$, one
can choose $\epsilon>-2\log(\rho)/\log(2)$ and then
$B(\log(r)/\epsilon)\equiv 1$ in \eqref{eq:piinsideMformula} for
$\rho<|z|<1$.  The estimate \eqref{eq:pinearcircleinsideestimate} then
follows by taking $m=k$, and using \eqref{eq:HmImWI} from
Proposition~\ref{prop:Hbound} in the case $p=0$.  Similarly, choosing
$\epsilon<-\log(\rho)/\log(2)$, we have $B(\log(r)/\epsilon)\equiv 0$
in \eqref{eq:piinsideMformula} for $|z|<\rho$.  Again taking $m=k$,
one obtains \eqref{eq:piawaycircleinsideestimate} by using
\eqref{eq:HmImWIaway} from Proposition~\ref{prop:Hbound} with $p=0$.
This completes the proof of Theorem~\ref{thm:piinside}.

The rest of this section will be devoted to the proof of
Theorem~\ref{thm:piinsidejumps}.  We begin with
\eqref{eq:piinsideMformula} for $m=k$, a formula that is valid for
all $z$ with $|z|<1$.
Using
\eqref{eq:HmImWI} from Proposition~\ref{prop:Hbound}, and keeping the
term corresponding to $\mathcal{W}_{k,\epsilon}^n\mathbb{I}$, we arrive
at the formula
\begin{equation}
\pi_n(z)=z^ne^{-V_0-\overline{N(1/\overline{z})}}e^{E_kV(r,\theta)}
B(\log(r)/\epsilon)
+\frac{e^{V_0+\overline{N(1/\overline{z})}}}{\pi}
\mathop{\int\int}_{2^{-\epsilon}<r'<1}\frac{W_{k,\epsilon,12}^n(r',\theta')}
{z'-z}\,dA' + O\left(\frac{\log(n)^2}{n^{2k}}\right)\,,
\label{eq:doubleinteqn}
\end{equation}
as $n\rightarrow\infty$, where the error term is uniform for $|z|<1$.
Here, $dA'=r'\,dr'\,d\theta'$ is an area element,
and in the integral $\theta'$ varies over $S^1$
(the support of $W_{k,\epsilon,12}^n(r,\theta)$ is the annulus
given by the inequalities $2^{-\epsilon}\le r\le 1$).  

In the annulus of support of $W_{k,\epsilon,12}^n(r,\theta)$, we have
from \eqref{eq:Wdef} that
\begin{equation}
W_{k,\epsilon,12}^n(r,\theta)=
\dbar\left[z^ne^{-2V_0-2\overline{N(1/\overline{z})}}
B(\log(r)/\epsilon)e^{E_kV(r,\theta)}\right]\,.
\label{eq:Wdbarformula}
\end{equation}

\begin{remark}
  The fact that $W_{k,\epsilon,12}^n(r,\theta)$ is in the range of
  $\dbar$ means that the double integral in \eqref{eq:doubleinteqn}
  can be reduced without approximation to the sum of an explicit
  contribution and a contour integral.  Indeed, by the inversion of
  the $\dbar$ operator, 
\begin{equation}
-\frac{1}{\pi}\mathop{\int\int}_{2^{-\epsilon}<r'<1}
\frac{W^n_{k,\epsilon,12}(r',\theta')}
{z'-z}\,dA' = F(z)+z^ne^{-2V_0-2\overline{N(1/\overline{z})}}
B(\log(r)/\epsilon)e^{E_kV(r,\theta)}\chi_{(2^{-\epsilon},1)}(r)\,,
\end{equation}
where $\chi_I$ denotes the characteristic function of an interval $I$,
and where $F(z)$ is function analytic except on the circles
$|z|=2^{-\epsilon}$ and $|z|=1$ bounding the support of
$W_{k,\epsilon,12}^n(r,\theta)$ that is chosen to make the right-hand
side continuous and decaying as $z\rightarrow\infty$.  These latter properties
uniquely identify $F(z)$ with the Cauchy integral
\begin{equation}
F(z)=-\frac{1}{2\pi i}\oint_{|s|=1}\frac{s^ne^{-V_0+i\Omega(\arg(s))}}{s-z}\,ds\,,
\end{equation}
where the contour of integration is oriented counterclockwise.  Note
that it is the presence of the bump function $B(\log(r)/\epsilon)$
that makes $F(z)$ continuous at $|z|=2^{-\epsilon}$.  Unfortunately,
this interesting formula, while apparently simpler than a double
integral, is not as useful for asymptotic analysis as the alternative
approach we now follow.
\end{remark}

Continuing our analysis, we carry out the differentiation in
\eqref{eq:Wdbarformula} in the region $|z|<1$ with the use of 
\eqref{eq:dbarTaylor}:
\begin{equation}
\begin{array}{rcl}
\displaystyle
W_{k,\epsilon,12}^n(r,\theta)&=&
\displaystyle
z^ne^{-2V_0-2\overline{N(1/\overline{z})}}
e^{E_kV(r,\theta)}\left[\dbar B(\log(r)/\epsilon) + 
B(\log(r)/\epsilon)\dbar E_kV(r,\theta)\right]\\\\
&=&\displaystyle
z^ne^{-2V_0-2\overline{N(1/\overline{z})}}
e^{E_kV(r,\theta)}\left[\dbar B(\log(r)/\epsilon) + 
\frac{ie^{i\theta}}{2r}\frac{V^{(k)}(\theta)}{(k-1)!}(-i\log(r))^{k-1}
B(\log(r)/\epsilon)\right]\,,
\end{array}
\end{equation}
and in the special case that $r=|z|>2^{-\epsilon/2}$, we have 
$B(\log(r)/\epsilon)\equiv 1$ and $\dbar B(\log(r)/\epsilon)\equiv 0$, so
\begin{equation}
W_{k,\epsilon,12}^n(r,\theta)=
z^ne^{-2V_0-2\overline{N(1/\overline{z})}}
e^{E_kV(r,\theta)}\cdot
\frac{ie^{i\theta}}{2r}\frac{V^{(k)}(\theta)}{(k-1)!}(-i\log(r))^{k-1}
\,,\hspace{0.2 in}\text{for $2^{-\epsilon/2}<r<1$}\,.
\label{eq:Wsimplerformula}
\end{equation}
The presence of the $z^n$ factor together with the absolute integrability
of the Cauchy kernel in two dimensions means that
\begin{equation}
-\frac{1}{\pi}\mathop{\int\int}_{2^{-\epsilon}<r'<1}\frac{W^n_{k,\epsilon,12}
(r',\theta')}{z'-z}\,dA' = 
-\frac{1}{\pi}\mathop{\int\int}_{2^{-\epsilon/2}<r'<1}\frac{W^n_{k,\epsilon,12}
(r',\theta')}{z'-z}\,dA' + O(2^{-\epsilon n/2})
\end{equation}
holds as $n\rightarrow\infty$ uniformly for all $z$ with $|z|<1$.  Therefore
using the simpler formula \eqref{eq:Wsimplerformula} in the integrand
and integrating over the smaller annulus $2^{-\epsilon/2}<r'<1$ introduces
an error that is uniformly exponentially small for $|z|<1$.  

Using the simple identity
\begin{equation}
\dbar (-\log(r))^k = -\frac{kz}{2r^2}(-\log(r))^{k-1}\,,
\label{eq:dbarlogidentity}
\end{equation}
we may rewrite \eqref{eq:Wsimplerformula} in the form
\begin{equation}
W_{k,\epsilon,12}^n(r,\theta)=-\frac{i^k}{k!}z^n
e^{-2V_0-2\overline{N(1/\overline{z})}}e^{E_kV(r,\theta)}V^{(k)}(\theta)
\dbar(-\log(r))^k\,,\hspace{0.2 in}\text{for $2^{-\epsilon/2}<r<1$}\,.
\end{equation}
Our subsequent analysis will be specialized to the case where $V^{(k)}(\theta)$
is piecewise continuous, with jump discontinuities at $\ell<\infty$ angles
$-\pi<\theta_1<\dots<\theta_\ell<\pi$, and is (at first) only Lipschitz between
the points of discontinuity.  Then, between the points of discontinuity,
$V^{(k+1)}(\theta)$ exists almost everywhere and may be identified with
a bounded function.  Under these circumstances, we may ``integrate by parts''
({\em i.e.} apply Stokes' Theorem) with the following formula
\begin{equation}
-\frac{1}{\pi}\int\int\frac{f(r',\theta')\dbar g(r',\theta')}{z'-z}\,dA' = 
-\frac{1}{\pi}\int\int\frac{\dbar[f(r',\theta')g(r',\theta')]}{z'-z}\,dA'
+\frac{1}{\pi}\int\int\frac{g(r',\theta')\dbar f(r',\theta')}{z'-z}\,dA'\,,
\label{eq:dbaribp}
\end{equation}
and the first integral on the right-hand side may be exchanged for a sum
of explicit terms and a contour integral as described in the above remark.
In \eqref{eq:dbaribp} and in the rest of the proof, whenever the operator
$\dbar$ appears in the integrand it acts on the primed variables.

To prepare to use this technique, we begin with
\begin{equation}
-\frac{1}{\pi}\mathop{\int\int}_{2^{-\epsilon/2}<r'<1}
\frac{W^n_{k,\epsilon,12}(r',\theta')}{z'-z}\,dA'=
-\frac{i^k}{k!}\sum_{j=1}^\ell\left[-\frac{1}{\pi}
\mathop{\mathop{\int\int}_{2^{-\epsilon/2}<r'<1}}_{\theta'\in I_j}
\frac{(z')^nh(r',\theta')
V^{(k)}(\theta')\dbar(-\log(r'))^k}{z'-z}\,dA'\right]\,,
\label{eq:Wintegralsum}
\end{equation}
where $h(r,\theta)$ is shorthand for the following terms:
\begin{equation}
h(r,\theta):=
e^{-2V_0-2\overline{N(1/\overline{z})}}e^{E_kV(r,\theta)}
\,,
\end{equation}
and $I_j$ refers to the interval in $S^1$ of initial angle $\theta_j$
and final angle equal to the point of next jump discontinuity as the
circle is traversed in the counterclockwise direction.  Using 
\eqref{eq:dbaribp}, we therefore find
\begin{equation}
-\frac{1}{\pi}\mathop{\int\int}_{2^{-\epsilon/2}<r'<1}
\frac{W^n_{k,\epsilon,12}(r',\theta')}{z'-z}\,dA'=
-J_k(r,\theta) + K_k(r,\theta)\,,
\end{equation}
where
\begin{equation}
\begin{array}{rcl}
J_k(r,\theta)&:=&\displaystyle\frac{i^k}{k!}\sum_{j=1}^\ell
\left[-\frac{1}{\pi}
\mathop{\mathop{\int\int}_{2^{-\epsilon/2}<r'<1}}_{\theta'\in I_j}
\frac{(z')^n\dbar[h(r',\theta')V^{(k)}(\theta')](-\log(r'))^k}{z'-z}\,dA'
\right]\\\\
K_k(r,\theta)&:=&\displaystyle\frac{i^k}{k!}\sum_{j=1}^\ell
\left[-\frac{1}{\pi}
\mathop{\mathop{\int\int}_{2^{-\epsilon/2}<r'<1}}_{\theta'\in I_j}
\frac{\dbar\left[(z')^nh(r',\theta')V^{(k)}(\theta')(-\log(r'))^k\right]}
{z'-z}\,dA'\right]\,.
\end{array}
\end{equation}
The double-integral expression $K_k(r,\theta)$ may be reduced to contour
integrals as follows:
\begin{equation}
K_k(r,\theta)=\frac{i^k}{k!}z^nh(r,\theta)V^{(k)}(\theta)(-\log(r))^k
\chi_{(2^{-\epsilon/2},1)}(r) +G_k(z)\,,
\label{eq:Kkrepn}
\end{equation}
where
\begin{equation}
\begin{array}{rcl}
G_k(z)&:=&\displaystyle
-\frac{i^k}{k!}\cdot\frac{1}{2\pi i}\oint_{|s|=2^{-\epsilon/2}}\frac{s^n
h(|s|,\arg(s))V^{(k)}(\arg(s))(-\log(|s|))^k}
{s-z}\,ds\\\\
&&\displaystyle \,\,\,+\,\,\,\frac{i^k}{k!}\sum_{j=1}^\ell
\frac{\Delta_j^{(k)}}{2\pi i}
\int_{2^{-\epsilon/2}e^{i\theta_j}}^{e^{i\theta_j}}
\frac{s^nh(|s|,\theta_j)(-\log(|s|))^k}{s-z}\,ds\,,
\end{array}
\end{equation}
and $\Delta_j^{(k)}:=V^{(k)}(\theta_j+)-V^{(k)}(\theta_j-)$.  Since
the $\||\cdot|\|_{\circ,1}$ norm of the numerator in the first Cauchy
integral is proportional to $2^{-\epsilon n/2}$, we may also write
$G_k(z)$ in the form
\begin{equation}
G_k(z)=
\frac{i^k}{k!}\sum_{j=1}^\ell
\frac{\Delta_j^{(k)}}{2\pi i}
\int_{2^{-\epsilon/2}e^{i\theta_j}}^{e^{i\theta_j}}
\frac{s^nh(|s|,\theta_j)(-\log(|s|))^k}{s-z}\,ds + O(2^{-\epsilon n/2})
\label{eq:Gkrepn}
\end{equation}
as $n\rightarrow \infty$ where the exponentially small error term is
uniform for $|z|<1$ (right up to the contour of integration; this
is a consequence of the Plemelj-Privalov Theorem \cite{Musk}).

Let us now consider $J_k(r,\theta)$.  Note that
\begin{equation}
\begin{array}{rcl}
\displaystyle
\dbar[h(r,\theta)V^{(k)}(\theta)]&=&\displaystyle
 \left[V^{(k)}(\theta)
\dbar E_kV(r,\theta) + \frac{iz}{2r^2}V^{(k+1)}(\theta)\right]
h(r,\theta)\\\\
&=&\displaystyle
\frac{iz}{2r^2}\left[\frac{V^{(k)}(\theta)^2}{(k-1)!}(-i\log(r))^{k-1}
+V^{(k+1)}(\theta)\right]h(r,\theta)
\,.
\end{array}
\end{equation}
Therefore, inserting this into the integrand for $J_k(r,\theta)$ and
applying Lemma~\ref{lem:keybound} with $\nu=2k$ to the integrals resulting
from the first term above, we
find
\begin{equation}
\begin{array}{rcl}
J_k(r,\theta)&=&\displaystyle
-\frac{i^{k+1}}{(k+1)!}\sum_{j=1}^\ell
\left[-\frac{1}{\pi}
\mathop{\mathop{\int\int}_{2^{-\epsilon/2}<r'<1}}_{\theta'\in I_j}
\left[-\frac{(k+1)z'}{2(r')^2}(-\log(r'))^k\right]
\frac{(z')^nh(r',\theta')V^{(k+1)}(\theta')
}
{z'-z}\,dA'\right]\\\\
&&\displaystyle\,\,\, +\,\,\, O\left(\frac{\log(n)}{n^{2k}}\right)\,,
\end{array}
\end{equation}
where the error is uniformly small as $n\rightarrow\infty$ for all
$|z|<1$.  With $V^{(k+1)}(\theta)$ bounded in each $I_j$, we could in
principle apply Lemma~\ref{lem:keybound} to the remaining integrals.
However, this would only give a bound of order $\log(n)/n^{k+1}$, and
as $G_k(z)$ will turn out to be (for most $z$) of size $1/n^{k+1}$, and
we will want to consider $G_k(z)$ to provide the dominant term, we need
to impose additional conditions on $V^{(k+1)}$ in each $I_j$ to see
that $J_k(r,\theta)$ is indeed subdominant.  Therefore, we first use
\eqref{eq:dbarlogidentity} to write $J_k(r,\theta)$ in the form
\begin{equation}
J_k(r,\theta)=-\frac{i^{k+1}}{(k+1)!}\sum_{j=1}^\ell
\left[-\frac{1}{\pi}
\mathop{\mathop{\int\int}_{2^{-\epsilon/2}<r'<1}}_{\theta'\in I_j}
\frac{(z')^nh(r',\theta')V^{(k+1)}(\theta')\dbar(\log(r'))^{k+1}}
{z'-z}\,dA'\right] + O\left(\frac{\log(n)}{n^{2k}}\right)\,.
\end{equation}
The double integral above is of the same form as the original double
integral we are trying to compute (see \eqref{eq:Wintegralsum}), but
with $k$ replaced by $k+1$ (everwhere it appears explicitly; in the
function $h(r,\theta)$, $k$ remains $k$).

To continue the analysis of $J_k(r,\theta)$, we therefore use the
further assumption that in each each of the intervals $I_j$,
$V^{(k+1)}(\theta')$ is a Lipschitz continuous function so that
$V^{(k+2)}(\theta')$ exists almost everywhere in $I_j$ and can be
identified there with a bounded function.  Repeating the above steps,
we find that
\begin{equation}
J_k(r,\theta)=-J_{k+1}(r,\theta)+K_{k+1}(r,\theta)\,,
\end{equation}
and with the help of the identity
\begin{equation}
\dbar[h(r,\theta)V^{(k+1)}(\theta)]=\frac{iz}{2r^2}\left[
\frac{V^{(k)}(\theta)V^{(k+1)}}{(k-1)!}(-i\log(r))^{k-1}
+V^{(k+2)}(\theta)\right]h(r,\theta)\,,
\end{equation}
we may apply Lemma~\ref{lem:keybound} to the integrals that result
from substituting the above into $J_{k+1}(r,\theta)$ with $\nu=2k+1$
(from the first term above) and with $\nu=k+2$ (from the second term above).
Consequently, 
\begin{equation}
J_{k+1}(r,\theta)=O\left(\frac{\log(n)}{n^{k+2}}\right)\,,
\end{equation}
where the error is uniformly small as $n\rightarrow\infty$ for all
$r<1$.  As before, $K_{k+1}(r,\theta)$ is given by \eqref{eq:Kkrepn}
with $G_{k+1}(z)$ given by \eqref{eq:Gkrepn} (note that in
substituting $k+1$ for $k$ in these formulae, one leaves $h(r,\theta)$
alone).

Combining these results, we have shown that
\begin{equation}
\begin{array}{rcl}\displaystyle
-\frac{1}{\pi}
\mathop{\int\int}_{2^{-\epsilon/2}<r'<1}\frac{W^n_{k,\epsilon,12}(r',\theta')}
{z'-z}\,dA' &=&\displaystyle 
\frac{i^k}{k!}\sum_{j=1}^\ell
\frac{\Delta_j^{(k)}}{2\pi i}
\int_{2^{-\epsilon/2}e^{i\theta_j}}^{e^{i\theta_j}}
\frac{s^nh(|s|,\theta_j)(-\log(|s|))^k}{s-z}\,ds \\\\
&&\displaystyle\,\,\,-\,\,\,
\frac{i^{k+1}}{(k+1)!}\sum_{j=1}^l
\frac{\Delta_j^{(k+1)}}{2\pi i}
\int_{2^{-\epsilon/2}e^{i\theta_j}}^{e^{i\theta_j}}
\frac{s^nh(|s|,\theta_j)(-\log(|s|))^{k+1}}{s-z}\,ds \\\\
&&\displaystyle\,\,\,+\,\,\, O\left(\frac{\log(n)}{n^{k+2}}\right) + 
O\left(r^n|\log(r)|^k\chi_{(2^{-\epsilon/2},1)}(r)\right)
\end{array}
\label{eq:Wintegralcombined}
\end{equation}
holds uniformly for $r=|z|<1$ under the assumptions in force on $V(\theta)$.

Now in addition to $|z|<1$ we choose arbitrarily a constant
$\sigma>0$ and consider those $z$ for which $\log(|z|)\le-
(k-\sigma)\log(n)/n$ (this is the interesting case, since
according to Corollary~\ref{cor:fixedzerofree} we are excluding a
zero-free annulus near the unit circle whenever
$\sigma>\delta$ where $\delta$ is the arbitrary positive
parameter in the statement of Corollary~\ref{cor:fixedzerofree}; note
that both $\delta$ and $\sigma$ may be taken to be arbitrarily
small).  Since the function $r^n|\log(r)|^k$ achieves its maximum value
when $|\log(r)|$ is proportional to $1/n$, we have
\begin{equation}
\max_{\log(r)\le-(k-\sigma)\log(n)/n}r^n|\log(r)|^k
= r^n|\log(r)|^k\Bigg|_{\log(r)=-(k-\sigma)\log(n)/n}
=(k-\sigma)^kn^{\sigma-2k}(\log(n))^k\,.
\end{equation}
Since we are assuming here that $k\ge 2$, $\sigma>0$ may be chosen
small enough that both error terms in \eqref{eq:Wintegralcombined}
may be replaced by $o(n^{-(k+1)})$ as $n\rightarrow\infty$ uniformly
for $\log(|z|)\le -(k-\sigma)\log(n)/n$.

It remains to evaluate the explicit integrals in \eqref{eq:Wintegralcombined}
by Laplace's method.  Letting $a=k$ or $a=k+1$ we consider
\begin{equation}
Y_{a,j}(z):=\int_{2^{-\epsilon/2}e^{i\theta_j}}^{e^{i\theta_j}}
\frac{s^nh(|s|,\theta_j)(-\log(|s|))^a}{s-z}\,ds = 
e^{in\theta_j}\int_{2^{-\epsilon/2}}^1\frac{x^nh(x,\theta_j)(-\log(x))^a}
{x-ze^{-i\theta_j}}\,dx\,.
\end{equation}
Expecting the dominant contribution to come from the neighborhood of $x=1$,
we write
\begin{equation}
Y_{a,j}(z)=\frac{e^{in\theta_j}h(1,\theta_j)}{1-ze^{-i\theta_j}}
\int_{2^{-\epsilon/2}}^1x^n(-\log(x))^a\,dx
+e^{in\theta_j}
\int_{2^{-\epsilon/2}}^1\left[\frac{h(x,\theta_j)}{x-ze^{-i\theta_j}}-
\frac{h(1,\theta_j)}{1-ze^{-i\theta_j}}\right]x^n(-\log(x))^a\,dx\,.
\end{equation}
Finding a common denominator and extracting from $Y_{a,j}(z)$ a factor
of $e^{in\theta_j}(1-ze^{-i\theta_j})^{-1}$, we see that
\begin{equation}
\begin{array}{rcl}
\displaystyle Y_{a,j}(z)&=&\displaystyle
\frac{e^{in\theta_j}}{1-ze^{-i\theta_j}}\Bigg[
h(1,\theta_j)\int_{2^{-\epsilon/2}}^1 x^n(-\log(x))^a\,dx  \\\\
&&\displaystyle\,\,\,+\,\,\,
\int_{2^{-\epsilon/2}}^1 
\frac{ze^{-i\theta_j}(h(1,\theta_j)-h(x,\theta_j))+
h(x,\theta_j)-xh(1,\theta_j)}{x-ze^{-i\theta_j}}x^n(-\log(x))^a\,dx\Bigg]\,.
\end{array}
\end{equation}
If $z$ does not approach the point $e^{i\theta_j}$, then the fraction
in the integrand of the second integral is easily seen to be bounded
by a multiple of $-\log(x)$ that is independent of $n$.  More
generally, if $z$ is allowed to approach the point $e^{i\theta_j}$ as
$n\rightarrow\infty$, then we note that the integrand is analytic in
$x$, and the path of integration may be deformed in such a way that
throughout the path of integration $|x-ze^{-i\theta_j}|$ is bounded
away from zero by a quantity that is proportional to $\log(n)/n$
because $\log(|z|)\le -(k-\sigma)\log(n)/n$, and therefore the
fraction in the integrand of the second integral is bounded by
a multiple of $-\log(x)$ that is proportional to $n/\log(n)$.  Therefore,
\begin{equation}
Y_{a,j}(z)=\frac{e^{in\theta_j}h(1,\theta_j)}{1-ze^{-i\theta_j}}
\frac{a!}{n^{a+1}}\left(1+O\left(\frac{1}{\log(n)}\right)\right)\,,
\end{equation}
where the error is uniformly small for $\log(|z|)\le
-(k-\sigma)\log(n)/n$.  Since $E_kV(1,\theta)=V(\theta)$, we have
$h(1,\theta)=e^{-V_0+i\Omega(\theta)}$.  It follows that
\begin{equation}
-\frac{1}{\pi}\mathop{\int\int}_{2^{-\epsilon/2}<r'<1}
\frac{W^n_{k,\epsilon,12}(r',\theta')}{z'-z}\,dA' = 
\frac{i^{k-1}e^{-V_0}}{2\pi n^{k+1}}\sum_{j=1}^\ell 
\Delta_j^{(k)}e^{i\Omega(\theta_j)}\frac{e^{i(n+1)\theta_j}}{e^{i\theta_j}-z}
+ o(n^{-(k+1)})\,,
\end{equation}
holds uniformly in the region $\log(|z|)<-(k-\sigma)\log(n)/n$, and
therefore by \eqref{eq:doubleinteqn} so does
\begin{equation}
\pi_n(z)=z^ne^{-V_0-\overline{N(1/\overline{z})}}e^{E_kV(r,\theta)}
B(\log(r)/\epsilon)+
\frac{i^{k+1}e^{\overline{N(1/\overline{z})}}}{2\pi n^{k+1}}\sum_{j=1}^\ell 
\Delta_j^{(k)}e^{i\Omega(\theta_j)}\frac{e^{i(n+1)\theta_j}}{e^{i\theta_j}-z}
+ o(n^{-(k+1)})\,.
\end{equation}
This concludes the proof of Theorem~\ref{thm:piinsidejumps}.

\section{Exponentially Varying Weights}
\label{sec:varying}
\subsection{Asymptotic formulae for $\pi_n(z)$ and $\gamma_n$ in the
varying weights case.}
\label{sec:asymptoticspids}
In this section, we consider weights of the form
\begin{equation}
\phi(\theta)=e^{-nV(\theta)}\,,\hspace{0.2 in}\text{for all $\theta\in S^1$}\,,
\label{eq:varyingweight}
\end{equation}
where $V:S^1\rightarrow\mathbb{R}$ is a given real-valued function of
period $2\pi$.  The weight \eqref{eq:varyingweight} varies
exponentially according to a parameter $n$, and for each $n$ we may
associate with $\phi$ the corresponding sequence of monic orthogonal
polynomials $\pi_0(z), \pi_1(z),\pi_2(z),\dots$ and normalization
constants $\gamma_0, \gamma_1, \gamma_2,\dots$.  Properly speaking,
these quantities depend on the parameter $n$, and we should invent
notation to express this dependence, such as $\pi_{j}^{(n)}(z)$.
We will not introduce this cumbersome notation.  However, the reader
should take note of this dependence.
The limit of interest here is to study the
behavior of the particular monic polynomial $\pi_n(z)$ of degree $n$
in this system along with its normalization constant $\gamma_n$, in
the limit $n\rightarrow\infty$.  Thus the large parameter $n$ enters
simultaneously into the degree of the polynomial and also into the
weight, and we are studying the asymptotic behavior along the diagonal
of a doubly indexed sequence.

The asymptotic behavior in this limit is governed by the
function $V$, along with some associated functions.  First, recall the
analytic function $N(z)$ defined by \eqref{eq:negfreq} for $|z|\ge 1$,
which is associated with the negative frequency component of the
Fourier series of $V(\theta)$.  Now define
\begin{equation}
\kappa(\theta):=\theta - 
i\left[N(e^{i\theta})-\overline{N(e^{i\theta})}\right]
=\theta+\Omega(\theta)\,,
\label{eq:kappadefine}
\end{equation}
where $\Omega(\theta)$ is the periodic function defined
in \eqref{eq:Omegadef}.
Both functions $\kappa$ and $\Omega$ are as smooth as $V$ is.

Our asymptotic results in this case are the following.
\begin{theorem}
Let $p\ge 0$ be a fixed integer.  Suppose that
$\phi(\theta)=e^{-nV(\theta)}$ where $V:S^1\rightarrow\mathbb{R}$ is
of class $C^{k-1,1}(S^1)$ for some $k\ge 2p+2$.  If $\kappa'(\theta)$
is strictly positive, then for each $\rho > 1$ there is a constant
$K_{p,\rho}>0$ such that 
\begin{eqnarray}
  \label{eq:varypiasympout2}
\sup_{|z|\ge \rho}\left|\frac{d^p}{dz^p}\left[
\pi_n(z)z^{-n}e^{-nN(z)}-1\right]\right|\le
K_{p, \rho}\frac{\log(n)}{n^{2(k-1)}}  
\end{eqnarray}
holds for all sufficiently large $n$. 
\label{thm:varypiasympoutbasic}
\end{theorem}
\begin{theorem}
Let $p\ge 0$ be a fixed integer.  Suppose that
$\phi(\theta)=e^{-nV(\theta)}$ where $V:S^1\rightarrow\mathbb{R}$ is
of class $C^{k-1,1}(S^1)$ for some $k\ge 2p+2$.  If $\kappa'(\theta)$
is strictly positive, then there is a constant $K_p>0$ such that
\begin{equation}
\sup_{|z|\ge 1}\left|\frac{d^p}{dz^p}\left[
\pi_n(z)z^{-n}e^{-nN(z)}-1\right]\right|\le
K_{p}\frac{\log(n)}{n^{k-2p-1}}
\label{eq:varyingpiasympoutside}
\end{equation}
holds for all sufficiently large $n$.
\label{thm:picirclevaryingweights}
\end{theorem}
\begin{remark}
As in the fixed weights case (see \eqref{eq:picircleestimate}), a special
case of the estimate \eqref{eq:varyingpiasympoutside} is the following
estimate (holding under the same conditions) 
\begin{equation}
\sup_{-\pi<\theta<\pi}\left|\left(-ie^{-i\theta}
\frac{d}{d\theta}\right)^p\left[\pi_n(e^{i\theta})e^{-in\theta}e^{-nN(e^{i\theta})}
- 1\right]\right|\le K_p\frac{\log(n)}{n^{k-2p-1}}.
\label{eq:varyingpiasympcircle}
\end{equation}
asymptotically characterizing $\pi_{n}(z)$ on the unit circle.
\end{remark}

Once again, the uniform nature of the convergence on the circle allows
us to prove the following mean result (the proof is the same as for
the analogous result in the fixed weight case).
\begin{corollary}
  Let $p\ge 0$ be a fixed integer.  Suppose that
  $\phi(\theta)=e^{-nV(\theta)}$ where $V:S^1\rightarrow\mathbb{R}$ is
  of class $C^{k-1,1}(S^1)$ with $k\ge 2p+2$.  If $\kappa'(\theta)$ is
  strictly positive, then
\begin{equation}
\lim_{n\rightarrow\infty}\frac{1}{n^p}\cdot\frac{\|\pi_n^{(p)}\|_\phi}
{\|\pi_n\|_\phi} =1 \,.
\end{equation}
\end{corollary}

The next result concerns the asymptotic behavior of the polynomial
$\pi_{n}(z)$ for $z$ inside a closed annular region whose outer
boundary is the unit circle.
  
\begin{theorem}
  Suppose that $\phi(\theta)=e^{-nV(\theta)}$ where
  $V:S^1\rightarrow\mathbb{R}$ is of class $C^{k-1,1}(S^1)$ for some
  $k\ge 2$.  If $\kappa'(\theta)$ is strictly positive, then for each
  $\rho$ satisfying $0<\rho<1$ there are constants $K^\pm_\rho>0$ such
  that the estimates
\begin{equation}
\sup_{\rho<|z|<1}\left|\pi_n(z)e^{-n\overline{N(1/\overline{z})}}-z^n
e^{inE_k\Omega(r,\theta)} \right|\le K^-_\rho\frac{\log(n)}{n^{k-1}}\,,
\label{eq:pivaryingarcircinsideestimate}
\end{equation}
and
\begin{equation}
\sup_{|z|<\rho}\left|\pi_n(z)
e^{-n\overline{N(1/\overline{z})}}\right|\le \frac{K^+_\rho}{n^{k-1}}
\label{eq:pivaryawaycircleinnerestimate}
\end{equation}
hold for all $n$ sufficiently large.
\label{thm:varyingcircleattract}
\end{theorem}

An immediate corollary is that there exists an annulus inside the unit circle
that asymptotically contains no zeros.  

\begin{corollary}[Zero-free regions]
Suppose that
$\phi(\theta)=e^{-nV(\theta)}$ where $V:S^1\rightarrow\mathbb{R}$ is
of class $C^{k-1,1}(S^1)$ for some $k\ge 2$, and suppose that $\kappa'(\theta)$
is strictly positive.  Let $\delta>0$ be an arbitrarily
small number.  Then there are no zeros
of $\pi_n(z)$ in the region 
\begin{eqnarray}
\label{eq:zerofreeset}
\left\{z=re^{i \theta}  \Bigg| \log(r)>-\left(\frac{k-1-\delta}{1 +
\Omega'(\theta)}\right) \frac{\log(n)}{n} \right\} 
\end{eqnarray}
as long as $n$ is sufficiently large.
\label{cor:varyzerofree}
\end{corollary}
\begin{remark}
  In the case of a fixed weight, Corollary \ref{cor:fixedzerofree}
  established the existence of a zero-free annulus with outer radius
  1, and whose inner radius depended explicitly on $k$, the degree of
  smoothness.  We subsequently considered a family of weights for
  which we could compute explicitly the behavior of the zeros, and
  showed that Corollary \ref{cor:fixedzerofree} is sharp.  Indeed, for
  explicit families of weights with $k$ degrees of smoothness, some
  zeros achieve a distance of $k \log(n)/n + \log(\log(n))/n$ from the
  unit circle.  In addition, for these examples, a majority of the
  zeros approach a circle whose radius is $1 - (k+1) \log(n)/n$.  For
  the case of varying weights, it is our belief that Corollary
  \ref{cor:varyzerofree} is similarly sharp.  Moreover, it is to be
  expected that for some canonical family of examples constructed such
  that $\Omega^{(k)}$ has jump discontinuities, it should be similarly
  possible to obtain a very detailed asymptotic description of the
  zeros.
\end{remark}

\begin{remark}
  The zero-free region is determined by observing that for there to be
  a zero of $\pi_{n}(z)$, $z^n e^{inE_k\Omega(r,\theta)}$ must be
  roughly the same size as $\log(n)/n^{k-1}$.  This is clearly not
  true for $z$ on the unit circle, and investigating the size of $z^n
  e^{inE_k\Omega(r,\theta)}$ relative to $\log(n)/n^{k-1}$ yields the
  result.  For $|z|<1$, $z^n e^{inE_k\Omega(r,\theta)}$ is small not
  only because $z^n$ is small, but also because (as will be clear in
  \S~\ref{sec:varyingdetails} below) we constructed
  $e^{inE_k\Omega(r,\theta)}$ so that for $r<1$ but $1-r$ sufficiently
  small,
\begin{equation}
e^{inE_k\Omega(r,\theta)} = O\left(r^{n\Omega'(\theta)}\right)\,.
\end{equation}
Thus our extension of the function $\Omega$ plays a role in
determining the zero-free region near the unit circle.
\end{remark}

\begin{remark}
  The quantity $1 + \Omega'(\theta)$ appearing in Corollary
  \ref{cor:varyzerofree} is strictly positive because $1 +
  \Omega'(\theta)=\kappa'(\theta)$, and $\kappa'(\theta)$ is strictly
  positive.  It is interesting to note that this condition also
  guarantees that there are no gaps in the support of the equilibrium
  measure (see Appendix~\ref{app:potential}).  The occurrence of a gap
  in the support of the equilibrium measure is heralded by the
  development of a zero of the function $\kappa'(\theta)$.  Although
  Corollary \ref{cor:varyzerofree} would not apply if
  $\kappa'(\theta)$ vanished at some $\theta_{0}$, an intuitive
  consideration of the set defined in \eqref{eq:zerofreeset} indicates
  that near $\theta_{0}$, the zeros are pushed away, further from the
  unit circle.  As a gap develops then, one might expect the zeros to
  accumulate on a contour approaching the unit circle, but with a gap
  aligned with the gap in the support of the equilibrium measure.  We
  will not carry out such an analysis here, but clearly the methods
  outlined here can be adapted in this direction.
\end{remark}

Just as Theorem \ref{thm:piinside} leads to Corollary \ref{cor:alpha},
\eqref{eq:pivaryawaycircleinnerestimate} from Theorem \ref{thm:varyingcircleattract}
yields the following result.  (Recall from the definition \eqref{eq:negfreq} that
$N(z)\rightarrow 0$ as $z\rightarrow \infty$, and so  
$\overline{N(1/\overline{z})}  \rightarrow 0$ as $z\rightarrow 0$.)  
\begin{corollary}[Varying recurrence coefficients]
Suppose that $\phi(\theta)=e^{-nV(\theta)}$ where
$V:S^1\rightarrow\mathbb{R}$ is of class $C^{k-1,1}(S^1)$ with $k\ge
2$.  If $\kappa'(\theta)$ is strictly positive, then there is a
constant $K>0$ such that the estimate
\begin{equation}
\label{eq:varyverblunskyas}
|\alpha_n|\le \frac{K}{n^{k-1}}
\end{equation}
holds for sufficiently large $n$.
\label{cor:varyingalpha}
\end{corollary}

Finally, we have the following result concerning the asymptotic
behavior of the normalization constant $\gamma_n>0$, defined  such that
$\|\gamma_n\pi_n(z)\|_\phi=1$, with $\phi(\theta)$ given in the
varying weights case by \eqref{eq:varyingweight}.
\begin{theorem}
Suppose that $V$ is a real function of class $C^{k-1,1}(S^1)$ for some
$k\ge 2$.  If $\kappa'(\theta)$ is strictly positive, then there is a
constant $K>0$ such that
\begin{equation}
\left|\gamma_{n-1}^2e^{-nV_0}-1\right|\le \frac{K\log(n)}{n^{2(k-1)}}
\label{eq:varynormingasym}
\end{equation}
holds for all $n$ sufficiently large.
\label{thm:varyinggamma}
\end{theorem}

\begin{remark}
At first glance, the asymptotic formulae \eqref{eq:varynormingasym}
and \eqref{eq:varyverblunskyas} may appear to be inconsistent with the
identity \eqref{eq:gammaalpha}.  However, in
\eqref{eq:varynormingasym}, $\gamma_{n-1}$ is the norming coefficient
of the $(n-1)$-st degree polynomial orthogonal with respect to the
$n$-dependent weight $e^{-nV}$.  To verify \eqref{eq:gammaalpha} in
the varying weights case, one would have to compute the asymptotics
for $\gamma_{n}$ as well, rather than merely replacing $n$ by $n+1$ in
\eqref{eq:varynormingasym}.  
\end{remark}

\subsection{The $\dbar$ steepest descent method for 
exponentially varying weights.}
\label{sec:varyingdetails}
One of the main points of this paper is that while the 
limit appropriate for exponentially varying weights of the form
\eqref{eq:varyingweight} lies for the most part beyond the reach of
classical techniques applicable for fixed weights (meaning primarily
the approximation of $\phi(\theta)^{-1}$ by positive trigonometric
polynomials), analysis of this limit by means of the $\dbar$ steepest
descent method presents almost no further difficulty beyond the
analysis carried out for fixed weights in \S~\ref{sec:strong}.  To
illustrate the ease with which many of the techniques carry over from
the fixed weights case, we now outline the analogous calculations for
the varying weight
\eqref{eq:varyingweight}.

As in the case of fixed weights, we begin with the matrix ${\bf
M}^n(z)$ solving Riemann-Hilbert Problem~\ref{rhp:M} and introduce a
sequence of explicit transformations.  In order that Riemann-Hilbert
Problem~\ref{rhp:M} indeed characterizes the polynomial of degree $n$
orthogonal to all lower degree polynomials with respect to the weight
\eqref{eq:varyingweight}, we make the following assumption.
\begin{assumption}
$V$ is a real continuous function on the circle that, for some exponent $\nu\in (0,1]$ and for some constant $K>0$, satisfies a uniform H\"older continuity condition $|V(\theta_2)-V(\theta_1)|\le K|\theta_2-\theta_1|^\nu$.  
\label{ass:kappa1}
\end{assumption}
\subsubsection{Conversion to an equivalent $\dbar$ problem.  Solution of
the $\dbar$ problem in terms of integral equations.}
\label{sec:varyingconvert}
Now, define ${\bf S}^n(z)$ in terms
of ${\bf M}^n(z)$ as follows:
\begin{equation}
{\bf S}^n(z):=\left\{
\begin{array}{ll}
\displaystyle e^{nV_0\sigma_3/2}{\bf M}^n(z)e^{-nV_0\sigma_3/2}
\left(\begin{array}{cc} z^{-n}e^{-nN(z)}
& 0\\\\
0 & z^ne^{nN(z)}\end{array}\right)\,,
&\hspace{0.2 in}\text{for $|z|>1$}\\\\
\displaystyle e^{nV_0\sigma_3/2}{\bf M}^n(z)e^{-nV_0\sigma_3/2}
\left(\begin{array}{cc} e^{-n\overline{N(1/\overline{z})}}&0\\\\
0 & e^{n\overline{N(1/\overline{z})}}\end{array}\right)\,,
&\hspace{0.2 in}\text{for $|z|<1$}\,.
\end{array}
\right.
\label{eq:MS}
\end{equation}
Since $N(z)\rightarrow 0$ as $z\rightarrow\infty$, it is easy to see
that as $z\rightarrow\infty$, ${\bf S}^n(z)\rightarrow\mathbb{I}$.
Moreover, ${\bf S}^n(z)$ is analytic for $|z|\neq 1$, and its boundary
values ${\bf S}^n_+(z)$ (respectively ${\bf S}^n_-(z)$) taken on the
circle $S^1$ from the inside (respectively outside) satisfy:
\begin{equation}
{\bf S}^n_+(e^{i\theta})={\bf S}^n_-(e^{i\theta})\left(
\begin{array}{cc} e^{in\kappa(\theta)} & 1 \\\\ 0 & e^{-in\kappa(\theta)}
\end{array}\right)\,.
\label{eq:Sjump}
\end{equation}
Note that according to Assumption~\ref{ass:kappa1}, $\kappa(\theta)$
also satisfies a uniform H\"older continuity condition with exponent $\nu$.

As in \S~\ref{sec:strongconvert}, we have available an algebraic factorization of this jump
condition:
\begin{equation}
{\bf S}^n_+(e^{i\theta})={\bf S}^n_-(e^{i\theta})
\left(\begin{array}{cc}
1 & 0 \\\\
\displaystyle e^{-in\kappa(\theta)} & 1\end{array}\right)
\left(\begin{array}{cc} 0 & 1 \\\\ -1 & 0\end{array}\right)
\left(\begin{array}{cc} 1 & 0 \\\\
\displaystyle e^{in\kappa(\theta)} & 1\end{array}\right)\,.
\label{eq:factorizationvarying}
\end{equation}
To take advantage of this factorization we need to extend the
functions $e^{\pm in\kappa(\theta)}$ from $S^1$ to an annulus
containing $S^1$.  To have continuity of these extensions, we make the
following assumption.
\begin{assumption}
The function $V$ is of class $C^{k-1}(S^1)$ for some $k=1,2,3,\dots$.  
\label{ass:kappa2}
\end{assumption}
With this assumption, which provides new information only if $k>1$, we
can extend $e^{\pm in\kappa(\theta)}$ as follows.  Writing $e^{\pm
in\kappa(\theta)}=e^{\pm in\theta}e^{\pm in\Omega(\theta)}$,
we extend the factor $e^{\pm in\theta}$ analytically as $z^{\pm n}$.
On the other hand, $\Omega:S^1\rightarrow\mathbb{R}$ is a
well-defined function on the circle, and therefore we may apply the
extension operator $E_m$ defined in
\eqref{eq:extensiondefine} to $\Omega$, resulting in an 
extension of $e^{\pm in\Omega(\theta)}$ to the domain
$\mathbb{R}^2\setminus\{0\}$ as a continuous function $e^{\pm
inE_m\Omega(r,\theta)}$ for any $m$ in the range $1\le m\le
k$.

Let $\epsilon>0$ be the radius parameter of the annular domains $A_\pm$
(see Figure~\ref{fig:domains}).  Recall the ``bump'' function $B$
defined in \S~\ref{sec:notation}.
We define a new matrix unknown ${\bf
T}^n_{m,\epsilon}(r,\theta)$ by setting
\begin{equation}
{\bf T}^n_{m,\epsilon}(r,\theta):=\left\{
\begin{array}{ll}
\displaystyle {\bf S}^n(z)\left(\begin{array}{cc}
1 & 0\\\\B(\log(r)/\epsilon)r^{-n}e^{-in\theta}
e^{-inE_m\Omega(r,\theta)} & 1\end{array}\right)\,,
&\hspace{0.2 in}\text{for $z=re^{i\theta}\in A_-$}\,,
\\\\
\displaystyle
{\bf S}^n(z)\left(\begin{array}{cc}
1 & 0 \\\\ -B(\log(r)/\epsilon)r^ne^{in\theta}
e^{inE_m\Omega(r,\theta)} & 1\end{array}\right)\,, &
\hspace{0.2 in}\text{for $z=re^{i\theta}\in A_+$}\,,
\\\\
\displaystyle {\bf S}^n(z)\,,&\hspace{0.2 in}\text{for $z=re^{i\theta}
\not\in\overline{A_+\cup A_-}$}\,.
\end{array}
\right.
\label{eq:ST}
\end{equation}
Note that the presence of the factor $B(\log(r)/\epsilon)$ ensures
that ${\bf T}^n_{m,\epsilon}(r,\theta)$ may be continuously extended
to the outer boundary of $A_-$ and the inner boundary of
$A_+$; that is, ${\bf T}^n_{m,\epsilon}(r,\theta)$ is a
continuous function for $r\neq 1$.

\begin{remark}
  Here our analysis rests upon extending $\Omega(\theta)$ from the
  circle into domains $A_\pm$.  By contrast in our treatment of the
  fixed weights case in \S~\ref{sec:dbarmethodfixed} it was more
  convenient to extend the function $V(\theta)$, which is related to
  $\Omega(\theta)$ by a Cauchy transform.
\end{remark}

If $m=1$, then it is easy to see that 
$E_m\Omega(r,\theta)=E_1\Omega(r,\theta)=
\Omega(\theta)$ for all $r>0$, so the off-diagonal matrix elements in 
${\bf S}^n(re^{i\theta})^{-1}{\bf T}_{1,\epsilon}^n(r,\theta)$ are bounded in 
magnitude by $e^{-n|\log(r)|}$. Therefore, in this case $\eqref{eq:ST}$ 
represents an exponentially near-identity
transformation when $n$ is large and $|\log(r)|$ is not too small.  We
want to ensure that a similar situation prevails when $m\ge 2$ as
well.  Since $1\le m\le k$, the consideration of $m\ge 2$ requires
that $k\ge 2$.  In this case Assumption~\ref{ass:kappa2} implies that
$\kappa(\theta)$ is continuously differentiable, and hence so is
$\Omega(\theta)$.  The crucial conditions for \eqref{eq:ST} to
be a near-identity transformation when $n$ is large are that
$\Omega'(\theta)+1$ is strictly positive and that the
parameter $\epsilon$ is chosen small enough, as the following lemma
shows.
\begin{lemma}
Suppose that $\Omega(\theta)$ is a real function of class
$C^{m-1}(S^1)$ for some $m\ge 2$, such that $\Omega'(\theta)+1$
is strictly positive.  Then there exist constants $\epsilon_0>0$ 
and $\mu>0$ such that whenever $0<\epsilon<\epsilon_0$,
\begin{equation}
\left|e^{-iE_m\Omega(r,\theta)}\right|
\le r^{1-\mu}\,,\hspace{0.2 in} 
\text{for all $\theta$ and for $1\le r \le 2^\epsilon$}\,,
\label{eq:impartboundm}
\end{equation}
and
\begin{equation}
\left|e^{iE_m\Omega(r,\theta)}\right|
\le r^{\mu-1}\,,\hspace{0.2 in} 
\text{for all $\theta$ and for $2^{-\epsilon}\le r\le 1$}\,.
\label{eq:impartboundp}
\end{equation}
\label{lemma:exponentcontrol}
\end{lemma}
\begin{proof}
By hypothesis, we have
\begin{equation}
\inf_{-\pi<\theta<\pi}\Omega'(\theta) = \tau-1\,,\hspace{0.2 in}
\text{where $\tau>0$}\,.
\end{equation}
Note that $\log|e^{\mp iE_m\Omega(r,\theta)}| = 
\pm\Im(E_m\Omega(r,\theta))$.  
From \eqref{eq:extensiondefine}, we have
\begin{equation}
\Im(E_m\Omega(r,\theta))=\sum_{p=0}^{P(m)}
(-1)^{p+1}\frac{\Omega^{(2p+1)}(\theta)}{(2p+1)!}\log(r)^{2p+1} = 
-\log(r)\left[\Omega'(\theta) +
\sum_{p=1}^{P(m)}(-1)^{p}\frac{\Omega^{(2p+1)}(\theta)}{(2p+1)!}
\log(r)^{2p}\right]\,,
\end{equation}
where $P(m)=(m-2)/2$ if $m$ is even and $P(m)=(m-3)/2$ if $m$ is odd.
Since $\Omega$ has $m-1$ continuous derivatives we have
\begin{equation}
\left|\sum_{p=1}^{P(m)}(-1)^p\frac{\Omega^{(2p+1)}(\theta)}{(2p+1)!}
\log(r)^{2p}\right|\le
\sum_{p=1}^{P(m)}
\frac{\displaystyle \sup_{-\pi<\theta<\pi}|\Omega^{(2p+1)}(\theta)|}{(2p+1)!}
(\epsilon\log(2))^{2p}
\end{equation}
for all $r$ satisfying $|\log(r)|\le \epsilon\log(2)$ where
$\epsilon>0$.  Therefore, $\epsilon_0>0$ may be chosen sufficiently
small that the inequality $0<\epsilon<\epsilon_0$ implies that
\begin{equation}
\Omega'(\theta) + \sum_{p=1}^{P(m)}
(-1)^p\frac{\Omega^{(2p+1)}(\theta)}{(2p+1)!}\log(r)^{2p}\ge
\frac{\tau}{2}-1
\end{equation}
holds whenever $|\log(r)|\le \epsilon\log(2)$.  Then, we have
\begin{equation}
\log\left|e^{-iE_m\Omega(r,\theta)}\right|\le 
-\log(r)\left[\frac{\tau}{2}-1\right]\,,\hspace{0.2 in}
\text{holding for $1\le r\le 2^\epsilon$}
\end{equation}
and
\begin{equation}
\log\left|e^{iE_m\Omega(r,\theta)}\right|\le 
\log(r)\left[\frac{\tau}{2}-1\right]\,,\hspace{0.2 in}
\text{holding for $2^{-\epsilon}\le r\le 1$}\,.
\end{equation}
The estimates
\eqref{eq:impartboundm} and \eqref{eq:impartboundp} thus both hold
with the choice $\mu:=\tau/2>0$.
\end{proof}

Thus, we are led to propose the following assumption.
\begin{assumption}
If $V$ is of class $C^{1}(S^1)$, then
$\kappa'(\theta)$ defined by \eqref{eq:kappadefine} is strictly positive.
\label{ass:kappa3}
\end{assumption}

This assumption is not an explicit assumption on $V(\theta)$, however it is possible to give conditions on $V(\theta)$ that are sufficient to make Assumption~\ref{ass:kappa3} hold.  For example, since for $V\in C^1(S^1)$,
we have
\begin{equation}
\kappa'(\theta)=1-2\sum_{j=1}^\infty
j|V_j|\cos(j\theta+\arg(V_j))\,,
\end{equation}
where we recall the Fourier coefficients of $V$ defined by
\eqref{eq:FourierSeriesV}, the condition
\begin{equation}
\sum_{j=-\infty}^\infty |j||V_j|< 1\,,
\label{eq:tcondition}
\end{equation}
is sufficient to guarantee that Assumption~\ref{ass:kappa3} holds.
If in fact $V\in C^2(S^1)$, then the convexity condition
\begin{equation}
\frac{d^2V}{d\theta^2}(\theta)>-\frac{1}{2}\,,
\label{eq:convexity}
\end{equation}
also guarantees that $\kappa'(\theta)$ is strictly positive.  Our
proof that the condition \eqref{eq:convexity} implies that
$\kappa'(\theta)$ is strictly positive makes use of certain aspects of
logarithmic potential theory and is given in
Appendix~\ref{app:potential}.  In Appendix~\ref{app:compare}, we show
that neither are the sufficient conditions (\ref{eq:tcondition}) and
(\ref{eq:convexity}) equivalent, nor does either condition imply the
other.  Rather, the two conditions are independent and thus complement
each other.

\begin{remark}
In the approach taken in \cite{tail}, where a particular case of a
varying weight in which the function $V(\theta)$, and hence
$\kappa(\theta)$, is an analytic function, the analytic extension
$E_\infty\Omega(r,\theta)$ is used, and the ``bump'' function
factor $B(\log(r)/\epsilon)$ is omitted.  The latter omission has the
effect of introducing exponentially small jump discontinuities in
${\bf T}^n_{m,\epsilon}(r,\theta)$ across the circles $\Sigma_\pm$.
\end{remark}

Since for $m=1,\dots,k$ we have 
$E_m\Omega(1,\theta)=\Omega(\theta)$, the matrix ${\bf
T}^n_m(r,\theta)$ satisfies the jump condition
\begin{equation}
\lim_{r\uparrow 1}{\bf T}^n_{m,\epsilon}(r,\theta)=
\lim_{r\downarrow 1}{\bf T}^n_{m,\epsilon}(r,\theta)\left(\begin{array}{cc}
0 & 1 \\\\ -1 & 0\end{array}\right)\,.
\end{equation}
To remove this jump discontinuity, introduce one further change of variables:
\begin{equation}
{\bf J}^n_{m,\epsilon}(r,\theta):=
\left\{\begin{array}{ll}
\displaystyle 
{\bf T}^n_{m,\epsilon}(r,\theta)\,,&\hspace{0.2 in} \text{for $r>1$}\\\\
\displaystyle {\bf T}^n_{m,\epsilon}(r,\theta)\left(\begin{array}{cc}
0 & -1 \\\\ 1 & 0\end{array}\right)\,,&\hspace{0.2 in}
\text{for $r<1$}\,.
\end{array}
\right.
\label{eq:JT}
\end{equation}

At this point, we can relate ${\bf J}^n_{m,\epsilon}(r,\theta)$ directly to 
${\bf M}^n(z)$.  Combining \eqref{eq:MS}, \eqref{eq:ST}, and \eqref{eq:JT},
we have by definition
\begin{equation}
{\bf J}^n_{m,\epsilon}(r,\theta):=\left\{\begin{array}{l}
\displaystyle e^{nV_0\sigma_3/2}
{\bf M}^n(z)e^{-nV_0\sigma_3/2}\left(\begin{array}{cc}
0 & -e^{-n\overline{N(1/\overline{z})}}\\\\
e^{n\overline{N(1/\overline{z})}} & e^{n\overline{N(1/\overline{z})}}
B(\log(r)/\epsilon)r^ne^{in\theta}
e^{inE_m\Omega(r,\theta)}
\end{array}\right)\,, \\\\ \hspace{4 in}
\text{for $0\le r <1$}\\\\
\displaystyle
e^{nV_0\sigma_3/2}
{\bf M}^n(z)e^{-nV_0\sigma_3/2}\left(\begin{array}{cc}
e^{-nN(z)}r^{-n}e^{-in\theta} &
0 \\\\
e^{nN(z)}B(\log(r)/\epsilon)e^{-inE_m\Omega(r,\theta)} & 
e^{nN(z)}r^ne^{in\theta}\end{array}\right)\,,\\\\\hspace{4 in}
\text{for $r>1$}\,.
\end{array}
\right.
\label{eq:JM}
\end{equation}
Given the assumptions in force, this matrix function is clearly continuous
throughout the plane.  To determine its deviation from being an analytic 
function in the regions $r<1$ and $r>1$, we need to control a derivative,
and consequently we make the following assumption.
\begin{assumption}
The function $V$ is of class $C^{k-1,1}(S^1)$ for some $k\ge 1$.  That is, 
$V^{(k-1)}(\theta)$ is Lipschitz continuous.
\label{ass:kappa4}
\end{assumption}
Then, for $1\le m\le k$, we see by direct calculation that
\begin{equation}
\dbar{\bf J}^n_{m,\epsilon}(r,\theta)=\left\{\begin{array}{l}
\displaystyle e^{nV_0\sigma_3/2}
{\bf M}^n(z)e^{-nV_0\sigma_3/2}\left(\begin{array}{cc}
0 & 0 \\\\
0 & e^{n\overline{N(1/\overline{z})}}r^ne^{in\theta}
\dbar\left[B(\log(r)/\epsilon)
e^{inE_m\Omega(r,\theta)}\right]
\end{array}\right)\,, \\\\ \hspace{4 in}
\text{for $0\le r <1$}\\\\
\displaystyle
e^{nV_0\sigma_3/2}
{\bf M}^n(z)e^{-nV_0\sigma_3/2}\left(\begin{array}{cc}
0 &
0 \\\\
e^{nN(z)}\dbar\left[B(\log(r)/\epsilon)
e^{-inE_m\Omega(r,\theta)}\right] & 0 \end{array}\right)\,,\\\\\hspace{4 in}
\text{for $r>1$}\,.
\end{array}
\right.
\label{eq:dbarJ}
\end{equation}
Eliminating ${\bf M}^n(z)$ in favor of ${\bf J}^n_{m,\epsilon}(r,\theta)$ 
using \eqref{eq:JM} again yields
\begin{equation}
\dbar{\bf J}^n_{m,\epsilon}(r,\theta) = 
{\bf J}^n_{m,\epsilon}(r,\theta){\bf X}^n_{m,\epsilon}(r,\theta)
\,,\hspace{0.2 in}\text{for $r\neq 1$ and almost all $\theta\in S^1$}\,,
\label{eq:dbarrelationJ}
\end{equation}
where
\begin{equation}
{\bf X}^n_{m,\epsilon}(r,\theta):=\left\{
\begin{array}{ll}
\displaystyle\left(\begin{array}{cc}
0 & r^ne^{in\theta}\dbar
\left[B(\log(r)/\epsilon)e^{inE_m\Omega(r,\theta)}\right]\\\\
0 & 0 \end{array}\right)\,,&\hspace{0.2 in}\text{for $0\le r<1$}\\\\
\displaystyle\left(\begin{array}{cc}
0 & 0  \\\\
r^{-n}e^{-in\theta}\dbar
\left[B(\log(r)/\epsilon)e^{-inE_m\Omega(r,\theta)}\right] & 0 \end{array}
\right)\,,&\hspace{0.2 in}\text{for $r>1$}\,.
\end{array}
\right.
\label{eq:Xdef}
\end{equation}
Then it is easy to see that ${\bf J}^n_{m,\epsilon}(r,\theta)$ satisfies the following
$\dbar$ problem.
\begin{dbp}
Find a $2\times 2$ matrix ${\bf U}(r,\theta)$ with the properties:
\begin{itemize}
\item[]{\bf Smoothness.}  ${\bf U}(r,\theta)$ is a Lipschitz continuous
function throughout $\mathbb{R}^2$.
\item[]{\bf Deviation From Analyticity.}
The relation
\begin{equation}
\dbar{\bf U}(r,\theta)={\bf U}(r,\theta){\bf X}^n_{m,\epsilon}(r,\theta)
\end{equation}
holds for all points in $\mathbb{R}^2$ with the exception of a set
of measure zero.  The matrix ${\bf X}^n_{m,\epsilon}(r,\theta)$ is defined almost everywhere by \eqref{eq:Xdef} and is essentially compactly supported.
\item[]{\bf Normalization.}
The matrix ${\bf U}(r,\theta)$ is normalized at $r=\infty$ as
follows:
\begin{equation}
\lim_{r\rightarrow\infty}{\bf U}(r,\theta)={\mathbb I}\,.
\label{eq:Unorm}
\end{equation}
\end{itemize}
\label{dbp:U}
\end{dbp}
In exactly the same way that Proposition~\ref{prop:dbarEsolve} was proved,
we have the following.
\begin{prop}
Suppose that $\phi=e^{-nV}$ where $V:S^1\rightarrow\mathbb{R}$ is of
class $C^{k-1,1}(S^1)$ for some $k=1,2,3,\dots$.  Then for all $n=0,1,2,3,\dots$, for $m=1,2,\dots,k$, and for all $\epsilon>0$, 
the matrix ${\bf X}_{m,\epsilon}^n(r,\theta)$ is well-defined almost everywhere by 
\eqref{eq:Xdef} and $\dbar$ Problem~\ref{dbp:U} has a unique solution,
namely ${\bf U}(r,\theta)={\bf J}^n_{m,\epsilon}(r,\theta)$.
\end{prop}
Also, the proof of Proposition~\ref{prop:dbarEinteqn} carries over to the
context of the varying weight \eqref{eq:varyingweight} in the following
form.
\begin{prop}
Suppose that $\phi=e^{-nV}$ where $V:S^1\rightarrow\mathbb{R}$ is of class $C^{k-1,1}(S^1)$ for some $k=1,2,3,\dots$.  Then for all $n=0,1,2,3,\dots$, for $m=1,2,\dots,k$, and for all $\epsilon>0$, 
the matrix ${\bf X}_{m,\epsilon}^n(r,\theta)$ is well-defined almost everywhere
by \eqref{eq:Xdef} and the corresponding solution ${\bf U}(r,\theta)=
{\bf J}^n_{m,\epsilon}(r,\theta)$ of $\dbar$ Problem~\ref{dbp:U} satisfies the
integral equation
\begin{equation}
{\bf U}(r,\theta)=\mathbb{I}-\frac{1}{\pi}\int\int
\frac{{\bf U}(r',\theta'){\bf X}_{m,\epsilon}^n(r',\theta')}{z'-z}\,dA'
\label{eq:inteqnU}
\end{equation}
where $z=re^{i\theta}$, $z'=r'e^{i\theta'}$, and $dA'$ is a positive area
element $dA'=r'\,dr'\,d\theta'$.  The integral is taken over the entire
plane.
\label{prop:inteqnU}
\end{prop}

\subsubsection{Asymptotic solution of the integral equation.  Estimates
of ${\bf J}_{m,\epsilon}^n(r,\theta)$ and its derivatives for large
$n$.}  As in \S~\ref{sec:stronginteqn}, it is possible to characterize
${\bf J}_{m,\epsilon}^n(r,\theta)$ by analyzing the integral equation
\eqref{eq:inteqnU} as long as $n$ is large enough.  However, in the
context of the varying weight \eqref{eq:varyingweight} we will require
the monotonicity condition expressed in Assumption~\ref{ass:kappa3},
and that the radius parameter $\epsilon$ be taken sufficiently small
for each admissable given $V$ that the conclusion of
Lemma~\ref{lemma:exponentcontrol} holds.  Moreover, the exponential
character of the varying weight
\eqref{eq:varyingweight} suggests that, by comparison with the
analysis presented in \S~\ref{sec:stronginteqn}, this approach is only
fruitful in giving the same degree of control on ${\bf
J}^n_{m,\epsilon}(r,\theta)$ as was achieved for ${\bf
H}_{m,\epsilon}^n(r,\theta)$ if $V(\theta)$ has more smoothness than
was required in Proposition~\ref{prop:Hbound}.  In particular, we will
require that the following assumption holds.
\begin{assumption}
The function $V$ is of class $C^{k-1,1}(S^1)$ for some $k\ge 2$.
\label{ass:kappa5}
\end{assumption}
In other words, to achieve the same convergence rates as in the fixed
weight case, one more derivative of $V$ will be required in the
exponentially varying weight case.

As in \S~\ref{sec:stronginteqn}, we need some bounds for ${\bf X}^n_{m,\epsilon}(r,\theta)$ and its derivatives.
\begin{prop}
Suppose that $V:S^1\rightarrow\mathbb{R}$ is a real function of class
$C^{k-1,1}(S^1)$ for some $k\ge 1$ for which
$\Omega'(\theta)+1$ is strictly positive, that $m$ is an
integer satisfying $1\le m\le k$, and that $\epsilon>0$ is
sufficiently small.  Let the integer $D$ be defined as
$D:=\min(k-m,m-1)$.  Then, the matrix function ${\bf
X}^n_{m,\epsilon}$ is of class
$C_0^{D-1,1}(\mathbb{R}^2\setminus\{0\})$ if $D>0$, and of class
$L_0^\infty(\mathbb{R}^2\setminus\{0\})$ if $D=0$.  Moreover, if
$\alpha$ and $\beta$ are nonnegative integers such that
$\alpha+\beta\le D$, then there are constants $\mu>0$ and
$C_{m,\epsilon}^{(\alpha,\beta)}>0$ such that for all $n$ the estimate
\begin{equation}
\left\|\frac{\partial^{\alpha+\beta}}{\partial r^\alpha\partial\theta^\beta}
{\bf X}^n_{m,\epsilon}(r,\theta)\right\|\le C_{m,\epsilon}^{(\alpha,\beta)}n^{1+\beta}
e^{-n\mu|\log(r)|}|\log(r)|^{m-1-\alpha}\sum_{p=0}^\alpha n^p|\log(r)|^p
\label{eq:derivXbound}
\end{equation}
holds throughout the region $|\log(r)|\le \epsilon\log(2)$ containing
the essential support of ${\bf X}^n_{m,\epsilon}(r,\theta)$.
\label{prop:Xbound}
\end{prop}

\begin{proof}
This proposition is proved in almost the same way as
Proposition~\ref{prop:Wbound}.  In this case, it is essential to
recall Lemma~\ref{lemma:exponentcontrol} which provides the constant
$\mu>0$ and thus the exponential decay of the term
$e^{-n\mu|\log(r)|}$ for $\epsilon>0$ sufficiently small.  Also, the
initial application of the $\dbar$ operator in the definition
\eqref{eq:Xdef} yields a factor of $n$ that does not appear in the
proof of Proposition~\ref{prop:Wbound}.
\end{proof}

\begin{remark}
This result should be compared with Proposition~\ref{prop:Wbound}.  The only important difference between the estimate \eqref{eq:derivXbound} and
\eqref{eq:derivWbound} is the presence of an additional factor of $n$.
\end{remark}

The proof of Proposition~\ref{prop:Hbound}, with references to
Proposition~\ref{prop:Wbound} replaced by references to
Proposition~\ref{prop:Xbound}, applies to the asymptotic estimation
of ${\bf J}_{m,\epsilon}^n(r,\theta)$, with the only important difference being
an additional factor of $n$.  This results in the following.
\begin{prop}
\label{prop:Jbound} 
Suppose that $\phi=e^{-nV}$ where $V:S^1\rightarrow\mathbb{R}$ is of
class $C^{k-1,1}(S^1)$ for some $k=2,3,4,\dots$.  Let the integer $m$
lie in the range $2\le m\le k$ and fix $\epsilon>0$ sufficiently
small. Define the integer $\tilde{D}:=\min(k-m,m-2)\ge 0$.  Then, for
all $n\ge 0$ the matrix ${\bf X}^n_{m,\epsilon}(r,\theta)$ is
well-defined almost everywhere by
\eqref{eq:Xdef}, and for all $n$ sufficiently large, ${\bf
J}^n_{m,\epsilon}(r,\theta)$ is given by a 
Neumann series
\begin{equation}
{\bf J}_{m,\epsilon}^n(r,\theta)=\mathbb{I} +
(\mathcal{X}_{m,\epsilon}^n\mathbb{I})(r,\theta) +
(\mathcal{X}_{m,\epsilon}^n\circ\mathcal{X}_{m,\epsilon}^n\mathbb{I})(r,\theta)
+ \cdots
\end{equation}
which converges in the norm $\||\cdot|\|_{\tilde D}$, where the double integral
operator $\mathcal{X}_{m,\epsilon}^n$ is defined by
\begin{equation}
(\mathcal{X}^n_{m,\epsilon}{\bf
F})(r,\theta):=-\frac{1}{\pi}\int\int\frac{{\bf F}(r',\theta') {\bf
X}^n_{m,\epsilon}(r',\theta')}{z'-z}\,dA'\,.
\label{eq:varyingCnmdef}
\end{equation}
In particular, if $\tilde{D}=0$ then ${\bf J}^n_{m,\epsilon}$ lies in
the space $L^\infty(\mathbb{R}^2)$, and if $\tilde{D}>0$ then ${\bf
J}^n_{m,\epsilon}$ lies in the space $C^{\tilde{D}-1,1}(\mathbb{R}^2)$
and $\||{\bf J}^n_{m,\epsilon}|\|_{\tilde{D}}$ is finite.  For all
integer $p$ in the range $0\le p\le \tilde{D}$, the following
estimates hold for sufficiently large $n$:
\begin{equation}
\||{\bf J}_{m,\epsilon}^n-\mathbb{I}|\|_p\le C_{m,\epsilon}^{(p)}
\frac{\log(n)}{n^{m-p-1}}\,,
\label{eq:varyingHmI}
\end{equation}
\begin{equation}
\||{\bf J}_{m,\epsilon}^n-\mathbb{I}-\mathcal{X}_{m,\epsilon}^n\mathbb{I}
|\|_p\le C_{m,\epsilon}^{(p)2}\frac{\log(n)^2}{n^{2m-2p-2}}\,,
\label{eq:varyingHmImWI}
\end{equation}
where $C_{m,\epsilon}^{(p)}$ is a positive constant.  Furthermore, for
each $\rho>2^\epsilon$ and for all integer $p$ in the range $0\le p\le
D$, the following estimates hold for sufficiently large $n$:
\begin{equation}
\sum_{\alpha+\beta\le p}
\mathop{\sup_{-\pi<\theta<\pi}}_{|\log(r)|\ge\log(\rho)}
\left\|\frac{\partial^{\alpha+\beta}}{\partial x^\alpha\partial y^\beta}
\left[{\bf J}^n_{m,\epsilon}(r,\theta)-\mathbb{I}\right]\right\|
\le \tilde{C}^{(p)}_{m,\rho}\frac{1}{n^{m-p-1}}\,,
\label{eq:varyingHmIaway}
\end{equation}
\begin{equation}
\sum_{\alpha+\beta\le p}
\mathop{\sup_{-\pi<\theta<\pi}}_{|\log(r)|\ge\log(\rho)}
\left\|\frac{\partial^{\alpha+\beta}}{\partial x^\alpha\partial y^\beta}
\left[{\bf J}^n_{m,\epsilon}(r,\theta)-\mathbb{I}-(\mathcal{X}_{m,\epsilon}^n\mathbb{I})(r,\theta)\right]\right\|
\le \tilde{C}^{(p)2}_{m,\rho}\frac{\log(n)}{n^{2m-2p-2}}\,,
\label{eq:varyingHmImWIaway}
\end{equation}
where $\tilde{C}^{(p)}_{m,\rho}$ is a positive constant.
\end{prop}

\begin{remark}
  Note that in \eqref{eq:varyingHmI}, due to the extra factor of $n$
  introduced into the estimates by Proposition~\ref{prop:Xbound},
  convergence of the Neumann series in the $\||\cdot|\|_p$ norm
  follows provided that $m-p-1>0$, and so we may only consider $m\ge
  2$.  In order to uniformly control $p$ derivatives of ${\bf
    J}_{m,\epsilon}^n(r,\theta)$, Proposition~\ref{prop:Jbound}
  requires that $m$ should lie in the range $2+p\le m\le k-p$, and
  therefore in order for there to exist suitable values of $m$,
  $V:S^1\rightarrow\mathbb{R}$ should be of class $C^{k-1,1}(S^1)$ for
  some $k\ge 2p+2$.
\end{remark}

\begin{remark}
  Assumption~\ref{ass:kappa3}, that $V(\theta)$ is such that
  $\kappa(\theta)$ is a strictly increasing function of the angle
  $\theta$, ties together two crucial aspects of our analysis.  First
  of all, as can be seen from the matrix factors involved in the
  change of variables \eqref{eq:ST} between ${\bf S}^n(z)$ and ${\bf
    T}^n_{m,\epsilon}(r,\theta)$, the inequality $\kappa'(\theta)>0$
  is precisely what makes this a near-identity change of variables in
  the regions $\Omega_\pm$.  On the other hand, it would not suffice
  to replace $\kappa'(\theta)$ by another unrelated positive quantity,
  because we require of whatever extension we choose some degree of
  vanishing of the $\dbar$ derivative at $|z|=1$ in order to control
  the $\dbar$ problem (that is, to sufficiently bound the operator
  $\mathcal{X}_{m,\epsilon}^n$).  This vanishing is built into the
  extension operators we have defined in \eqref{eq:extensiondefine} by
  the key property \eqref{eq:dbarTaylor}.  Moreover, {\em any} smooth
  extension of $\kappa(\theta)$ that has a vanishing $\dbar$
  derivative on the unit circle will behave similarly near the unit
  circle.  So we conclude that while the analyticity of
  $\kappa(\theta)$ is not important, the monotonicity of this function
  is crucial.  To make an analogy with the asymptotic analysis of
  oscillatory exponential integrals, the $\dbar$ steepest descent
  method resembles Kelvin's method of stationary phase more than it
  does the saddle point method.  On the other hand, it is important to
  note that the $\dbar$ steepest descent method remains fundamentally
  a method of deformation into the complex plane, and is not based on
  integration by parts\footnote{Properly speaking, we do not rely on
    integration by parts (or more generally, Stokes' Theorem for the
    $\dbar$ operator in the plane --- in as much as this can be
    considered a generalization of the standard $\dbar$ inversion
    formula) to establish the existence of an asymptotic expansion.
    However, such methods are useful in the detailed analysis of
    individual terms in the expansion.  This technique was used, for
    example, in the proof of Theorem~\ref{thm:piinsidejumps}.}, the
  Riemann-Lebesgue Lemma, or related arguments of harmonic analysis.
  For an approach based on the latter, see Varzugin \cite{varzugin}.
\end{remark}

\subsection{Proofs of theorems stated in \S~\ref{sec:asymptoticspids}.}
The proofs are generally based on expressing $\pi_n(z)=M_{11}^n(z)$ in terms 
of explicit functions and ${\bf J}_{m,\epsilon}^n(r,\theta)$ by means of 
\eqref{eq:JM}.  Then, one applies
Proposition~\ref{prop:Jbound} to control
${\bf J}_{m,\epsilon}^n(r,\theta)-\mathbb{I}$ and its derivatives.

\subsubsection{Asymptotic behavior of $\pi_n(z)$ for $|z|\ge 1$.
Proof of Theorems~\ref{thm:varypiasympoutbasic} and
\ref{thm:picirclevaryingweights}.} 
From \eqref{eq:JM} we have the following exact representation for $\pi_n(z)$ 
valid for $|z|>1$:
\begin{equation}
\pi_n(z)=e^{nN(z)}\left[z^nJ_{m,\epsilon,11}^n(r,\theta)-B(\log(r)/\epsilon)
e^{-inE_m\Omega(r,\theta)}J_{m,\epsilon,12}^n(r,\theta)\right]\,,
\hspace{0.2 in}|z|>1\,.
\end{equation}
Here $\epsilon$ should be taken to be sufficiently small.  It follows that
\begin{equation}
\pi_n(z)z^{-n}e^{-nN(z)}-1 = \left[J_{m,\epsilon,11}^n(r,\theta)-1\right]
-B(\log(r)/\epsilon)z^{-n}e^{-inE_m\Omega(r,\theta)}
J_{m,\epsilon,12}^n(r,\theta)\,,\hspace{0.2 in}|z|>1\,.
\end{equation}
If $\rho>1$ is fixed, then perhaps by making $\epsilon>0$ smaller yet,
it can be arranged that $B(\log(r)/\epsilon)\equiv 0$ whenever
$|z|\ge\rho$. 
In this case, we have
\begin{equation}
\frac{d^p}{dz^p}\left[\pi_n(z)z^{-n}e^{-nN(z)}\right] = \partial^p
\left[J_{m,\epsilon,11}^n(r,\theta)-1\right]\,,\hspace{0.2 in}
|z|\ge\rho>1\,.
\end{equation}
Using the estimate \eqref{eq:varyingHmImWIaway} from 
Proposition~\ref{prop:Jbound} in the case $m=k-p$, and the fact that 
${\bf X}^n_{m,\epsilon}(r,\theta)$ is off-diagonal, we see that for some 
constant $C>0$,
\begin{equation}
\mathop{\sup_{-\pi<\theta<\pi}}_{r\ge\rho>1}\left|\partial^p\left[
J_{k-p,\epsilon,11}^n(r,\theta)-1\right]\right|\le 
C\frac{\log(n)}{n^{2(k-2p-1)}}\,.
\label{eq:varyJ11outest}
\end{equation}
Now since $J^{n}_{k-p, \epsilon, 11}(z) - 1 $ is analytic for $|z|>\rho$, 
and tends to zero like $z^{-1}$ as $z \rightarrow \infty$, Cauchy's Theorem 
for an exterior domain, together with \eqref{eq:varyJ11outest} for $p=0$, 
proves \eqref{eq:varypiasympout2}.

To prove \eqref{eq:varyingpiasympcircle}, we fix $\epsilon$ sufficiently small 
and consider $1\le r\le 2^{\epsilon/2}$ in which case 
$B(\log(r)/\epsilon)\equiv 1$ so that
\begin{equation}
\frac{d^p}{dz^p}\left[\pi_n(z)z^{-n}e^{-nN(z)}-1\right] =
\partial^p\left[J_{m,\epsilon,11}^n(r,\theta)-1\right]-
\partial^p\left[z^{-n}e^{-inE_m\Omega(r,\theta)}J_{m,\epsilon,12}^n(r,\theta)
\right]\,,
\hspace{0.2 in} 1\le r\le 2^{\epsilon/2}\,.
\label{eq:intermediate3}
\end{equation}
Again, since ${\bf X}^n_{m,\epsilon}(r,\theta)$ is an off-diagonal matrix,
we see from \eqref{eq:varyingHmImWI} that there is a constant $C>0$ such that
\begin{equation}
\sup_{\mathbb{R}^2}\left|\partial^p\left[J_{k-p,\epsilon,11}^n(r,\theta)-1\right]\right|\le C\frac{\log(n)^2}{n^{2(k-2p-1)}}\,.
\end{equation}
As in the fixed weights case, the dominant contribution comes from the remaining terms on the right-hand side of \eqref{eq:intermediate3}.  
Using \eqref{eq:varyingHmI} from Proposition~\ref{prop:Jbound} to
see that $n^j\partial^{p-j}J_{m,\epsilon,12}^n(r,\theta)$ is of order
$\log(n)/n^{m-p-1}$ and taking the best case of $m=k-p$, we then find that
\begin{equation}
\mathop{\sup_{-\pi<\theta<\pi}}_{1\le r\le 2^{\epsilon/2}}
\left|\frac{d^p}{dz^p}\left[\pi_n(z)z^{-n}e^{-nN(z)}-1\right]\right|
\le K\frac{\log(n)}{n^{k-2p-1}}\,.
\end{equation}
This proves \eqref{eq:varyingpiasympcircle}.  
Now the maximum modulus principle applied to 
$\pi_{n}(z) z^{-n} e^{-n N(z)}$ implies \eqref{eq:varyingpiasympoutside}.

\subsubsection{Asymptotic behavior of $\pi_n(z)$ for $|z|<1$ and of 
$\gamma_{n-1}$.  Proof
of Theorem~\ref{thm:varyingcircleattract} and
Theorem~\ref{thm:varyinggamma}.}   
\label{sec:varyingbehaviorinside}
The proof of Theorem~\ref{thm:varyinggamma} is based on the identity
$\gamma_{n-1}^2 = -M_{21}^n(0)$.  Using \eqref{eq:JM}, we see that
\begin{equation}
\gamma_{n-1}^2e^{-nV_0} = J_{m,\epsilon,22}^n(0,\theta)\,,
\end{equation}
and thus Theorem~\ref{thm:varyinggamma} is proved by applying the estimate
\eqref{eq:varyingHmImWIaway} from Proposition~\ref{prop:Jbound} in the case of $p=0$ and $m=k$.

Proving Theorem~\ref{thm:varyingcircleattract} begins with the exact formula
\begin{equation}
\pi_n(z)=-e^{n\overline{N(1/\overline{z})}}J_{m,\epsilon,12}^n(r,\theta)
\,,\hspace{0.2 in} 0\le r\le 2^{-\epsilon}\,,
\end{equation}
which follows from \eqref{eq:JM} using $\pi_n(z)=M_{11}^n(z)$ and
the fact that $B(\log(r)/\epsilon)\equiv 0$ for $r<2^{-\epsilon}$.
Using \eqref{eq:varyingHmIaway} from Proposition~\ref{prop:Jbound} 
with $p=0$ and $m=k$ then completes the proof of
\eqref{eq:pivaryawaycircleinnerestimate}.

For the proof of \eqref{eq:pivaryingarcircinsideestimate}, one begins
with the following formula for $\pi_{n}$,
\begin{eqnarray}
\pi_{n}(z) =
-e^{n\overline{N(1/\overline{z})}}J_{m,\epsilon,12}^n(r,\theta)
+ e^{n\overline{N(1/\overline{z})}} e^{- i n E_{k}\Omega(r,
  \theta)} J_{m,\epsilon,11}^n(r,\theta),
\end{eqnarray}
which again follows from \eqref{eq:JM}, and the fact that
$B(\log(r)/\epsilon)\equiv 1$  for $r>2^{-\epsilon/2}$.
Using \eqref{eq:varyingHmI} from Proposition~\ref{prop:Jbound} 
with $p=0$ and $m=k$ then completes the proof of
\eqref{eq:pivaryingarcircinsideestimate}, and this completes the proof
of Theorem~\ref{thm:varyingcircleattract}.

\section{Acknowledgements}
During our preparation of this work, we benefited from useful
conversations with Percy Deift, Andrei Martinez-Finkelshtein, Paul
Nevai, and Walter van Assche.  We also thank Barry Simon for letting
us see a preliminary version of his forthcoming monograph
\cite{simon}.

The research of K. T.-R. McLaughlin was supported in part by the National
Science Foundation under grants DMS-9970328 and DMS-0200749.
The work of P. D. Miller was also supported in part by the 
National Science Foundation, under grant DMS-0103909, and by a grant from the 
Alfred P. Sloan Foundation.

\appendix
\section{Logarithmic Potential Theory of Orthogonal Polynomials 
on the Unit Circle}
\label{app:potential}
A more general and systematic strategy for extracting asymptotics of
${\bf M}^n(z)$ in the exponentially varying weight case when
$\phi(\theta)$ is of the form \eqref{eq:varyingweight} is the
following.  First introduce a function $g(z)$ (to be determined)
satisfying $g(z)\sim \log(z)$ as $z\rightarrow\infty$ and that
$e^{g(z)}$ is analytic for $|z|\neq 1$, taking continuous boundary
values on the unit circle $\Sigma$.  Then one converts Riemann-Hilbert
Problem~\ref{rhp:M} for ${\bf M}^n(z)$ into one with identity
asymptotics as $z\rightarrow\infty$ by the change of variables
\begin{equation}
{\bf M}^n(z)={\bf
Y}^n(z)e^{ng(z)\sigma_3}\,,
\end{equation}
resulting in a new unknown matrix ${\bf Y}^n(z)$.  In addition
to the normalization condition $\lim_{z\rightarrow\infty}{\bf
Y}^n(z)={\mathbb I}$, ${\bf Y}^n(z)$ satisfies the jump condition
\begin{equation}
{\bf Y}^n_+(e^{i\theta})={\bf
Y}^n_-(e^{i\theta})\left(\begin{array}{cc}
e^{-n(g_+(e^{i\theta})-g_-(e^{i\theta}))} &
e^{n(-V(\theta)-i\theta+g_+(e^{i\theta})+g_-(e^{i\theta}))}\\\\ 0 &
e^{n(g_+(e^{i\theta})-g_-(e^{i\theta}))}
\end{array}\right)\,,\hspace{0.2 in}
\text{for $\theta\in S^1$}\,,
\label{eq:Yjump}
\end{equation}
where the subscript ``$+$'' indicates the boundary value taken from within
the circle, and ``$-$'' indicates the boundary value taken from outside.
Without any loss of generality, we may represent $g(z)$ in the form of
a complex logarithmic potential
\begin{equation}
g(z)=\int_{-\pi}^\pi L_{\theta'}(z)\psi(\theta')\,d\theta'
\label{eq:logtransform}
\end{equation}
where for each $\theta\in S^1$ we consider the function
$L_\theta(z):=\log(z-e^{i\theta})$ to be real for $z-e^{i\theta}$
sufficiently large and positive real, with the branch cut from the
point $z=e^{i\theta}$ in the clockwise direction along the unit circle
to the negative real axis, and then along the negative real axis to
$z=-\infty$. We also require that
\begin{equation}
\int_{-\pi}^\pi\psi(\theta')\,d\theta' = 1
\label{eq:probability}
\end{equation}
in order to satisfy the required normalization condition $g(z)\sim\log(z)$ as
$z\rightarrow\infty$.

Additional conditions may now be placed on $g$, or equivalently on
$\psi$, in order to make the jump condition \eqref{eq:Yjump}
for ${\bf Y}^n(z)$ asymptotically tractable.  One key condition is
that the function $\psi$ should be real-valued.  From this
condition it follows that $g_+(z)+g_-(z)-\log(z)$ has a constant
imaginary part for $|z|=1$, as the following argument shows.  When
$|z|=e^{i\theta}$, the identity $d/d\theta = izd/dz$ yields
\begin{equation}
\frac{d}{d\theta}\Im\left(g_+(z)+g_-(z)-\log(z)\right) = \Re\left(zg'_+(z)+zg'_-(z)-1\right)\,.
\end{equation}
Using \eqref{eq:probability} and differentiating \eqref{eq:logtransform}
under the integral with respect to $z$ for $|z|\neq 1$, we obtain
\begin{equation}
\frac{d}{d\theta}\Im\left(g_+(z)+g_-(z)-\log(z)\right)=
-\lim_{\epsilon\downarrow 0}\Re\left[\int_{-\pi}^\pi
\frac{e^{2i\theta'}-(1+\epsilon^2)e^{2i\theta}}{(e^{i\theta'}-e^{i\theta})^2-
\epsilon^2e^{2i\theta}}
\psi(\theta')\,d\theta'\right]\,.
\end{equation}
Assuming reality of $\psi$ we can pass the real part under
the integral and thus arrive at
\begin{equation}
\frac{d}{d\theta}\Im
\left(g_+(z)+g_-(z)-\log(z)\right)=-\lim_{\epsilon\downarrow 0}
\int_{-\pi}^\pi\frac{\epsilon^4+2\epsilon^2+2\epsilon^2\cos(\zeta)-4\epsilon^2\cos^2(\zeta)}
{4+\epsilon^4 + (4\epsilon^2-8)\cos(\zeta) + (4-4\epsilon^2)\cos^2(\zeta)}
\,\psi(\theta')\,d\theta'
\label{eq:zetaintegral}
\end{equation}
where $\zeta=\theta'-\theta$.  The numerator of the fraction in the
integrand is clearly uniformly bounded by a quantity of order
$\epsilon^2$, and the minimum value of the denominator is achieved
when $\cos(\zeta)=(2-\epsilon^2)/(2-2\epsilon^2)\in (-1,1)$ yielding a
minimum value that has the asymptotic expansion
$\epsilon^2-\epsilon^4+O(\epsilon^6)$ as $\epsilon\downarrow 0$.
Consequently the fraction is uniformly bounded independent of
$\epsilon$.  Moreover, the fraction is easily seen to tend to zero
pointwise in the limit $\epsilon\downarrow 0$ for $\zeta\neq 0$.  A
dominated convergence argument therefore shows that the limit on the
right-hand side of \eqref{eq:zetaintegral} is zero.  See also
\cite{longest}, where it is also shown that the constant value is
exactly $\Im(g_+(z)+g_-(z)-\log(z))=\pi$.  Another consequence of assuming
that $\psi$ is real comes from noting that
\begin{equation}
g_+(e^{i\theta})-g_-(e^{i\theta})=2\pi i\int_{\theta}^{\pi}\psi(\theta')\,
d\theta'\,.
\end{equation}
It therefore follows that $g_+(z)-g_-(z)$ is purely imaginary for
$|z|=1$.

Since
\begin{equation}
\Re(g_+(e^{i\theta})+g_-(e^{i\theta}))=2\int_{-\pi}^\pi
\log\left|e^{i\theta}-e^{i\theta'}\right|\psi(\theta')\,d\theta'\,,
\end{equation}
one is led to seek $\psi$ so that the circle is split into intervals
of two different types:
\begin{itemize}
\item[]{\bf Bands:} For $\theta$ in a band $I$, $\psi(\theta)$ is real
and positive, and
\begin{equation}
2\int_{-\pi}^\pi\log\left|e^{i\theta}-e^{i\theta'}\right|
\psi(\theta')\,d\theta' -V(\theta)=\ell\,,\hspace{0.2 in}
\text{for $\theta$ in a band}\,,
\label{eq:bandcondition}
\end{equation}
where $\ell$ is a real constant (the same constant for all bands ---
in fact if there is only one band it turns out that $\ell=-V_0$).
Thus, in a band the jump condition \eqref{eq:Yjump} for ${\bf Y}^n(z)$
takes the form
\begin{equation}
{\bf Y}^n_+(z)={\bf Y}^n_-(z)\left(\begin{array}{cc}
e^{in(\kappa_I(\theta)-\pi)} & e^{n(\ell+i\pi)}\\\\
0 & e^{-in(\kappa_I(\theta)-\pi)}\end{array}\right)\,,
\hspace{0.2 in}\text{for $\theta$ in a band $I$}\,,
\end{equation}
where $\kappa_I(\theta)$ is a strictly increasing real function.
\item[]{\bf Gaps:}  For $\theta$ in a gap $\Gamma$, $\psi(\theta)\equiv 0$ 
and we have the
strict inequality
\begin{equation}
2\int_{-\pi}^\pi\log\left|e^{i\theta}-e^{i\theta'}\right|\psi(\theta')
\,d\theta'-V(\theta)<\ell\,,\hspace{0.2 in}
\text{for $\theta$ in a gap}\,.
\label{eq:gapcondition}
\end{equation}
Thus, in a gap, the jump condition \eqref{eq:Yjump} for ${\bf Y}^n(z)$
takes the form
\begin{equation}
{\bf Y}^n_+(z)={\bf Y}^n_-(z)\left(\begin{array}{cc}
e^{in(\kappa_\Gamma-\pi)} & e^{n(\ell+i\pi)}\cdot\mbox{exponentially small}\\\\
0 & e^{-in(\kappa_\Gamma-\pi)}
\end{array}\right)\,,\hspace{0.2 in}\text{for $\theta$ in a gap $\Gamma$}\,,
\end{equation}
where $\kappa_\Gamma$ is a real constant.
\end{itemize}
The alternative conditions \eqref{eq:bandcondition} and 
\eqref{eq:gapcondition} are exactly the Euler-Lagrange conditions for the minimization of the weighted logarithmic energy
\begin{equation}
E[\psi]:=\int_{-\pi}^\pi\int_{-\pi}^\pi\log\left(\frac{1}{\left|e^{i\theta}-e^{i\theta'}\right|}\right)\psi(\theta')\psi(\theta)\,d\theta'\,d\theta
+ \int_{-\pi}^\pi V(\theta')\psi(\theta')\,d\theta'
\end{equation}
over all probability measures $\psi(\theta)\,d\theta>0$ supported on the
unit circle. The constant $\ell$ is the Lagrange multiplier
introduced to enforce the constraint \eqref{eq:probability}.  The
minimizing measure with density $\psi(\theta)$ is called the {\em equilibrium
measure}.  The connection of this extremal problem with the jump condition
\eqref{eq:Yjump} through the Euler-Lagrange variational conditions suggests
that logarithmic potential theory plays an important role in the
asymptotic theory of orthogonal polynomials for general exponentially
varying weights of the form \eqref{eq:varyingweight}.  The use of
equilibrium measures for the asymptotic analysis of the
Riemann-Hilbert problem associated to orthogonal polynomials on
${\mathbb R}$ with analytic exponentially varying weights was carried
out in
\cite{op1,op2}.  For orthogonal polynomials on the unit
circle with analytic exponentially varying weights of a specific form,
this was done in \cite{longest}.

The transformation
\begin{equation}
{\bf Z}^n(z)=\left\{\begin{array}{ll}
-e^{-n\ell\sigma_3/2}{\bf Y}^n(z)e^{n\ell\sigma_3/2}\,, &|z|<1\\\\
e^{-n\ell\sigma_3/2}{\bf Y}^n(z)e^{n\ell\sigma_3/2}\,, & |z|>1
\end{array}\right.
\end{equation}
leads to jump conditions of the form
\begin{equation}
{\bf Z}^n_+(z)={\bf Z}^n_-(z)\left(\begin{array}{cc}
e^{in\kappa_I(\theta)} & 1 \\\\ 0 & e^{-in\kappa_I(\theta)}
\end{array}\right)\,,\hspace{0.2 in}
\text{for $z=e^{i\theta}$ in a band $I$}\,,
\end{equation}
which is of exactly the same form as \eqref{eq:Sjump}, and
\begin{equation}
{\bf Z}^n_+(z)={\bf Z}^n_-(z)\left(\begin{array}{cc}
e^{in\kappa_\Gamma} & \mbox{exponentially small}\\\\
0 & e^{-in\kappa_\Gamma}\end{array}\right)\,,\hspace{0.2 in}
\text{for $z=e^{i\theta}$ in a gap $\Gamma$}\,.
\end{equation}
Clearly, we also have the normalization condition ${\bf
Z}^n(z)\rightarrow{\mathbb I}$ as $z\rightarrow\infty$. 

Now, if the function $V(\theta)$ is such that the whole unit circle
consists of a single band, or equivalently there are no gaps in the
support of the equilibrium measure, then in fact ${\bf Z}^n(z)={\bf
S}^n(z)$ as defined in \S~\ref{sec:varyingconvert} and the analysis
proceeds as in the main body of this paper.  However, in the more
general context --- when it is only true that the support of the
equilibrium measure consists of a finite number of disjoint intervals
on the unit circle --- a modification of the $\dbar$ steepest descent
method described in \S~\ref{sec:varying} is required.  A more general
method may still be based on the algebraic factorization
\eqref{eq:factorizationvarying} of the jump matrix in each band;
however the annuli $A_\pm$ must be replaced with
a system of lens-shaped regions $A_\pm^I$ adjacent to each band $I$.
Since the variational inequalities become less effective near the band edges,
a local analysis must be supplied to control the error.
When such a method is developed, the conditions
\eqref{eq:tcondition} or \eqref{eq:convexity} on the function
$V(\theta)$ can be dropped as long as it is known that there are only
a finite number of gaps in the support of the equilibrium measure.  In
\cite{otherpaper} it is shown how to carry out such a program in the
context of nonanalytic exponentially varying weights on the real line,
where a convexity condition is known to guarantee the existence of a single
isolated band $[\alpha,\beta]\subset\mathbb{R}$.

Let us show how the logarithmic potential theory described briefly
above leads to the condition \eqref{eq:convexity} guaranteeing that
the support of the equilibrium measure is the entire unit circle.  If
$V$ is of class $C^2(S^1)$,
then for $\theta$ in any gap in the support of the
equilibrium measure the real function
\begin{equation}
\Phi(\theta):=2\int_{-\pi}^\pi\log\left|e^{i\theta}-e^{i\theta'}\right|
\psi(\theta')\,d\theta'-V(\theta)
\end{equation}
is also twice differentiable.
In fact, a calculation shows that
\begin{equation}
\Phi''(\theta) = -2\int_{-\pi}^\pi\frac{\psi(\theta')\,d\theta'}
{\left|e^{i\theta}-e^{i\theta'}\right|^2} - V''(\theta)\,.
\end{equation}
Note that this is not a singular integral since $e^{i\theta}$ is
assumed to lie outside the support of the equilibrium measure.  Let us
assume both the condition \eqref{eq:convexity}, and also the existence
of a gap in the support of the equilibrium measure; we will then derive a
contradiction.  Now, $\Phi(\theta)$ is continuous for
$-\pi<\theta\le\pi$.  Therefore, since $\Phi(\theta)=\ell$ at both
endpoints of the gap (according to (\ref{eq:bandcondition})), and
$\Phi(\theta)<\ell$ strictly in the interior of the gap (according to
(\ref{eq:gapcondition})), there must be a point $\theta$ in the gap at
which $\Phi''(\theta)>0$.  However, since $\left|e^{i\theta}-e^{i\theta'}\right|\le 2$ for
all angles $\theta'$, we see that
\begin{equation}
\Phi''(\theta)\le -\frac{1}{2} -V''(\theta)
\end{equation}
which is negative in view of the assumption \eqref{eq:convexity}, thus
establishing the desired contradiction.

\section{Comparison of the Conditions (\ref{eq:tcondition}) and (\ref{eq:convexity})}
\label{app:compare}
In this appendix, we illustrate by concrete examples how the two conditions
\eqref{eq:tcondition} and \eqref{eq:convexity}, 
while both sufficient for the prevention of gaps in the support of the
equilibrium measure, are completely independent.  Thus neither
condition implies the other.

First, consider the example $V(\theta)=A\cos(k\theta)$ where
$k=1,2,3,\dots$ and $A\in{\mathbb R}$ are parameters.  A direct
calculation shows that
\begin{equation}
\begin{array}{rcl}
\mbox{Condition (\ref{eq:tcondition})} &\Leftrightarrow &
\displaystyle |A|<\frac{1}{k}\,,\\\\
\mbox{Condition (\ref{eq:convexity})} &\Leftrightarrow & \displaystyle
|A|<\frac{1}{2k^2}\,.
\end{array}
\end{equation}
Since $2k^2>k$ for all $k=1,2,3,\dots$, we see that if
\begin{equation}
\frac{1}{2k^2}<|A|<\frac{1}{k}
\end{equation}
then condition (\ref{eq:tcondition}) is satisfied but
(\ref{eq:convexity}) is not.

Next, consider for $M>0$ and $\epsilon>0$ the example
\begin{equation}
V(\theta)=\left\{\begin{array}{ll}
\displaystyle
-\frac{M}{6}|\theta|^3 + \left(\frac{M\epsilon}{2}-\frac{M\epsilon^2}{4\pi}
\right)\theta^2 +\frac{M\epsilon^4}{4\pi}-\frac{M\epsilon^3}{3}\,, &
|\theta|\le\epsilon
\\\\
\displaystyle
-\frac{M\epsilon^2}{4\pi}\theta^2 + \frac{M\epsilon^2}{2}|\theta| +
\frac{M\epsilon^4}{4\pi}-\frac{M\epsilon^3}{2}\,, &\epsilon<|\theta|<\pi\,.
\end{array}\right.
\end{equation}
On one hand, we have
\begin{equation}
V''(\theta)=\left\{\begin{array}{ll}
\displaystyle M(\epsilon-|\theta|) -\frac{M\epsilon^2}{2\pi}\,, &|\theta|\le
\epsilon\\\\ \displaystyle
-\frac{M\epsilon^2}{2\pi}\,, & \epsilon<|\theta|<\pi\,,
\end{array}\right.
\end{equation}
and so
\begin{equation}
\mbox{Condition (\ref{eq:convexity})}\,\,\,\Leftrightarrow\,\,\,
M<\frac{\pi}{\epsilon^2}\,.
\end{equation}
On the other hand, direct calculation of the Fourier coefficients gives
\begin{equation}
V_k=\left\{\begin{array}{ll}
\displaystyle-\frac{2M}{\pi k^4}\sin^2\left(\frac{k\epsilon}{2}\right)\,, &
k\neq 0\,,\\\\
\displaystyle \frac{7M\epsilon^4}{48\pi}-\frac{M\epsilon^3}{4}+
\frac{M\pi\epsilon^2}{12}\,, & k=0\,.
\end{array}\right.
\end{equation}
Thus, we may estimate the sum in the condition (\ref{eq:tcondition}) as follows:
\begin{equation}
\sum_{k=-\infty}^\infty |kV_k|=\sum_{k=1}^\infty\frac{4M}{\pi k^3}\sin^2
\left(\frac{k\epsilon}{2}\right)\ge\sum_{k=1}^{\lceil \pi/\epsilon\rceil}
\frac{4M}{\pi k^3}\sin^2\left(\frac{k\epsilon}{2}\right)\ge
\sum_{k=1}^{\lceil \pi/\epsilon\rceil}
\frac{M\epsilon^2}{4k}
\end{equation}
since $|\sin(x)|\ge|x|/2$ for $|x|\le\pi/2$.  If the product
$M\epsilon^2$ is held fixed this lower bound can be made
arbitrarily large (and in particular larger than one) simply by
taking $\epsilon$ sufficiently small due to the divergence of the
harmonic series.  We therefore see that if $0<c\le
M\epsilon^2<\pi$ and if $\epsilon$ is sufficiently small, then
condition (\ref{eq:convexity}) is satisfied but condition
(\ref{eq:tcondition}) is not.

\end{document}